\documentclass[12pt]{amsart}
\usepackage{amssymb}
\usepackage{amsbsy}
\usepackage{amscd}
\usepackage{verbatim}
\makeatletter
\def\Bbb{\mathbb}
\def\frak{\mathfrak}

\newenvironment{pf*}[1]{\proof[#1]}{\endproof}
\newcommand{\rom}{\textup}
\makeatother
\newenvironment{aenume}{%
  \begin{enumerate}%
  }{\end{enumerate}}
\newenvironment{enume}{%
  \begin{enumerate}%
  }{\end{enumerate}}
\newtheorem{Theorem}[equation]{Theorem}
\newtheorem{Corollary}[equation]{Corollary}
\newtheorem{Lemma}[equation]{Lemma}
\newtheorem{Proposition}[equation]{Proposition}

\theoremstyle{definition}
\newtheorem{Definition}[equation]{Definition}
\newtheorem{Example}[equation]{Example}

\newtheorem{Notation}[equation]{Notation}
\newtheorem{Convention}[equation]{Convention}
\newtheorem{Conjecture}[equation]{Conjecture}

\theoremstyle{remark}
\newtheorem{Remark}[equation]{Remark}
\newtheorem*{Claim}{Claim}
\numberwithin{equation}{section}
\newcommand{\thmref}[1]{Theorem~\ref{#1}}
\newcommand{\secref}[1]{\S\ref{#1}}
\newcommand{\lemref}[1]{Lemma~\ref{#1}}
\newcommand{\propref}[1]{Proposition~\ref{#1}}
\newcommand{\corref}[1]{Corollary~\ref{#1}}
\newcommand{\subsecref}[1]{\S\ref{#1}}

\newcommand{\defref}[1]{Definition~\ref{#1}}
\newcommand{\remref}[1]{Remark~\ref{#1}}
\newcommand{\defeq}{\overset{\operatorname{\scriptstyle def.}}{=}}
\newcommand{\C}{{\Bbb C}}
\newcommand{\Z}{{\Bbb Z}}
\newcommand{\Q}{{\Bbb Q}}
\newcommand{\R}{{\Bbb R}}
\newcommand{\GL}{\operatorname{GL}}

\newcommand{\Hom}{\operatorname{Hom}}
\newcommand{\Ext}{\operatorname{Ext}}
\newcommand{\Ker}{\operatorname{Ker}}

\newcommand{\Ima}{\operatorname{Im}}

\newcommand{\rank}{\operatorname{rank}}

\newcommand{\tr}{\operatorname{tr}}

\newcommand{\id}{\operatorname{id}}

\newcommand{\vin}{\operatorname{in}} % incoming vertex
\newcommand{\vout}{\operatorname{out}} % outgoing vertex
\newcommand{\bM}{{\mathbf M}} % the vector space
\newcommand{\M}{{\frak M}} % the moduli space

\newcommand{\La}{{\frak L}} % the Lagrangian variety
\newcommand{\dslash}{/\!\!/} % for algebro-geometric quotient (double slash)
\newcommand{\bv}{{\mathbf v}} % the vector v
\newcommand{\bw}{{\mathbf w}} % the vector w
\newcommand{\bC}{{\mathbf C}} % the Cartan matrix
\newcommand{\codim}{\operatorname{codim}} % codimension
\newcommand{\topdeg}{\operatorname{top}} % top degree
\newcommand{\Uq}{{\mathbf U}_q(\mathfrak g)} % the QUE algebra

\newcommand{\Ua}{{\mathbf U}_q(\widehat{\mathfrak g})} % the quantum affine
\newcommand{\Uai}{{\mathbf U}_q^\Z(\widehat{\mathfrak g})} % the
\newcommand{\Ul}{{\mathbf U}_q({\mathbf L}{\mathfrak g})} % the quantum
\newcommand{\Ui}{{\mathbf U}^{\Z}_q(\mathfrak g)}

\newcommand{\Uli}{{\mathbf U}^{\Z}_q({\mathbf L}\mathfrak g)}
\newcommand{\Ule}{{\mathbf U}_{\varepsilon}({\mathbf L}\mathfrak g)}
\newcommand{\Pa}{{\frak P}} % Hecke correspondence
\newcommand{\HomE}{\operatorname{E}}
\newcommand{\HomL}{\operatorname{L}}

\newcommand{\shfO}{\mathcal O}

\newcommand{\Wedge}{{\textstyle \bigwedge}}

\newcommand{\Mw}{\M(\bw)}
\newcommand{\Law}{\La(\bw)}
\newcommand{\Zl}{\widetilde{\mathfrak Z}}
\newcommand{\Zm}{\mathfrak Z}
\newcommand{\Irr}{\operatorname{Irr}}
\newcommand{\B}{\mathcal B}
\newcommand{\wt}{\operatorname{wt}}
\newcommand{\wH}{\widetilde H}

\begin{document}
\author{Hiraku Nakajima}
\title[quiver varieties and tensor products]
{Quiver varieties and tensor products
%\\ {\rm \small (Preliminary version: \today)}
}
\address{Department of Mathematics, Kyoto University, Kyoto 606-8502,
Japan}
\email{nakajima@kusm.kyoto-u.ac.jp}
\thanks{Supported by the Grant-in-aid
for Scientific Research (No.11740011), the Ministry of Education,
Japan.}
%
%\subjclass{Primary 17B37;
%Secondary 14D21, 14L30, 16G20, 33D80}
%
\begin{abstract}
In this article, we give geometric constructions of tensor products
in various categories using quiver varieties.
More precisely, we introduce a lagrangian subvariety $\Zl$ in a quiver
variety, and show the following results:
(1) The homology group of $\Zl$ is a representation of a symmetric
Kac-Moody Lie algebra $\mathfrak g$, isomorphic to the tensor product
$V(\lambda_1)\otimes\dots\otimes V(\lambda_N)$ of integrable highest
weight modules.
(2) The set of irreducible components of $\Zl$ has a structure of a
crystal, isomorphic to that of the $q$-analogue of
$V(\lambda_1)\otimes\dots\otimes V(\lambda_N)$.
(3) The equivariant $K$-homology group of $\Zl$ is isomorphic to the
tensor product of universal standard modules of the quantum loop 
algebra $\Ul$, when $\mathfrak g$ is of type $ADE$.
We also give a purely combinatorial description of the crystal of (2).
This result is new even when $N=1$.
\end{abstract}

\maketitle
\tableofcontents
%%%%%%%%%%%%%%%%%%%%%%%%%%%%%%%%%%%%%%%%%%%%%%%%%%%%%%%%%%%%%%%%%%%%%%%%

\section*{Introduction}\label{sec:intro}

Let $\mathfrak g$ be a symmetric Kac-Moody Lie algebra and
let $\Uq$ be the quantum enveloping algebra of Drinfeld-Jimbo attached 
to $\mathfrak g$.
For each dominant weight $\bw$ of $\mathfrak g$, the author associated
a nonsingular variety $\M(\bw)$ (called a {\it quiver variety\/}),
containing a half dimensional subvariety $\La(\bw)$
\cite{Na-quiver,Na-alg}. It is related to the representation theory of
$\mathfrak g$ and $\Uq$ as follows:
\begin{enumerate}
\item The top degree homology group
$H_{\topdeg}(\La(\bw),\C)$ has a structure of a $\mathfrak g$-module,
simple with highest weight $\bw$. (\cite{Na-alg})

\item The set $\Irr\La(\bw)$ of irreducible components of $\La(\bw)$
has a structure of a crystal, isomorphic to the crystal of the simple
$\Uq$-module with highest weight $\bw$. (Kashiwara-Saito~\cite{KS,Saito})
\end{enumerate}

Suppose $\mathfrak g$ is of type $ADE$. Let ${\mathbf L}{\mathfrak g}
= \mathfrak g[x,x^{-1}]$ be the loop algebra of $\mathfrak g$
and let $\Ul$ be the quantum loop algebra (the quantum affine algebra
without centeral extension and degree operator) with the integral form 
$\Uli$ generated by $q$-divided powers. In \cite{Na-qaff}, the author
showed the following:
\begin{enumerate}
\setcounter{enumi}{2}
\item The equivariant $K$-homology group
$K^{H_{\bw}\times\C^*}(\La(\bw))$ has a structure of
\(
    \Uli[x_1^\pm,\dots, x_L^\pm]
\)-module, which is simple and has various nice properties.
\end{enumerate}
Here $L = \sum_k \langle h_k, \bw\rangle$ and $H_{\bw}$ is a torus of
dimension $L$ acting on $\M(\bw)$ and $\La(\bw)$.
We call it a {\it universal standard module}, and denote it by
$M(\bw)$.
(In the main body of this article, we replace $H_{\bw}$ by a product
of general linear groups, and
\(
    \Uli[x_1^\pm,\dots, x_N^\pm]
\)
by the invariant part of a product of symmetric groups.)
It was shown that any simple $\Ul$-module is obtained as a
quotient of a specialization of $M(\bw)$. This specialization is called 
a {\it standard module}.
Moreover, the multiplicities of simple modules in standard modules are
computable by a combinatorial algorithm \cite{Na-qchar}.

In this article, we generalize results~(1),(2),(3) to the case of
tensor products of simple modules. For given dominant weights $\bw^1$,
$\bw^2$, \dots, $\bw^N$ we introduce another half-dimensional
subvarity $\Zl$ of $\M(\bw)$, containing $\La(\bw)$
($\bw=\bw^1+\dots+\bw^N$), and show the followings:
\begin{enume}
\item $H_{\topdeg}(\Zl,\C)$ has a structure of a $\mathfrak g$-module,
isomorphic to the tensor product of simple $\mathfrak g$-modules with 
highest weights $\bw^1$, \dots, $\bw^N$. (\thmref{thm:main2})

\item The set $\Irr\Zl$ of irreducible components of $\Zl$ has a
structure of a crystal, isomorphic to the tensor product of
the crystals $\Irr\La(\bw^1)$, \dots, $\Irr\La(\bw^N)$. (\thmref{thm:main1})

\item When $\mathfrak g$ is of type $ADE$, the equivariant
$K$-homology group $K^{H_\bw\times\C^*}(\Zl)$ is isomorphic to the
tensor product of universal standard modules $M({\bw^1})$, \dots,
$M({\bw^N})$. (\thmref{thm:main4})
\end{enume}

The result (2)' means that $\{ \Irr\La(\bw)\mid \bw\in P^+
\}$ is a closed family of highest weight normal crystals (see
\subsecref{subsec:crystal} for definition). This property
characterizes crystals of simple highest weight modules. Thus we
obtain a new proof of (2).

As an application of (3)', we give a new proof of the main result of
Varagnolo-Vasserot \cite{VV} (see \corref{cor:VV}). It says that a
standard module is a tensor product of fundamental representations
in an appropriately chosen order.

We also give a combinatorial description of the crystal $\Irr\La(\bw)$
in \secref{sec:comb}. It is essentially the same as a combinatorial
description of $V(\bw)$, provided by the embedding theorem of the
crystal $\B(\infty)$ of the lower part of the quantized enveloping
algebra (\cite{Kas-Dem}). We shall discuss this description further,
relating it with the theory of $q$--characters in a separate
publication.

The varieties $\Zl$, $\Zm$ are defined as attracting sets of a
$\C^*$-action for some one parameter subgroup $\lambda\colon \C^*\to
H_\bw$. The relation between the $\C^*$-action and tensor products has
been known to many people after the author related the quiver
varieties to Lusztig's construction of canonical bases
\cite{Na-quiver}, since Lusztig defined the comultiplcation in a
geometric way (see \cite{Lu-book}).
Moreover, when $\mathfrak g$ is of type $A_n$ and $\bw$ is a multiple
of the fundamental weight corresponding to the vector representation,
the quiver variety $\M(\bw)$ is isomorphic to the cotangent bundle of
the $n$-step partial flag variety. In this case the comultiplication
was constructed by Ginzburg-Reshetkhin-Vasserot~\cite{GRV}. The result
was also mentioned without detail in an earlier paper by
Grojnowski~\cite{Gr_copr}.
And the author used the $\C^*$-action to compute Betti numbers of
$\M(\bw)$ (when $\mathfrak g$ is of type $A$), and checked that the
generating function of the Euler numbers is the character of the
tensor product \cite{Na-hom}.
For general $\mathfrak g$, Grojnowski mentioned, in his
`advertisement' \cite{Gr} of his book, that the coproduct is defined
by the localization to the fixed point set of the
$\C^*$-action. However the details of the construction of $\Ul$-module
structures were not given.
The details were given in \cite{Na-qaff}, and by the localization,
tensor products were studied, but only for generic parameters [14.1.2,
loc.\ cit.] (see also \lemref{lem:general_point}).
Tensor products for arbitrary parameters need further study, and it
was first done by Varagnolo-Vasserot~\cite{VV}. The result~(3)' above
is motivated by their study, although the author knew $\Zl$ before
their paper appeared.

The variety $\Zl$ is an analogue of subvarieties of cotangent bundles
of flag varieties, introduced by Lusztig \cite{Lu-Base2}. (Our
notation is taken from his.)
In his picture, $\M(\bw)$ corresponds to Slodowy's variety, $\La(\bw)$
to the Springer fiber $\mathcal B_x$. The equivariant $K$-homology
group of $\mathcal B_x$ is a module of the affine Hecke algebra
$\widehat H_q$. The equivariant $K$-homology of $\Zl$ is an induced
module of a module of above type for a smaller affine Hecke algebra. 

One of motivations of \cite{Lu-Base2} was a conjectural construction
of a base in the equivariant $K$-homology of $\mathcal B_x$. Lusztig
pointed out a possibility of a similar construction for quiver
varieties \cite{Lu-rem}. Combining a result in \cite{Na-qaff} with a
recent result by Kashiwara \cite{Kas2}, the universal standard module
$M(\bw)$ has a global crystal base. It is interesting to compare his
base with Lusztig's (conjectural) base.

Finally we comment that there is a geometric construction of tensor
products of {\it two\/} simple $\mathfrak g$-modules by Lusztig
\cite{Lus-tensor}. His variety is a subvariety of a product of quiver
varieties. In fact, an open subvariety of an analogue of Steinberg
variety, which will be denoted by $Z(\bw)$ in this article.
The relation between his variety and $\Zl$ is not clear, although
there is an example where a close relation can be found (see
\secref{sec:example}).
And it seems difficult to generalize his varities to the case of
tensor products of several modules. (Compare \secref{sec:general}).

After this work was done, we were informed that Malkin also defined
the variety $\Zl$ and obtained the result (2)' above \cite{Malkin}.

\subsection*{Acknowledgement}
This work originated in G.~Lusztig's question about analogue of his
variety in quiver varieties, asked at the Institute for Advanced
Study, 1998 winter. It is a great pleasure for me to answer his
question after two years.
I would like to thank I.~Grojnowski for explaining me his construction
of tensor products for generic parameters during 2000 winter.
I also express my sincere gratitude to M.~Kashiwara for interesting
discussions about crystal. 

%\tableofcontents

\section{Preliminaries (I) -- algebraic part}

\subsection{Quantized enveloping algebra}\label{subsec:QUE}
We briefly recall the definition of the quantized enveloping algebra
in this subsection. See \cite{Kas} for further detail.

A {\it root datum\/} consists of
\begin{enumerate}
  \item $P$ : free $\Z$-module (weight lattice),
  \item $P^* = \Hom_{\Z}(P, \Z)$ with a natural pairing 
    $\langle\ , \ \rangle\colon P\otimes P^*\to \Z$,
    \item a finite set $I$ (index set of simple roots)
  \item $\alpha_k\in P$ ($k \in I$) (simple root),
  \item $h_k \in P^*$ ($k \in I$) (simple coroot),
  \item a symmetric bilinear form $(\ ,\ )$ on $P$.
\end{enumerate}
Those are required to satisfy the followings:
\begin{aenume}
\item $\langle h_k, \lambda\rangle 
= 2(\alpha_k, \lambda)/(\alpha_k,\alpha_k)$ for  
  $k \in I$ and $\lambda\in P$,
\item $\bC \defeq (\langle h_k, \alpha_l \rangle)_{k,l}$ is a symmetrizable
  generalized Cartan matrix, i.e., $\langle h_k, \alpha_k\rangle = 2$,
  and 
  $\langle h_k, \alpha_l\rangle\in\Z_{\le 0}$ and
  $\langle h_k, \alpha_l\rangle = 0 \Longleftrightarrow
   \langle h_l, \alpha_k\rangle = 0$ for $k\ne l$,
\item $(\alpha_k,\alpha_k)\in 2\Z_{> 0}$,
\item $\{\alpha_k\}_{k\in I}$ are linearly independent,
\item there exists $\Lambda_k\in P$ ($k \in I$) such that 
  $\langle h_l, \Lambda_k\rangle = \delta_{kl}$
  (fundamental weight).
\end{aenume}
Let $\mathfrak g$ be the symmetrizable Kac-Moody Lie algebra
corresponding to the generalized Cartan matrix $\bC$ with the Cartan
subalgebra $\mathfrak h = P^*\otimes_\Z \Q$. 
Let $Q = \bigoplus_k \Z \alpha_k\subset P$ be the root lattice.
Let $P^+$ be the semigroup of dominant weights, i.e.,
$P^+ = \{ \lambda\mid \langle h_k,\lambda\rangle\ge 0 \}$.
Let $Q^+ = \sum_k \Z_{\ge 0} \alpha_k$.

Let $q$ be an indeterminate. For nonnegative integers $n\ge r$, define 
\begin{equation*}%\label{eq:q-binom}
  [n]_q \defeq \frac{q^n - q^{-n}}{q - q^{-1}}, \quad
  [n]_q ! \defeq 
  \begin{cases}
   [n]_q [n-1]_q \cdots [2]_q [1]_q &(n > 0),\\
   1 &(n=0),
  \end{cases}
  \quad
  \begin{bmatrix}
  n \\ r
  \end{bmatrix}_q \defeq \frac{[n]_q !}{[r]_q! [n-r]_q!}.
\end{equation*}

The quantized universal enveloping algebra $\Uq$ of the Kac-Moody
algebra $\mathfrak g$ is the $\Q(q)$-algebra generated by $e_k$, $f_k$
($k \in I$), $q^h$ ($h\in P^*$) with relations
{\allowdisplaybreaks
\begin{subequations}
\begin{gather}
  q^0 = 1, \quad q^h q^{h'} = q^{h+h'},\\
  q^h e_k q^{-h} = q^{\langle h, \alpha_k\rangle} e_k,\quad
  q^h f_k q^{-h} = q^{-\langle h, \alpha_k\rangle} f_k,\\
  e_k f_l - f_l e_k = \delta_{kl}
   \frac{t_k-t_k^{-1}}{q_k - q_k^{-1}},\\
  \sum_{p=0}^{b}(-1)^p e_k^{(p)} e_l e_k^{(b-p)} =
  \sum_{p=0}^{b}(-1)^p f_k^{(p)} f_l f_k^{(b-p)} = 0 \quad
  \text{for $k\ne l$,}
\end{gather}
\end{subequations}
where $q_k = q^{(\alpha_k,\alpha_k)/2}$,
$t_k = q^{(\alpha_k,\alpha_k)h_k/2}$,
$b = 1 - \langle h_k,\alpha_l\rangle$,
$e_k^{(p)} = e_k^p/[p]_{q_k}!$, $f_k^{(p)} = f_k^p/[p]_{q_k}!$.}

Let $\Uq^+$ (resp.\ $\Uq^-$) be the $\Q(q)$-subalgebra of $\Uq$
generated by elements $e_k$'s (resp.\ $f_k$'s).  Let $\Uq^0$ be the
$\Q(q)$-subalgebra generated by elements $q^h$ ($h\in P^*$). We have
the triangular decomposition $\Uq \cong \Uq^+\otimes \Uq^0 \otimes
\Uq^-$.

Let $\Ui$ be the $\Z[q,q^{-1}]$-subalgebra of $\Uq$ generated by
elements $e_k^{(n)}$, $f_k^{(n)}$, $q^h$ for $k\in I$, $n\in\Z_{> 0}$,
$h\in P^*$.

In this article, we take the comultiplication $\Delta$ on $\Uq$ given
by
\begin{equation}\label{eq:comul}
\begin{gathered}
   \Delta q^h = q^h \otimes q^h, \quad
   \Delta e_k = e_k\otimes q^{-h_k} + 1 \otimes e_k,
\\
   \Delta f_k = f_k\otimes 1 + q^{h_k} \otimes f_k.
\end{gathered}
\end{equation}
Note that this is different from one in \cite{Lu-book}, although there
is a simple relation between them \cite[1.4]{Kas}.
The results in \cite{Na-qaff} hold for either comultiplication
(tensor products appear in (1.2.19) and (14.1.2)).
In \cite[\S2]{Na-qchar} another comultiplication was used. If we reverse
the order of the tensor product, the results hold.

For each dominant weight $\lambda\in P^+$, there is unique simple
module $V(\lambda)$ with highest weight $\lambda$. The highest weight
vector is denoted by $b_\lambda$.

Later we use the classical counter part of $\Uq$. We just erase $q$ in 
the above notation, e.g., $\Uq\Rightarrow\mathbf U(\mathfrak g)$, etc.
Simple highest weight modules are denoted by the same notation
$V(\lambda)$.

\subsection{Crystal}\label{subsec:crystal}
Let us review the notion of crystals briefly. See \cite{Kas,KS} for
detail.

\begin{Definition}\label{def:crystal}
A {\it crystal\/} $\B$ associated with a root datum in
\subsecref{subsec:QUE} is a set together with maps
\(
  \wt\colon \B\to P
\),
\(
  \varepsilon_k,\varphi_k\colon \B\to\Z\sqcup\{-\infty\}
\),
\(
  \widetilde e_k,\widetilde f_k\colon \B\to \B\sqcup\{ 0\}
\)
($k\in I$) satisfying the following properties
{\allowdisplaybreaks
\begin{subequations}
\begin{align}
   & \varphi_k(b) = \varepsilon_k(b) + \langle h_k, \wt(b)\rangle,
\\
   & \wt(\widetilde e_k b) = \wt(b) + \alpha_k,\ 
   \varepsilon_k(\widetilde e_k b) = \varepsilon_k(b)-1,\ 
   \varphi_k(\widetilde e_k b)=\varphi_k(b)+1,
   \quad\text{if $\widetilde e_k b\in\B$},
\\
   & \wt(\widetilde f_k b) = \wt(b) - \alpha_k,\ 
   \varepsilon_k(\widetilde f_k b) = \varepsilon_k(b)+1,\ 
   \varphi_k(\widetilde f_k b)=\varphi_k(b)-1,
   \quad\text{if $\widetilde f_k b\in\B$},
\\
   & b' = \widetilde f_k b \Longleftrightarrow
     b = \widetilde e_k b' \quad\text{for $b$, $b'\in\B$},
\\
   & \text{if $\varphi_k(b)=-\infty$ for $b\in\B$, then $\widetilde e_k b
   = \widetilde f_k b = 0$}
\end{align}
\end{subequations}
}
\end{Definition}
We set $\wt_k(b) = \langle h_k,\wt(b)\rangle$.

The crystal was introduced by abstracting the notion of crystal bases
constructed by Kashiwara~\cite{Kas}. Thus we have the following
examples of crystals.

\begin{Notation}\label{not:Blambda}
(1) Let $\B(\infty)$ denote the crystal associated with
${\mathbf U}_q(\mathfrak g)^-$.

(2) For $\lambda\in P^+$, let $\B(\lambda)$ denote the crystal
associated with the simple $\Uq$-module $V(\lambda)$ with highest
weight $\lambda$.
\end{Notation}

We also have the following examples.
\begin{Example}\label{exp:crystal}
{\allowdisplaybreaks
(1) For all $k\in I$, we define the crystal $\B_k$ as follows:
\begin{gather*}
  \B_k = \{ b_k(n) \mid n\in \Z \},
\\
   \wt(b_k(n)) = n\alpha_k, \quad \varphi_k(b_k(n)) = n, \quad
   \varepsilon_k(b_k(n)) = - n,
\\
   \varphi_l(b_k(n)) = \varepsilon_l(b_k(n)) = -\infty \quad (l\neq k),
\\
   \widetilde e_k (b_k(n)) = b_k(n+1),\quad
   \widetilde f_k (b_k(n)) = b_k(n-1),
\\
   \widetilde e_l(b_k(n)) = \widetilde f_l (b_k(n)) = 0 \quad (l\neq k).
\end{gather*}
}

(2) For $\lambda\in P^+$, we define the crystal $T_\lambda$ by
\begin{gather*}
  T_\lambda = \{ t_\lambda \},
\\
  \wt(t_\lambda) = \lambda, \quad
  \varphi_k(t_\lambda) = \varepsilon_k(t_\lambda) = -\infty,
\\
  \widetilde e_k(t_\lambda) = \widetilde f_k (t_\lambda) = 0.
\end{gather*}
\end{Example}

A crystal $\B$ is called {\it normal\/} if
\begin{equation*}
   \varepsilon_k(b) = \max\{n\mid {\widetilde e_k}^n b\neq 0\}, \quad
   \varphi_k(b) = \max\{n\mid {\widetilde f_k}^n b\neq 0\}.
\end{equation*}

It is known that $\B(\lambda)$ is normal.

For given two crystals $\B_1$, $\B_2$, a {\it morphism\/} $\psi$ of
crystal from $\B_1$ to $\B_2$ is a map $\B_1\sqcup\{0\}\to \B_2\sqcup
\{0\}$ satisfying $\psi(0) = 0$ and the following conditions for all
$b\in \B_1$, $k\in I$:
{\allowdisplaybreaks
\begin{subequations}\label{eq:mor}
\begin{align}
   & \wt(\psi(b)) = \wt(b),\; \varepsilon_k(\psi(b)) = \varepsilon_k(b), \;
   \varphi_k(\psi(b)) = \varphi_k(b) \quad
   \text{if $\psi(b)\in \B_2$},
\\
   & \widetilde e_k \psi(b) = \psi(\widetilde e_k b)\quad
   \text{if $\psi(b)\in \B_2$, $\widetilde e_k b\in \B_1$},
\\
   & \widetilde f_k \psi(b) = \psi(\widetilde f_k b)\quad
   \text{if $\psi(b)\in \B_2$, $\widetilde f_k b\in \B_1$}.
\end{align}
\end{subequations}
}

A morphism $\psi$ is called {\it strict\/} if $\psi$ commutes with
$\widetilde e_k$, $\widetilde f_k$ for all $k\in I$ without any
restriction.
A morphism $\psi$ is called an {\it embedding\/} if $\psi$ is an
injective map from $B_1\sqcup\{0\}$ to $B_2\sqcup\{0\}$.

\begin{Definition}
The {\it tensor product\/} $\B_1\otimes \B_2$ of crystals $\B_1$ and
$\B_2$ is defined to be the set $\B_1\times\B_2$ with maps defined by
{\allowdisplaybreaks
\begin{subequations}
\begin{align}
   & \wt(b_1\otimes b_2) = \wt(b_1)+\wt(b_2),\label{eq:wt_tensor}
\\
   & \varepsilon_k(b_1\otimes b_2) 
   = \max(\varepsilon_k(b_1),\varepsilon_k(b_2)-\wt_k(b_1)),\label{eq:ep}
\\
   & \varphi_k(b_1\otimes b_2) 
   = \max(\varphi_k(b_2),\varphi_k(b_1)+ \wt_k(b_2)), \label{eq:phi}
\\
   & \widetilde e_k(b_1\otimes b_2)
   = 
   \begin{cases}
      \widetilde e_k b_1 \otimes b_2 &
        \text{if $\varphi_k(b_1)\ge \varepsilon_k(b_2)$},
   \\
      b_1\otimes \widetilde e_k b_2 & \text{otherwise},
   \end{cases}\label{eq:e}
\\
   & \widetilde f_k(b_1\otimes b_2)
   = 
   \begin{cases}
      \widetilde f_k b_1 \otimes b_2 &
        \text{if $\varphi_k(b_1) > \varepsilon_k(b_2)$},
   \\
      b_1\otimes \widetilde f_k b_2 & \text{otherwise}.
   \end{cases}\label{eq:f}
\end{align}
\end{subequations}
Here $(b_1,b_2)$ is denoted by $b_1\otimes b_2$ and $0\otimes b_2$,
$b_1\otimes 0$ are identified with $0$.
}\end{Definition}
It is easy to check that these satisfy the axioms in
\defref{def:crystal}. It is also easy to check that the tensor product 
of two normal crystals is again normal.

It is easy to check $(B_1\otimes B_2)\otimes B_3 = B_1\otimes
(B_2\otimes B_3)$. We denote it by $B_1\otimes B_2\otimes
B_3$. Similarly we can define $B_1\otimes\dots\otimes B_n$.
Then using
\(
    B_1\otimes\cdots\otimes B_n
    = (B_1\otimes\cdots\otimes B_{n-1})\otimes B_n
\),
we can determine the tensor product of more than two crystals
inductively as follows.
For given $b_1\otimes\cdots\otimes b_n$, we define
\begin{equation}\label{eq:vep}
%\begin{split}
%   &
   \varepsilon_k^p
   \defeq \varepsilon_k(b_p) - \sum_{q: q<p} \wt_k(b_q),
\qquad
%\\
%   &
   \varphi_k^p
   \defeq \varphi_k(b_p) + \sum_{q: q>p} \wt_k(b_q).
%\end{split}
\end{equation}
It is easy to show the following by induction:
{\allowdisplaybreaks
\begin{subequations}\label{eq:sub}
\begin{align}
%\begin{gather}
   & \wt(b_1\otimes\dots\otimes b_n) = \sum_p \wt(b_p),
\\
   & \varepsilon_k(b_1\otimes\dots\otimes b_n) 
   = \max_{1\le p\le n} \varepsilon_k^p,
\\
   &\varphi_k(b_1\otimes\dots\otimes b_n) 
   = \max_{1\le p\le n} \varphi_k^p,
\\
%   &
\begin{split}
   & \widetilde e_k(b_1\otimes\dots\otimes b_n)
   = b_1\otimes\dots\otimes \widetilde e_k b_p\otimes\dots\otimes b_n
\\
   & \qquad\qquad\qquad\text{where $p=\min\{q\mid \varepsilon_k^q =
                       \varepsilon_k(b_1\otimes\dots\otimes b_n)\}$},
\end{split}
\\
\begin{split}
   & \widetilde f_k(b_1\otimes\dots\otimes b_n)
   = b_1\otimes\dots\otimes \widetilde f_k b_p\otimes\dots\otimes b_n
\\
   & \qquad\qquad\qquad\text{where $p=\max\{q\mid \varphi_k^q =
                    \varphi_k(b_1\otimes\dots\otimes b_n)\}$}.
\end{split}
%\end{gather}
\end{align}
\end{subequations}
}

\begin{Definition}
A crystal $\B$ is said to be of {\it highest weight $\lambda$\/} if
the following conditions are satisfied:
\begin{enumerate}
\item there exists $b_\lambda\in\B$ with $\wt(b_\lambda)=\lambda$ such
that $\widetilde e_k(b_\lambda) = 0$ for all $k\in I$,
\item $\B$ is generated by $b_\lambda$, i.e., any element in $\B$ is
obtained from $b_\lambda$ by applying $\widetilde f_k$ successively.
\end{enumerate}
\end{Definition}

Note that $b_\lambda$ is unique if it exists.

\begin{Definition}
Suppose that a family $\{ \mathcal D(\lambda) \mid \lambda\in P^+ \}$ of
highest weight normal crystals $\mathcal D(\lambda)$ of highest weight
$\lambda$ with $b_\lambda\in\mathcal D(\lambda)$ satisfying the above
properties is given.
It is called {\it closed\/} if the crystal generated by
$b_\lambda\otimes b_\mu$ in $\mathcal D(\lambda)\otimes\mathcal
D(\mu)$ is isomorphic to $\mathcal D(\lambda+\mu)$.
\end{Definition}

We have the following characterization of $\B(\lambda)$ in
\ref{not:Blambda}.
\begin{Proposition}[\protect{\cite[6.4.21]{Joseph}}]\label{prop:closed}
If $\{ \mathcal D(\lambda) \mid \lambda\in P^+ \}$ is a closed family
of highest weight normal crystals, then $\mathcal D(\lambda)$ is
isomorphic to $\B(\lambda)$ as a crystal for any $\lambda\in P^+$.
\end{Proposition}

\subsection{Quantum loop algebra}\label{subsec:qloop}

We briefly recall the notion of quantum loop algebra. See
\cite{CP3,Na-qaff} for detail.

Suppose that a root datum $P$, $P^*$, etc as in \subsecref{subsec:QUE} 
is given. 
Let $\mathbf L\mathfrak g$ be the loop algebra $\mathfrak
g\otimes_{\Q}[z,z^{-1}]$ of the symmetrizable Kac-Moody Lie algebra
$\mathfrak g$.
We define the {\it quantum loop algebra\/} $\Ul$ as a $\Q(q)$-algebra
generated by $e_{k,r}, f_{k,r}$ ($k\in I$, $r\in\Z$), $q^h$ ($h \in P^*$),
$h_{k,m}$ ($k\in I$, $m\in \Z\setminus\{0\}$) with the following
defining relations
{\allowdisplaybreaks[4]
\begin{subequations}
\begin{gather}
%  \text{$q^{\pm c/2}$ is central,} \label{eq:relCcent}\\
%
  q^0 = 1, \quad q^h q^{h'} = q^{h+h'}, \quad
  [q^h, h_{k,m}] = 0, \quad
  [h_{k,m}, h_{l,n}] = 0,
\label{eq:relHH2}\\
  q^h e_{k,r} q^{-h} = q^{\langle h,\alpha_k\rangle} e_{k,r},
  \quad
  q^h f_{k,r} q^{-h} = q^{-\langle h,\alpha_k\rangle} f_{k,r},
  \label{eq:relHE'}
\\
  ( z - q^{\pm \langle h_l,\alpha_k\rangle} w)
    \psi_k^s(z) x_l^\pm(w) =
  ( q^{\pm\langle h_l,\alpha_k\rangle} z -  w)
    x_l^\pm(w) \psi_k^s(z), \quad
    \label{eq:relHE}
\\
%  ( q^{\pm sc/2} z - q^{\pm \langle h_k,\alpha_l\rangle} w)
%    \psi_l^s(z) x_k^\pm(w) =
%  ( q^{\pm\langle h_k,\alpha_l\rangle} q^{\pm sc/2} z -  w)
%    x^\pm_k(w) \psi_l^s(z), \quad
%    \label{eq:relHE}
%\\
%
%  q^h e_k q^{-h} = q^{\langle h, \alpha_k\rangle} e_k, \quad
%  q^h f_k q^{-h} = q^{-\langle h, \alpha_k\rangle} f_k \label{eq:relQHE}\\ 
%
%  \left[h_{l,m}, e_{k,r}\right] = 
%    \frac{[m\langle h_l, \alpha_k\rangle]_{q_l}}{m} e_{k,r+m},\quad
%  \left[h_{l,m}, f_{k,r}\right] =
%    -\frac{[m\langle h_l, \alpha_k\rangle]_{q_l}}{m} f_{k,r+m},
%    \label{eq:relHE}\\
%
  \left[x_{k}^+(z), x_{l}^-(w)\right] =
  \frac{\delta_{kl}}{q_k  - q_k^{-1}}
  \left\{\delta\left(\frac{w}{z}\right)\psi^+_k(w) -
        \delta\left(\frac{z}{w}\right)\psi^-_k(z)\right\},
    \label{eq:relEF}
\\
%  [e_{k,r}, f_{l,s}] = \delta_{kl}
%  \frac{\psi^+_{k,r+s} - \psi^-_{k,r+s}}{q_k  - q_k^{-1}},
%    \label{eq:relEF}
%\\
%
   (z - q^{\pm 2} w) x_{k}^\pm(z) x_{k}^\pm(w)
   = (q^{\pm 2} z - w) x_{k}^\pm(w) x_{k}^\pm(z),
   \label{eq:relexE2}
\\
%   (z - q^2 w) e_{k}(z) e_{k}(w) = (q^2 z - w) e_{k}(w) e_{k}(z),
%   \label{eq:relexE2}
%\\
%
  \prod_{p=1}^{-\langle \alpha_k,h_l\rangle}
   (z - q^{\pm(b'-2p)} w) x_{k}^\pm(z) x_{l}^\pm(w)
 = \prod_{p=1}^{-\langle \alpha_k, h_l\rangle}
   (q^{\pm(b'-2p)} z - w) x_{l}^\pm(w) x_{k}^\pm(z),
   \quad\text{if $k\neq l$},
   \label{eq:relexE1}
\\
%  \prod_{p=1}^{b'} (z - q^{b'-2p} w) e_{k}(z) e_{l}(w)
%  = \prod_{p=1}^{b'} (q^{b'-2p} z - w) e_{l}(w) e_{k}(z),
%   \quad\text{if $k\neq l$}, \label{eq:relexE1}
%\\
%
%   (z - q^{-2} w) f_{k}(z) f_{k}(w) = (q^{-2} z - w) f_{k}(w) f_{k}(z),
%   \label{eq:relexF2}
%\\
%
%  \prod_{p=1}^{b'} (z - q^{2p - b'} w) f_{k}(z) f_{l}(w)
%  = \prod_{p=1}^{b'} (q^{2p - b'} z - w) f_{l}(w) f_{k}(z),
%   \quad\text{if $k\neq l$}, \label{eq:relexF1}
%\\
%
%  e_{k,r+1} e_{l,s} - q^{c_{kl}}e_{l,s}e_{k,r+1}
%  = q^{c_{kl}}e_{k,r} e_{l,s+1} - e_{l,s+1} e_{k,r},\label{eq:relexE}\\
%  f_{k,r+1} f_{l,s} - q^{-c_{kl}}f_{l,s}f_{k,r+1}
%  = q^{-c_{kl}}f_{k,r} f_{l,s+1} - f_{l,s+1} f_{k,r},\label{eq:relexF}\\
%
  \sum_{\sigma\in S_b}
   \sum_{p=0}^{b}(-1)^p 
   \begin{bmatrix} b \\ p\end{bmatrix}_{q_k}
   x_{k}^\pm(z_{\sigma(1)})\cdots x_{k}^\pm(z_{\sigma(p)})
   x_{l}^\pm(w)
   x_{k}^\pm(z_{\sigma(p+1)})\cdots x_{k}^\pm(z_{\sigma(b)}) = 0,
   \quad \text{if $k\neq l$,}
 \label{eq:relDS}
%\\
%  \sum_{\sigma\in S_b}
%   \sum_{p=0}^{b}(-1)^p 
%   \begin{bmatrix} b \\ p\end{bmatrix}_{q_k}
%   e_{k}(z_{\sigma(1)})\cdots e_{k}(z_{\sigma(p)}) e_{l}(w)
%   e_{k}(z_{\sigma(p+1)})\cdots e_{k}(z_{\sigma(b)}) = 0,
%   \quad \text{if $k\neq l$,} \label{eq:relDS_E}
%\\
%
%  \sum_{\sigma\in S_b}
%   \sum_{p=0}^{b}(-1)^p 
%   \begin{bmatrix} b \\ p\end{bmatrix}_{q_k}
%   e_{k,r_{\sigma(1)}}\cdots e_{k, r_{\sigma(p)}} e_{l,s}
%   e_{k, r_{\sigma(p+1)}}\cdots e_{k,r_{\sigma(b)}} = 0,
%   \quad \text{if $k\neq l$,} \label{eq:relDS_E}\\ 
%%
%  \sum_{\sigma\in S_b}
%   \sum_{p=0}^{b}(-1)^p
%   \begin{bmatrix} b\\ p\end{bmatrix}_{q_k}
%   f_{k}(z_{\sigma(1)})\cdots f_{k}(z_{\sigma(p)}) f_{l}(w)
%   f_{k}(z_{\sigma(p+1)})\cdots f_{k}(z_{\sigma(b)}) = 0,
%   \quad \text{if $k\neq l$,} \label{eq:relDS_F} 
\end{gather}
\end{subequations}
where}
$s = \pm$,
$b = 1-\langle h_k, \alpha_l\rangle$, 
$b' = - (\alpha_k,\alpha_l)$,
and 
$S_b$ is the symmetric group of $b$ letters.
Here $\delta(z)$, $x_k^+(z)$, $x_k^-(z)$, $\psi^{\pm}_{k}(z)$ are
generating functions defined by
{\allowdisplaybreaks[4]
\begin{gather*}
   \delta(z) \defeq \sum_{r=-\infty}^\infty z^{r}, \qquad
   x_k^+(z) \defeq \sum_{r=-\infty}^\infty e_{k,r} z^{-r}, \qquad
   x_k^-(z) \defeq \sum_{r=-\infty}^\infty f_{k,r} z^{-r}, \\
   \psi^{\pm}_k(z)
%  = \sum_{r=0}^\infty \psi^{\pm}_{k,\pm r} z^{\mp r}
  \defeq t_k^\pm
   \exp\left(\pm (q_k-q_k^{-1})\sum_{m=1}^\infty h_{k,\pm m} z^{\mp m}\right).
\end{gather*}
%We set $\psi^+_{k,r} = 0$ (resp.\ $\psi^-_{k,r} = 0$) if $r < 0$
%(resp.\ $r > 0$).
We} also need the following generating function
\begin{equation*}
   p_k^\pm(z) \defeq 
   \exp\left(
     - \sum_{m=1}^\infty \frac{h_{k,\pm m}}{[m]_{q_k}} z^{\mp m}
   \right).
\end{equation*}
We have
\(
  \psi^{\pm}_k(z) = t_k^\pm
  p_k^\pm(q_k z)/p_k^\pm(q_k^{-1} z).
\)

There is a homomorphism $\Uq\to \Ul$ defined by
\begin{equation*}
   q^h \mapsto q^h, \quad
   e_k \mapsto e_{k,0}, \quad f_k \mapsto f_{k,0}.
\end{equation*}

Let $\Ul^+$ (resp.\ $\Ul^-$) be the $\Q(q)$-subalgebra of $\Ul$
generated by elements $e_{k,r}$'s (resp.\ $f_{k,r}$'s).  Let $\Ul^0$
be the $\Q(q)$-subalgebra generated by elements $q^h$, $h_{k,m}$.
We have $\Ul = \Ul^+\cdot\Ul^0\cdot\Ul^-$.

Let $e_{k,r}^{(n)} \defeq e_{k,r}^n / [n]_{q_k}!$,
$f_{k,r}^{(n)} \defeq f_{k,r}^n / [n]_{q_k}!$.
Let $\Uli$ be the $\Z[q,q^{-1}]$-subalgebra generated by
$e_{k,r}^{(n)}$, $f_{k,r}^{(n)}$ and $q^h$
for $k\in I$, $r\in \Z$, $h\in P^*$.
The specialization $\Uli\otimes_{\Z[q,q^{-1}]}\C$ with respect to the
homomorphism $\Z[q,q^{-1}]\ni q\mapsto \varepsilon\in\C^*$ is denoted
by $\Ule$ for $\varepsilon\in\C^*$.

Let $\Uli^+$ (resp.\ $\Uli^-$) be $\Z[q,q^{-1}]$-subalgebra generated
by $e_{k,r}^{(n)}$ (resp.\ $f_{k,r}^{(n)}$) for $k\in I$, $r\in \Z$,
$n\in Z_{> 0}$.
Let $\Uli^0$ be the $\Z[q,q^{-1}]$-subalgebra generated by $q^h$, the
coefficients of $p_k^\pm(z)$ and
\begin{equation*}
   \begin{bmatrix}
      q^{h_k}; n \\ r
   \end{bmatrix}
   \defeq
   \prod_{s=1}^r \frac{t_k q_k^{n-s+1} - t_k^{-1} q_k^{-n+s-1}}
   {q_k^s - q_k^{-s}}
\end{equation*}
for all $h\in P$, $k\in I$, $n\in \Z$, $r\in \Z_{> 0}$. We have
$\Uli = \Uli^+\cdot \Uli^0\cdot \Uli^-$ if $\mathfrak g$ is of finite
type (\cite[6.1]{CP3}).

Suppose that a datum
$\left(\lambda,\left(\Psi^{\pm}_{k}(z)\right)_{k\in I}\right)\in
P\times \Q(q)[[z^{\mp}]]^I$ satisfying 
\[
   \Psi_k^+(\infty) = (\alpha_k,\alpha_k)\langle\lambda,h_k\rangle/2,
\quad
   \Psi_k^-(0) = -(\alpha_k,\alpha_k)\langle\lambda,h_k\rangle/2
\]
is given.
We say a $\Ul$-module $M$ is an {\it l--highest weight module\/} (`{\it
l\/}' stands for the loop) with {\it l--highest weight\/}
$\left(\lambda,\left(\Psi^{\pm}_{k}(z)\right)_{k\in I}\right)$
if there exists a vector $m_0\in M$ such that
\begin{subequations}
\begin{align}
  & e_{k,r} m_0 = 0, \qquad \Ul^- m_0 = M, \\
  & q^h m_0 = q^{\langle h, \lambda\rangle} m_0
    \quad\text{for $h\in P^*$}, \qquad
    \psi^{\pm}_{k}(z) m_0 = \Psi^{\pm}_{k}(z) m_0
    \quad\text{for $k\in I$}.
\end{align}
\end{subequations}

A $\Ul$-module $M$ is said to be {\it l--integrable} if the following
two conditions are satisfied.
\begin{aenume}
\item $M$ has a weight space decomposition as a $\Uq$-module:
\begin{equation*}
   M = \bigoplus_{\mu\in P} M_\mu, \quad
   M_\mu \defeq \{ m\mid \text{$
   q^h\cdot v = q^{\langle h, \mu\rangle} m$ for any $h\in P^*$}\}.
\end{equation*}
And $\dim M_\mu < \infty$.
\item For any $m\in M$, there exists $n_0\ge 1$ such that
$e_{k,r_1}\cdots e_{k,r_n} \ast m = 
f_{k,r_1}\cdots f_{k,r_n} \ast m = 0$
for all $r_1,\dots,r_n\in\Z$, $k\in I$ and $n\ge n_0$.
\end{aenume}

We say $\left(\lambda,\left(\Psi^{\pm}_{k}(z)\right)_{k\in I}\right)$
is {\it l--dominant}, if $\lambda\in P^+$ and 
there exists a $I$-tuple of polynomials $P(u) = (P_k(u))_k\in
\Q(q)[u]^I$ with $P_k(0) = 1$ satisfying
\begin{equation}
  \Psi^{\pm}_{k}(z) = q_k^{\deg P_k}
  \left(\frac{P_k(q_k^{-1}/z)}{P_k(q_k/z)}\right)^{\pm},
\end{equation}
where $\left(\ \right)^{\pm} \in\Q(q)[[z^{\mp}]]$ denotes the
expansion at $z = \infty$ and $0$ respectively.

The simple {\it l\/}--highest weight module $M$ is {\it
l\/}--integrable if and only if its {\it l\/}--highest weight
$\left(\lambda,\left(\Psi^{\pm}_{k}(z)\right)_{k\in I}\right)$ is {\it
l\/}--dominant, provided $\mathfrak g$ is of finite type (\cite{CP2})
or $\mathfrak g$ is symmetric (\cite{Na-qaff}). In this case, $P(u)$
is called a {\it Drinfeld polynomial\/} of $M$.
Since the simple {\it l\/}--highest weight module is determined by
$\lambda$ and $P$, we denote it by $L(\lambda,P)$.
%If $\mathfrak g$ is of finite type, then $\Lambda$ is determined from
%$P$, so we simply denote it by $L(P)$.

Let $M$ be a $\Ul$-module with the weight space decomposition $M =
\bigoplus_{\mu\in P} M_\mu$ as a $\Uq$-module such that 
$\dim M_\mu < \infty$.
Since the commutative subalgebra $\Ul^0$ preserves each $M_\mu$, we
can further decompose $M$ into a sum of generalized simultaneous
eigenspaces for $\Ul^0$:
\begin{equation}\label{eq:gen_wt}
   M = \bigoplus M_{\Psi^\pm},
\end{equation}
where $\Psi^\pm(z)$ is a pair $(\mu,(\Psi_k^\pm(z))_k)$ and
\begin{equation*}
   M_{\Psi^\pm} \defeq \left\{ m\in M \left|\,
   \begin{aligned}[c]
     & \text{$q^h \ast m = q^{\langle h,\mu\rangle} m$
     for $h\in P^*$} \\
     & \text{$(\psi_k^\pm(z) - \Psi_k^\pm(z)\operatorname{Id})^N m
     = 0$ for $k\in I$ and sufficiently large $N$}
   \end{aligned}
   \right\}\right. .
\end{equation*}
If $M_{\Psi^\pm}\neq 0$, we call $M_{\Psi^\pm}$ an {\it l--weight
space}, and the corresponding $\Psi^\pm(z)$ an {\it l--weight} of $M$.
This is a refinement of the weight space decomposition.

Let $\lambda\in P^+$ and $\lambda_k = \langle
h_k,\lambda\rangle\in\Z_{\ge 0}$.
Let $G_{\lambda} = \prod_{k\in I} \GL(\lambda_k,\C)$.
Its representation ring
$R(G_{\lambda})$ is the invariant part of the Laurant polynomial ring:
\begin{equation*}
   R(G_{\lambda})
   = \Z[x_{1,1}^\pm,\dots, x_{1,\lambda_1}^\pm]^{\mathfrak S_{\lambda_1}}
   \otimes \Z[x_{2,1}^\pm,\dots,
     x_{2,\lambda_2}^\pm]^{\mathfrak S_{\lambda_2}}
   \otimes\cdots\otimes
     \Z[x_{n,1}^\pm,\dots, x_{n,\lambda_n}^\pm]^{\mathfrak S_{\lambda_n}},
\end{equation*}
where we put a numbering $1,\dots,n$ to $I$.
In \cite{Na-qaff}, when $\mathfrak g$ is symmetric, we constructed a
$\Uli\otimes_{\Z} R(G_{\lambda})$-module $M(\lambda)$ such that it is
{\it l\/}--integrable and has a vector $[0]_{\lambda}$ satisfying
\begin{subequations}
\begin{gather}
   e_{k,r}[0]_\lambda = 0\quad\text{for any $k\in I$, $r\in \Z$},
\\
   M(\lambda) = 
   \left(\Uli^-\otimes_{\Z} R(G_\lambda)\right)[0]_{\lambda},
   \label{eq:span}
\\
   q^h [0]_\lambda = q^{\langle h,\lambda\rangle} [0]_\lambda,
\\
   \psi_k^\pm(z)[0]_\lambda
   = q^{w_k} \left( \prod_{i=1}^{w_k}
     \frac{1-q^{-1}x_{k,i}/z}{1-q x_{k,i}/z}\right)^\pm
   [0]_\lambda.
\end{gather}
\end{subequations}
We call this a {\it universal standard module}. Its construction,
via quiver varieties, will be explained briefly in \secref{sec:Ul}.

If an $I$-tuple of monic polynomials $P(u) =
(P_k(u))_{k\in I}$ with $\deg P_k = \langle h_k,\lambda\rangle$ are given,
then we define a {\it standard module\/} by the specialization
\begin{equation*}
   M(\lambda,P) = M(\lambda)\otimes_{R(G_\lambda)[q,q^{-1}]} \C(q),
\end{equation*}
where the algebra homomorphism $R(G_\lambda)[q,q^{-1}]\to \C(q)$ sends
$x_{k,1},\dots, x_{k,w_k}$ to roots of $P_k$.
The simple module $L(\lambda,P)$ is the simple quotient of
$M(\lambda,P)$.

It is known that $M(\lambda)$ is free as an $R(G_\lambda)$-module
(\cite[7.3.5]{Na-qaff}). Thus $M(\lambda,P)$ depends on $P$ {\it
continuously}, while $L(\lambda,P)$ depends {\it discontinuously}.

\subsection{Drinfeld realization}
Now we assume that $\mathfrak g$ is of finite type, i.e., $(\ ,\ )$ is
positive definite. We take a root datum so that $\rank P = \rank
\mathfrak g$. It is well-known that the untwisted affine Lie algebra
$\widehat{\mathfrak g} = \mathbf L\mathfrak g\oplus\Q c\oplus \Q d$ is
a Kac-Moody Lie algebra with a root datum
\begin{gather*}
   \widehat P^* = P^*\oplus \Z c\oplus \Z d, \quad
   \widehat I = I\cup \{ 0\},
\end{gather*}
and certain $\alpha_0$, $h_0$, $(\ ,\ )$.
The {\it quantum affine algebra\/} $\Ua$ is the quantum enveloping
algebra associated with this root datum. Let $\Ua'$ be the subalgebra
of $\Ua/(q^c - 1)$ generated by $e_k$, $f_k$ ($k\in \widehat I$) and
$q^h$ ($h\in P$). Then a comultiplication $\Delta$ is defined on
$\Ua'$ by the same formula \eqref{eq:comul}. The integral form $\Uai'$ 
is defined in the same way.

In \cite{Drinfeld} Drinfeld observed that there exists an isomorphism
$\Ul\to \Ua'$ of $\Q(q)$-algebras. Its explicit form was given by Beck 
\cite{Beck}:
\begin{equation*}
   e_{k,r} = o(k)^r T_{\check\omega_k}^{-r}(e_k), \quad
   f_{k,r} = o(k)^r T_{\check\omega_k}^{r}(f_k),
\end{equation*}
where $o\colon I\to \{\pm 1\}$ is an orientation of $I$ such that
$o(k) = -o(l)$ if $a_{kl}\neq 0$ for $k\neq l$, and
$T_{\check\omega_k}$ is an automorphism of $\Ua'$ defined by Lusztig
(the braid group action).
Since $T_{\check\omega_k}$ preserves $\Uai'$, the isomorphism induces
an isomorphism $\Uli\to\Uai'$.

We shall identify $\Ua'$ with $\Ul$ when $\mathfrak g$ is of finite
type hereafter. In particular, we have a comultiplication $\Delta$ on
$\Ul$. We need the following asymptotic formula, which can be deduced
from \cite{Da,JKK}.

\begin{Lemma}\label{lem:comult}
\rom{(1)} On finite dimensional $\Ul$-modules, we have
\begin{equation*}
   \Delta(h_{k,\pm m}) 
   = h_{k,\pm m}\otimes 1 + 1 \otimes h_{k,\pm m}
   + \text{a nilpotent term}.
\end{equation*}

\rom{(2)} Let $V$ and $W$ be finite dimensional
$\Ul$-modules. Suppose that $V$ has a vector $b$ such that
$e_{k,r}b = 0$ for all $k\in I$, $r\in\Z$. Then
a subspace $\{b\}\otimes W\subset V\otimes W$ is invariant under
$\Ul^+$, and the map $W\ni x\mapsto b\otimes x\in \{b\}\otimes W$
respects $\Ul^+$-module structure.
\end{Lemma}

\begin{Remark}
The property (2) depends on our choice of the comultiplication
\eqref{eq:comul}. If we take one in \cite{Lu-book}, then the property
holds after we exchange $V$ and $W$.
\end{Remark}

Before we close this section, we give an algebraic characterization of
the standard module.
We assume $\mathfrak g$ is of type $ADE$. For each fundamental weight
$\Lambda_k$, the universal standard module $M({\Lambda_k})$ is a
$\Uli\otimes_{\Z[q,q^{-1}]}R(\C^*\times\C^*) \cong
\Uli[x,x^{-1}]$-module. We set $W({\Lambda_k}) =
M({\Lambda_k})/(x-1)M({\Lambda_k})$. It is called an {\it
  l--fundamental representation}.

\begin{Theorem}\label{thm:alg_std}
Put a numbering $1,\dots,n$ on $I$.
Let $\lambda_k = \langle h_k,\lambda\rangle$.
The universal standard module $M(\lambda)$ is the submodule of
\[
   W({\Lambda_1})^{\otimes \lambda_1}
%   \otimes W({\Lambda_2})^{\otimes \lambda_2}
   \otimes \cdots\otimes
   W({\Lambda_n})^{\otimes \lambda_n}
   \otimes
   \Z[q,q^{-1},x_{1,1}^\pm,\dots, x_{1,\lambda_1}^\pm,
%    x_{2,1}^\pm,\dots, x_{2,\lambda_2}^\pm,
    \cdots,
    x_{n,1}^\pm,\dots, x_{n,\lambda_n}^\pm]^{\mathfrak S_{\lambda_1}\times
    \cdots\times\mathfrak S_{\lambda_n}}
\]
\rom(the tensor product is over $\Z[q,q^{-1}]$\rom)
generated by
\(
   \bigotimes_{k\in I} [0]_{\Lambda_k}^{\otimes \lambda_k}
\).
\rom(The result holds for the tensor product of any order.\rom)
\end{Theorem}

This theorem was not explicitly stated in \cite{Na-qaff}, but can be
deduced by above properties as follows:
First we make tensor products of both modules by the quotient field of
$R(G_\lambda)[q,q^{-1}]$. Then both modules are simple and have the
same Drinfeld polynomial (cf.\ proof of
\lemref{lem:general_point}). Thus there exists a unique isomorphism
from $M(\lambda)\otimes\text{(quotient field)}$ to the tensor product,
sending $[0]_\lambda$ to
\(
   \bigotimes_{k\in I} [0]_{\Lambda_k}^{\otimes \lambda_k}
\).
Now by the property \eqref{eq:span} and the freeness of $M(\lambda)$,
we have the assertion.

\begin{Remark}
Varagnolo-Vasserot \cite[\S7]{VV} conjectured that $M(\lambda)$ is
isomorphic to a module studied by Kashiwara~\cite{Kas2} ($V(\lambda)$
in his notation), after tensoring $\Q$ and forgetting the symmetric
group invariance.
When $\lambda$ is a fundamental weight, both modules are isomorphic,
since they are simple and have the same Drinfeld polynomial.
Kashiwara conjectures that his module $V(\lambda)$ has the property in
\thmref{thm:alg_std} [loc.\ cit., \S13]. Thus Varagnolo-Vasserot's
conjecture is equivalent to Kashiwara's conjecture.
On the other hand, Kashiwara shows that the submodule above has a
global crystal base [loc.\ cit., Theorem 8.5]. ($N_\Q$ in his
notation.)
Probably, this base coincides with the (conjectural) canonical base
considered in \cite{Lu-rem} as analogue of \cite{Lu-Base,Lu-Base2}.
\end{Remark}

\section{Preliminaries~(II) -- geometric part}

In this paper, all varieties are defined over $\C$.

\subsection{$K$-homology groups}
Let $X$ be a quasi-projective variety. Its integral Borel-Moore
homology group of degree $k$ is denoted by $H_k(X,\Z)$. Set $H_*(X,\Z)
= \bigoplus_k H_k(X,\Z)$.
When a linear algebraic group $G$ acts algebraically on
$X$, we denote by $K^G(X)$ the Grothendieck group of the abelian
category of $G$-equivariant coherent sheaves on $X$ (the $K$-homology
group). It is a module over $R(G)$, the representation ring of $G$. We 
shall use several operations on $H_*(X,\Z)$ and $K^G(X)$ in this
article, but we do not review them here. See \cite[\S6]{Na-qaff} and
\cite{Gi-book}.

\subsection{Quiver variety}
We briefly review the notion of quiver varieties. The reference to
results can be found in \cite[\S2]{Na-quiver}, unless refered explictly.

Suppose a root datum is given. We assume that it is symmetric, i.e.,
$(\alpha_k,\alpha_k) = 2$ for all $k\in I$. Then the generalized
Cartan matrix $\bC$ is equal to $((\alpha_k,\alpha_l))_{k,l\in I}$ and
symmetric.
To the root datum, we associate a finite graph $(I,E)$ as follows (the
Dynkin diagram). The set of vertices is identified with $I$, and we
drow $(\alpha_k,\alpha_l)$ edges between vertices $k$ and $l$ ($k\neq
l$). We give no edge loops, edge joining a vertex with itself.
Conversely a finite graph without edge loops gives a {\it symmetric\/} 
generalized Cartan matrix.

Let $H$ be the set of pairs consisting of an edge together with its
orientation. For $h\in H$, we denote by $\vin(h)$ (resp.\ $\vout(h)$)
the incoming (resp.\ outgoing) vertex of $h$.
For $h\in H$ we denote by $\overline h$ the same edge as $h$ with the
reverse orientation.
Choose and fix a numbering $1,2,\dots, n$ on $I$. We then define a
subset $\Omega\subset H$ so that $h\in\Omega$ if $\vout(h) < \vin(h)$.
Then $\Omega$ is an orientation of the graph, i.e.,
$\overline\Omega\cup\Omega = H$, $\Omega\cap\overline\Omega =
\emptyset$. (The numbering and the orientation will not play a role
until \secref{sec:comb}.)
The pair $(I,\Omega)$ is called a {\it quiver}.

If $V^1 = \bigoplus_k V^1_k$ and $V^2 = \bigoplus_k V^2_k$ are
$I$-graded vector spaces, we define vector spaces by
%\begin{gather*}
%  \HomL(V^1, V^2) \defeq
%  \bigoplus_{k\in I} \Hom(V^1_k, V^2_k), \quad
%  \HomE(V^1, V^2) \defeq
%  \bigoplus_{h\in H} \Hom(V^1_{\vout(h)}, V^2_{\vin(h)})
%, \\
%  \HomE_\Omega(V^1, V^2) \defeq
%  \bigoplus_{h\in \Omega} \Hom(V^1_{\vout(h)}, V^2_{\vin(h)}), \quad
%  \HomE_{\overline\Omega}(V^1, V^2) \defeq
%  \bigoplus_{h\in \overline\Omega} \Hom(V^1_{\vout(h)}, V^2_{\vin(h)}).
%\end{gather*}
\begin{equation}\label{eq:LE}
  \HomL(V^1, V^2) \defeq
  \bigoplus_{k\in I} \Hom(V^1_k, V^2_k), \quad
  \HomE(V^1, V^2) \defeq
  \bigoplus_{h\in H} \Hom(V^1_{\vout(h)}, V^2_{\vin(h)})
\end{equation}

For $B = (B_h) \in \HomE(V^1, V^2)$ and 
$C = (C_h) \in \HomE(V^2, V^3)$, let us define a multiplication of $B$
and $C$ by
\[
  CB \defeq (\sum_{\vin(h) = k} C_h B_{\overline h})_k \in
  \HomL(V^1, V^3).
\]
Multiplications $ba$, $Ba$ of $a\in \HomL(V^1,V^2)$, $b\in\HomL(V^2,
V^3)$, $B\in \HomE(V^2, V^3)$ are defined in obvious manner. If
$a\in\HomL(V^1, V^1)$, its trace $\tr(a)$ is understood as $\sum_k
\tr(a_k)$.

If $V$ and $W$ are $I$-graded vector spaces, we consider the vector
spaces
\begin{equation}\label{def:bM}
  \bM \equiv \bM(V, W) \defeq
  \HomE(V, V) \oplus \HomL(W, V) \oplus \HomL(V, W),
\end{equation}
where we use the notation $\bM$ unless we want to specify $V$, $W$.
The above three components for an element of $\bM$ is denoted by $B$,
$i$, $j$ respectively.

\begin{Convention}\label{convention:weight}
When quiver varieties will be related the representation theory, we
will choose $V$ and $W$ corresponding to a pair $(\bv,\bw)\in
Q^+\times P^+$.
The rule is
\(
    \dim V_k = v_k
\),
\(
   \dim W_k = \langle h_k,\bw\rangle
\),
where $\bv = \sum_k v_k \alpha_k$.
Conversely $V$ determines $\bv$, while $W$ determines $\bw$ modulo an
element $*$ such that $\langle h_k,*\rangle = 0$ for all $h_k$.
But the action of $\mathfrak g$ on simple highest weight modules
$V(\bw)$ and $V(\bw+*)$ differ only by scalars (see
\cite[9.10]{Kac}). So essentially there is no umbiguity.
\end{Convention}

For an $I$-graded subspace $S = \bigoplus_k S_k$ of subspaces $V$ and
$B\in\HomE(V, V)$, we say $S$ is {\it $B$-invariant\/} if
$B_h(S_{\vout(h)}) \subset S_{\vin(h)}$.

Fix a function $\varepsilon\colon H \to \C^*$ such that
$\varepsilon(h) + \varepsilon(\overline{h}) = 0$ for all $h\in H$.
%
%In \cite{Na-quiver,Na-alg}, it was assumed that $\varepsilon$ takes
%its value $\pm 1$, but this assumption is not necessary as remarked by 
%Lusztig~\cite{Lu-qv}.
%
For $B\in\HomE(V^1, V^2)$, let us denote
by $\varepsilon B\in \HomE(V^1, V^2)$ data given by $(\varepsilon B)_h 
= \varepsilon(h) B_h$ for $h\in H$.

Let us define a symplectic form $\omega$ on $\bM$ by
\begin{equation}
        \omega((B, i, j), (B', i', j'))
        \defeq \tr(\varepsilon B\, B') + \tr(i j' - i' j).
\label{def:symplectic}\end{equation}

Let
\(
   G_V \defeq \prod_k \GL(V_k).
\)
It acts on $\bM$ by
\begin{equation}\label{eq:Gaction}
  (B, i, j) \mapsto g\cdot (B, i, j) \defeq (g B g^{-1}, g i, j g^{-1})
\end{equation}
preserving the symplectic form $\omega$. The moment map
$\mu\colon\bM\to\HomL(V, V)$ vanishing at the origin is given by
\begin{equation}\label{eq:mu}
  \mu(B, i, j) = \varepsilon B B + i j,
\end{equation}
where the dual of the Lie algebra of $G_V$ is identified with the Lie
algebra via the trace.
Let $\mu^{-1}(0)$ be an affine algebraic variety (not necessarily
irreducible) defined as the zero set of $\mu$.
%The equation $\mu = 0$ will be called the {\it ADHM equation\/}.

\begin{Definition}\label{def:stable}
A point $(B, i, j) \in \mu^{-1}(0)$ is said to be {\it stable\/} if 
the following condition holds:
\begin{itemize}
\item[] if an $I$-graded subspace $S = \bigoplus_k S_k$ of $V$ is
$B$-invariant and contained in $\Ker j$, then $S = 0$.
\end{itemize}
Let us denote by $\mu^{-1}(0)^{\operatorname{s}}$ the set of stable points.
\end{Definition}
Clearly, the stability condition is invariant under the action of
$G_V$. Hence we may say an orbit is stable or not.

Let
\(
   \M \equiv \M(\bv,\bw)
   \defeq \mu^{-1}(0)^{\operatorname{s}}/G_V.
\)
We use the notation $\M$ unless we need to specify dimensions of $V$
and $W$.
It is known that the $G_V$-action is free on
$\mu^{-1}(0)^{\operatorname{s}}$ and $\M$ is a nonsingular
quasi-projective variety, having a symplectic form induced by
$\omega$. A $G_V$-orbit though $(B,i,j)$, considered as a point of
$\M$ is denoted by $[B,i,j]$.
Since the action is free, $V$ and $W$ can be considered as $I$-graded
vector bundles over $\M$. We denote them by the same notation.
We consider $\HomE(V,V)$, $\HomL(W,V)$, $\HomL(V, W)$ as vector
bundles defined by the same formula as in \eqref{eq:LE}. By the
definition, $B$, $i$, $j$ can be considered as sections of those
bundles.

Let us consider three-term sequence of vector bundles over
$\M(\bv,\bw)$ given by
\begin{equation}\label{eq:quiver_tangent}
        \HomL(V, V)
        \overset{\iota}{\longrightarrow}
        \HomE(V, V) \oplus \HomL(W, V) \oplus \HomL(V,W)
        \overset{d \mu}{\longrightarrow}
        \HomL(V, V),
\end{equation}
where $d\mu$ is the differential of $\mu$ at $(B,i,j)$, i.e.,
\begin{equation*}
     d\mu(C,I,J) = \varepsilon B C + \varepsilon C B + iJ + I j,
\end{equation*}
and $\iota$ is given by
\begin{equation*}
        \iota(\xi) = (B \xi - \xi B) \oplus
         (-\xi i) \oplus j \xi.
\end{equation*}
Then $\iota$ is injective and $d\mu$ is surjective, and 
the tangent bundle of $\M(\bv,\bw)$ is identified with
$\Ker d\mu/\Ima\iota$.

Let
\(
   \M_0 \equiv \M_0(\bv,\bw)
   \defeq \mu^{-1}(0)\dslash G_V,
\)
where $\dslash$ is the affine algebro-geometric quotient, i.e., the
coordinate ring of $\mu^{-1}(0)\dslash G_V$ is $G_V$-invariant
polynomials on $\mu^{-1}(0)$. It is an affine algebraic variety, and
identified with the set of closed $G_V$-orbits in $\mu^{-1}(0)$ as a 
set.

There exists a projective morphism
\(
   \pi\colon \M \to \M_0,
\)
sending $[B,i,j]$ to the unique closed orbit contained in the closure
of the orbit $G_V\cdot(B,i,j)$.

Let $\La \equiv \La(\bv,\bw) \defeq \pi^{-1}(0)$. It is a lagrangian
subvariety in $\M(\bv,\bw)$.

If $\bv'-\bv\in Q^+$, then $\M_0(\bv,\bw)$ can be identified with a
closed subvariety of $\M_0(\bv',\bw)$. We consider the direct limit
\(
   \M_0(\infty,\bw) \defeq \bigcup_\bv \M_0(\bv,\bw)
\).
If the graph is of finite type, $\M_0(\bv,\bw)$ stabilizes at some
$\bv$.
This is {\it not\/} true in general. However, it has no harm in
this paper. We use $\M_0(\infty,\bw)$ to simplify the notation, and do
not need any structures on it.
We can always work on individual $\M_0(\bv,\bw)$,
not on $\M_0(\infty,\bw)$.

We set
\(
    \Mw \defeq \bigsqcup_{\bv}\M(\bv,\bw),
\)
\(
   \Law \defeq \bigsqcup_{\bv}\La(\bv,\bw).
\)
They may have infinitely many components, but ho harm as above.

Let $\Delta(\bv,\bw)$ denote the diagonal in
$\M(\bv,\bw)\times\M(\bv,\bw)$. It is a lagrangian subvariety, if we
endow a symplectic form $\omega\times(-\omega)$ on
$\M(\bv,\bw)\times\M(\bv,\bw)$.

For $n\in\Z_{>0}$, we define $\Pa_k^{(n)}(\bv,\bw)$ by
\begin{equation}
\label{eq:Pa^{(n)}}
   \Pa_k^{(n)}(\bv,\bw) \defeq
   \{ (B,i,j,S)\mid
   \text{$(B,i,j)\in\bM(V,W)$, $S\subset V$ as below}\}
   /G_V,
\end{equation}
\begin{aenume}
\item $(B, i, j)\in\mu^{-1}(0)^{\operatorname{s}}$,
\item $S$ is a $B$-invariant subspace containing the image of $i$ with
         $\dim S = \bv - n\alpha_k$.
\end{aenume}
When $n = 1$, we simply denote it by $\Pa_k(\bv,\bw)$.
If we set $\bv' = \bv-n\alpha_k$, we have a natural morphism
$\Pa_k^{(n)}(\bv,\bw)\to \M(\bv',\bw)\times\M(\bv,\bw)$ by
\begin{align*}
  & [B,i,j,S] \longmapsto ([B',i',j'], [B,i,j]),\\
  & \qquad \text{where $(B',i',j')$ is the restriction
                               of $(B,i,j)$ to $S$}.
\end{align*}
Then $\Pa_k^{(n)}(\bv,\bw)$ is a nonsingular closed lagrangian
subvariety in $\M(\bv',\bw)\times\M(\bv,\bw)$
(\cite[11.2.3]{Na-qaff}). The quotient $V/V'$ defines a rank $n$
vector bundle over $\Pa_k^{(n)}(\bv,\bw)$.

Let $G_W = \prod_{k\in I} \GL(W_k)$. It acts naturally on $\bM$,
$\M$ and $\M_0$.
We define $\C^*$-actions on $\M$ and $\M_0$ by
\begin{equation*}%\label{eq:C*action}
    B_{h} \mapsto t^{m(h)+1} B_{h}, \quad
    i \mapsto t i, \quad j \mapsto t j \qquad
    (t\in \C^*),
\end{equation*}
where $m\colon H\to \Z$ is a certain function determined by a
numbering on edges joining common vertices (see
\cite[2.7]{Na-qaff}). When the root datum is simply-laced, i.e.,
symmetric and $(\alpha_k,\alpha_l) \in \{0,1\}$ for $k\neq l$, we have
$m\equiv 0$.
We denote this action by 
$(B,i,j)\mapsto g\ast (B,i,j)$ or $[B,i,j]\mapsto g\ast [B,i,j]$ for
$g\in G_W\times \C^*$.
The vector bundles $V$, $W$, $\HomE(V,V)$, $\HomL(W,V)$, $\HomL(V, W)$
are $G_W\times\C^*$-equivariant bundles, and $B$, $i$, $j$ are
equivariant sections.

Let $L(m)$ be the $1$-dimensional $\C^*$-module defined by $t\mapsto
t^m$ for $m\in\Z$. For a $\C^*$-module $V$, $L(m)\otimes V$ is
denoted by $q^m V$. That is $L(1)$ is identified with $q$.

We consider the following $G_W\times\C^*$-equivariant complex 
$C_k^\bullet$ over $\M$:
\begin{equation}\label{eq:taut_cpx}
C_k^\bullet\equiv C_k^\bullet(\bv,\bw)\colon
\begin{CD}
  q^{-1} V_k
  @>{\sigma_k}>>
 \displaystyle{\bigoplus_{l:k\neq l}}
     [-\langle h_k,\alpha_l\rangle]_q V_l
    \oplus W_k
  @>{\tau_k}>>
  qV_k,
\end{CD}
\end{equation}
where
\begin{equation*}
\sigma_k = \bigoplus_{\vin(h)=k} B_{\overline h} \oplus j_k, \qquad
\tau_k = \sum_{\vin(h)=k} \varepsilon(h) B_h + i_k.
\end{equation*}
We assign degree $0$ to the middle term.

Let $Q_k(\bv,\bw)$ the degree $0$ cohomology of the
complex~\eqref{eq:taut_cpx}, i.e.,
\begin{equation*}
    Q_k(\bv,\bw) \defeq 
    \Ker \tau_k / \Ima \sigma_k
%    \Ker\left(\sum_{\vin(h) = k} \varepsilon(h) B_h + i_k\right)\Biggm/
%    \Ima\left(\bigoplus_{\vout(h) = k} B_h\oplus j_k\right)
.
\end{equation*}

We introduce the following subsets of $\M(\bv,\bw)$:
\begin{equation}\label{taut:subset}
\begin{gathered}
        \M_{k;n}(\bv,\bw) \defeq
        \left\{ [B,i,j]\in \M(\bv,\bw) \Biggm| \codim_{V_k}
%        \left( \sum_{\vin(h) = k} \Ima B_h + \Ima i_k\right) 
        \Ima \tau_k
        = n \right\} \\
        \M_{k;\le n}(\bv,\bw) \defeq
        \bigcup_{m\le n} \M_{k;m}(\bv,\bw), \qquad
        \M_{k;\ge n}(\bv,\bw) \defeq
        \bigcup_{m\ge n} \M_{k;m}(\bv,\bw).
\end{gathered}
\end{equation}
Since $\M_{k;\le n}(\bv,\bw)$ is an open subset of $\M(\bv,\bw)$,
$\M_{k;n}(\bv,\bw)$ is a locally closed subvariety.
The restriction of $Q_k(\bv,\bw)$ to $\M_{k;n}(\bv,\bw)$ is a
$G_W\times\C^*$-equivariant vector bundle of rank $\langle h_k, \bw -
\bv\rangle + n$.

Replacing $V_k$ by $\Ima\tau_k$, we have a natural map
\begin{equation}\label{taut:induct}
  p\colon \M_{k;n}(\bv,\bw) \to \M_{k;0}(\bv - n\alpha_k,\bw).
\end{equation}
Note that the projection $\pi\colon \M(\bv,\bw)\to \M_0(\bv,\bw)$
factors through $p$. In particular, the fiber of $\pi$ is preserved
under $p$.

\begin{Proposition}\label{prop:taut_fib}
Let $G(n,Q_k(\bv - n\alpha_k,\bw)|_{\M_{k;0}(\bv-n\alpha_k,\bw)})$ be
the Grassmann bundle of $n$-planes in
the vector bundle 
%\rom(of rank $\langle h_k, \bw - \bC\bv\rangle + 2n$\rom)
obtained by restricting $Q_k(\bv - n\alpha_k,\bw)$ to
$\M_{k;0}(\bv-n\alpha_k,\bw)$.
Then
% both $\Pa_k^{(n)}(\bv,\bw)\cap
%(\M(\bv-n\alpha_k,\bw)\times\M_{k;\le n}(\bv,\bw))$ and
%$G(n,Q_k(\bv - n\alpha_k,\bw)|_{\M_{k;0}(\bv-n\alpha_k,\bw)})$
%are isomorphic to $\M_{k;n}(\bv,\bw)$.
%
%More precisely
we have the following diagram:
\begin{equation*}
\begin{CD}
   G(n,Q_k(\bv - n\alpha_k,\bw)|_{\M_{k;0}(\bv-n\alpha_k,\bw)})
     @>\Pi>> \M_{k;0}(\bv - n\alpha_k,\bw) \\
   @VV{\cong}V @| \\
   \M_{k;n}(\bv,\bw) @>p>> \M_{k;0}(\bv - n\alpha_k,\bw),
\end{CD}
\end{equation*}
where $\Pi$ is the natural projection.
The kernel of the natural surjective homomorphism
\(
   p^* Q_k(\bv - n\alpha_k,\bw) \to Q_k(\bv,\bw)
\)
is isomorphic to the tautological vector bundle of the Grassmann
bundle of the first row.
\end{Proposition}

The following formula will play an important role later.
\begin{equation}\label{eq:dim_fiber}
   \dim \text{fiber of $p$}
   = n(\langle h_k, \bw-\bv\rangle + n)
   = \frac 12 \left(\dim \M(\bv,\bw) - \dim \M(\bv-n\alpha_k,\bw)\right).
\end{equation}

\section{Varieties $\Zm$, $\Zl$ and their equivariant
   $K$-theories}\label{sec:vardef}
 
The main body of this article starts from this section. For the sake
of space, we only consider the case of tensor products of {\it
two\/} modules. The arguments can be generalized to the case of $N$
modules in a straightforward way. We will mention in
\secref{sec:general}.
 
Let $\bw, \bw^1, \bw^2\in P^+$ be dominant weights such that $\bw =
\bw^1 + \bw^2$. These will be fixed until \secref{sec:general}.

Let us fix a direct sum decomposition $W = W^1\oplus W^2$ of
$I$-graded vector spaces with $\langle h_k,\bw^1\rangle = \dim W^1_k$,
$\langle h_k, \bw^2\rangle = \dim W^2_k$. Set
$G_{W^1} = \prod_k \GL(W_k^1)$, $G_{W^2} = \prod_k \GL(W_k^2)$.

We define a three-term sequence of vector bundles
over $\M(\bv^1,\bw^1)\times \M(\bv^2, \bw^2)$ by
\begin{equation}\label{eq:hecke_complex}
%\begin{CD}
  \HomL(V^1, V^2)
  \overset{\alpha^{21}}{\longrightarrow}
    \HomE(V^1,V^2) \oplus 
    \HomL(W^1, V^2) \oplus q \HomL(V^1,W^2)
%  \left(\begin{matrix}
%    \HomE(V^1,V^2) \\
%    \oplus \\
%    \HomL(W^1, V^2) \oplus \HomL(V^1,W^2)
%  \end{matrix}\right)
  \overset{\beta^{21}}{\longrightarrow}
  \HomL(V^1, V^2),
%\end{CD}
\end{equation}
where
\begin{equation*}
\begin{split}
        \alpha^{21}(\xi) & = (B^2 \xi - \xi B^1) \oplus
         (-\xi i^1) \oplus j^2 \xi, \\
        \beta^{21}(C\oplus I\oplus J)
         &= \varepsilon B^2 C + \varepsilon C B^1 + i^2 J + I j^1.
\end{split}
\end{equation*}
This is a complex, that is $\beta^{21}\alpha^{21} = 0$, thanks to the
equation $\varepsilon B^pB^p+i^pj^p = 0$ ($p=1,2$). It is $\C^*\times
G_{W^1}\times G_{W^2}$-equivariant.

By the same argument as in \cite[3.10]{Na-alg}, $\alpha^{21}$ is
injective and $\beta^{21}$ is surjective. Thus the quotient
$\Ker\beta^{21}/\Ima\alpha^{21}$ is a vector bundle over
$\M(\bv^1,\bw^1)\times \M(\bv^2, \bw^2)$ with rank
\(
   (\bv^1, \bw^2) + (\bw^1, \bv^2) - (\bv^1, \bv^2)
\).

We define a one parameter subgroup $\lambda\colon \C^*\to G_W$ by
\begin{equation*}
   \lambda(t) = \id_{W^1} \oplus t \id_{W^2}
   \in G_{W^1}\times G_{W^2} \subset G_W.
\end{equation*}

\begin{Lemma}\label{lem:fixed}
The fixed point set of $\lambda(\C^*)$ in $\M(\bv,\bw)$ is isomorphic
to $\bigsqcup_{\bv^1+\bv^2=\bv} \M(\bv^1,\bw^1)\times \M(\bv^2,\bw^2)$.
\end{Lemma}

\begin{proof}
A point $[B,i,j]\in\M(\bv,\bw)$ is fixed by $\lambda(\C^*)$ if and
only if there exists a one parameter subgroup $\rho\colon \C^*\to
G_V$ such that
\begin{equation*}
   \lambda(t)\ast (B,i,j) = \rho(t)^{-1} \cdot (B,i,j).
\end{equation*}
Let $V^1$ (resp.\ $V^2$) be the eigenspace of $V$ with eigenvalue
$1$ (resp.\ $t$). Let $V'$ be the sum of other eigenspaces. The above
equation implies that 
\begin{enumerate}
\item $B(V^1) \subset V^1$, $B(V^2)\subset V^2$, $B(V')\subset V'$,
\item $i(W^1)\subset V^1$, $i(W^2)\subset V^2$,
\item $j(V^1)\subset W^1$, $j(V^2)\subset W^2$, $j(V')= 0$.
\end{enumerate}
The stability condition implies that $V' = 0$. Thus we have $V =
V^1\oplus V^2$, and $[B,i,j]$ decomposes into a sum
$[B^1,i^1,j^1]\in\M(\bv^1,\bw^1)$ and $[B^2,i^2,j^2]\in\M(\bv^2,\bw^2)$,
where $\bv^1 = \dim V^1$, $\bv^2 = \dim V^2$.

Conversely $([B^1,i^1,j^1],
[B^2,i^2,j^2])\in\M(\bv^1,\bw^1)\times\M(\bv^2,\bw^2)$ defines a fixed
point. Thus we have a surjective morphism
\begin{equation*}
   \bigsqcup_{\bv^1+\bv^2=\bv} \M(\bv^1,\bw^1)\times \M(\bv^2,\bw^2)
   \to
   \M(\bv,\bw)^{\lambda(\C^*)}.
\end{equation*}
By the freeness of the $G_V$-action on
$\mu^{-1}(0)^{\operatorname{s}}$, this is injective.

Let us identify the tangent bundle of $\M(\bv,\bw)$ with $\Ker
d\mu/\Ima\iota$ in \eqref{eq:quiver_tangent}. Its restriction to
the fixed point set $\M(\bv^1,\bw^1)\times \M(\bv^2,\bw^2)$ decomposes 
as
\begin{equation*}
   \left(\Ker d\mu^1/\Ima\iota^1\right)
   \oplus\left(\Ker d\mu^2/\Ima\iota^2\right)
   \oplus \left(\Ker \beta^{21}/\Ima\alpha^{21}\right)
   \oplus \left(\Ker \beta^{12}/\Ima\alpha^{12}\right),
\end{equation*}
where 
$\iota^p$, $d\mu^p$ are as in \eqref{eq:quiver_tangent} with
$\M(\bv,\bw)$ replaced by $\M(\bv^p,\bw^p)$ ($p=1,2$),
$\alpha^{21}$, $\beta^{21}$ are as above, and $\alpha^{12}$,
$\beta^{12}$ are defined by exchanging the role of $V^1$, $W^1$ and
$V^2$, $W^2$.
We let $\lambda(\C^*)$ acts on $V^1$ (resp.\ $V^2$) with weight $0$
(resp.\ $1$). Then this identification respects the
$\lambda(\C^*)$-action. Thus the tangent space of the fixed point
component, which is the $0$-weight space of the whole tangent space,
is identified with
\(
   \left(\Ker d\mu^1/\Ima\iota^1\right)
   \oplus\left(\Ker d\mu^2/\Ima\iota^2\right)
\).
It is isomorphic to the tangent space of
$\M(\bv^1,\bw^1)\times\M(\bv^2,\bw^2)$.
Therefore the map is isomorphism.
\end{proof}

We move all $\bv$, $\bv^1$, $\bv^2$. We get
\[
   \M(\bw^1)\times\M(\bw^2) \cong \M(\bw)^{\lambda(\C^*)}.
\]

There are a natural $R(G_{W^1}\times
G_{W^2}\times\C^*)$-homomorphism
\begin{multline}
   \boxtimes\colon
   K^{G_{W^1}\times\C^*}(\M(\bv^1,\bw^1))\otimes_{R(\C^*)}
    K^{G_{W^2}\times\C^*}(\M(\bv^2,\bw^2))\ni
   (E,F)
\\
   \longmapsto E\boxtimes F
   \in K^{G_{W^1}\times
    G_{W^2}\times\C^*}(\M(\bv^1,\bw^1)\times\M(\bv^2,\bw^2))
\end{multline}
and a similar homomorphisms for $\La(\bv^1,\bw^1)$ and
$\La(\bv^2,\bw^2)$. We call them {\it K\"unneth homomorphisms}.

\begin{Theorem}\label{thm:Kunneth}
The two K\"unneth homomorphisms are isomorphisms. If $A$ is an abelian 
reductive subgroup of $G_{W^1}\times G_{W^2}\times\C^*$, the same
holds for the fixed point set
$\M(\bv^1,\bw^1)^A\times\M(\bv^2,\bw^2)^A$,
$\La(\bv^1,\bw^1)^A\times\La(\bv^2,\bw^2)^A$.
\end{Theorem}

\begin{proof}
In \cite[\S7]{Na-qaff}, the following was shown: varieties
$\M(\bv,\bw)$, $\La(\bv,\bw)$ and fixed point sets have
$\alpha$-partitions (see [loc.\ cit., 7.1] for definition) such that
each piece is an affine space bundle over a nonsingular projective
manifold, which is a locally equivariant vector bundle. Moreover, each
base manifold has a decomposable diagonal class as in [loc.\ cit., 7.2.1].

Note that K\"unneth homomorphisms are defined for arbitary varieties,
and have obvious functorial properties. By the arguments in [loc.\ 
cit., \S7.1], it is enough to show that the K\"unneth homomorphism is
an isomorphism for each base manifold.
By the proof of [loc.\ cit., 7.2.1], each base manifold has this property.
\end{proof}

Let us define subsets of $\M(\bw)$ by
\begin{gather*}
   \Zm \equiv \Zm(\bw^1;\bw^2)
   \defeq \left\{ [B,i,j]\in \M(\bw) \left| \;
   \text{$\lim_{t\to 0} \lambda(t)\ast[B,i,j]$ exists} \right\}\right.,
\\
   \Zl \equiv \Zl(\bw^1;\bw^2)
   \defeq \left\{ [B,i,j]\in \M(\bw)  \left| \;
   \lim_{t\to 0} \lambda(t)\ast[B,i,j]\in \La(\bw^1)\times\La(\bw^2)
   \right\}\right..
\end{gather*}
We use the symbol $\Zm(\bw^1;\bw^2)$, $\Zl(\bw^1;\bw^2)$ when we want
to emphasize the dimensions.
These subsets are invariant under the action of $G_{W^1}\times
G_{W^2} \times \C^*$.

Since $\pi\colon\M(\bw)\to \M_0(\bw)$ is a projective morphism, the
following is clear.
\begin{Lemma}\label{lem:saturated}
We have
\begin{gather*}
   \Zm 
   = \left\{ [B,i,j]\in \M(\bw) \left| \;
   \text{$\lim_{t\to 0} \lambda(t)\ast\pi([B,i,j])$ exists} \right\}\right.,
\\
   \Zl 
   = \left\{ [B,i,j]\in \M(\bw)  \left| \;
   \lim_{t\to 0} \lambda(t)\ast\pi([B,i,j]) = 0\right\}\right..
\end{gather*}
In particular, $\Zm$, $\Zl$ are $\pi$-saturated, i.e.,
$\pi^{-1}(\pi(\Zm)) = \Zm$, $\pi^{-1}(\pi(\Zl)) = \Zl$.
\end{Lemma}

\begin{Lemma}\label{lem:closed}
$\Zm$ and $\Zl$ are closed subvarieties of
$\M(\bw)$.
\end{Lemma}

\begin{proof}
It is enough to show that $\pi(\Zm)$ is a closed subvariety of
$\M_0(\infty,\bw)$. By \cite[1.3]{Lu-qv} the coordinate ring
$\M_0(\bv,\bw)$ (the ring of $G_V$-invariant polynomials in
$\mu^{-1}(0)\subset\bM(V,W)$) is generated by the following two
types of functions:
\begin{aenume}
\item $\tr(B_{h_N}B_{h_{N-1}}\cdots B_{h_1}\colon
V_{\vout(h_1)} \to V_{\vout(h_1)})$,
where $h_1$, \dots, $h_N$ is a cycle in our graph, i.e., $\vin(h_1) =
\vout(h_2)$, $\vin(h_2) = \vout(h_3)$, \dots, $\vin(h_{N-1}) =
\vout(h_N)$, $\vin(h_N) = \vout(h_1)$.
\item $\chi(j_{\vin(h_N)} B_{h_N} B_{h_{N-1}} \cdots B_{h_1}
i_{\vout(h_1)})$,
where $h_1$, \dots, $h_N$ is a path in our graph, i.e., $\vin(h_1) =
\vout(h_2)$, $\vin(h_2) = \vout(h_3)$, \dots, $\vin(h_{N-1}) =
\vout(h_N)$, and $\chi$ is a linear form on\linebreak[4]
$\Hom(W_{\vout(h_1)}, W_{\vin(h_N)})$.
\end{aenume}
Functions of the first type are invariant under the
$\lambda(\C^*)$-action. Functions of the second type are of weight
$1$, $-1$, $0$ if $\chi$ is the extension (by $0$) of a linear form of
$\Hom(W^1_{\vout(h_1)}, W^2_{\vin(h_N)})$,
$\Hom(W^2_{\vout(h_1)}, W^1_{\vin(h_N)})$,
$\Hom(W^1_{\vout(h_1)}, W^1_{\vin(h_N)})
\oplus \Hom(W^1_{\vout(h_1)}, W^2_{\vin(h_N)})$ respectively.

Thus $[B,i,j]\in\M_0(\bv,\bw)$ is contained in $\pi(\Zm)$ if and only
if $j_{\vin(h_N)} B_{h_N} B_{h_{N-1}} \cdots B_{h_1} i_{\vout(h_1)}$
maps $W_{\vout(h_1)}^2$ into $W_{\vin(h_N)}^2$ for any
path $h_1$, \dots, $h_N$.
Similarly $[B,i,j]\in\M_0(\bv,\bw)$ is contained in $\pi(\Zl)$ if and
only if functions of the first type vanishes, and $j_{\vin(h_N)}
B_{h_N} B_{h_{N-1}} \cdots B_{h_1} i_{\vout(h_1)}$ maps
$W_{\vout(h_1)}^2$ into $0$, and $W_{\vout(h_1)}^1$ to
$W_{\vin(h_N)}^2$ for any path $h_1$, \dots, $h_N$.
Now the assertions are clear from these descriptions.
\end{proof}

The limit $\lim_{t\to 0} \lambda(t)\ast[B,i,j]$ must be contained in
the fixed point set $\M(\bw)^{\lambda(\C^*)}$ if it exists. Thus we
have the decomposition
\begin{equation*}
\begin{split}
   & \Zm(\bw^1;\bw^2) = \bigsqcup_{\bv^1,\bv^2}
                     \Zm(\bv^1,\bw^1;\bv^2,\bw^2),
\\
   & \qquad\text{where }
   \Zm(\bv^1,\bw^1;\bv^2,\bw^2) \defeq
   \left\{ [B,i,j]  \left| \;
   \lim_{t\to 0} \lambda(t)\ast[B,i,j]\in
   \M(\bv^1,\bw^1)\times\M(\bv^2,\bw^2) \right\}\right..
\end{split}
\end{equation*}
Similarly we have
\(
   \Zl(\bw^1;\bw^2) = \bigcup_{\bv^1,\bv^2}
                     \Zl(\bv^1,\bw^1;\bv^2,\bw^2)
\)
defined exactly in the same way.
These are the Bialynicki-Birula decomposition of $\Zm$, $\Zl$. Thus we 
have
\begin{Proposition}\label{prop:BB}
\rom{(1)} $\Zm(\bv^1,\bw^1;\bv^2,\bw^2)$,
$\Zl(\bv^1,\bw^1;\bv^2,\bw^2)$ are nonsingular locally closed
subvarieties of $\M(\bv^1+\bv^2, \bw^1+\bw^2)$.

\rom{(2)} The map
\[
   \Zm(\bv^1,\bw^1;\bv^2,\bw^2)\ni [B,i,j]
   \longmapsto \lim_{t\to 0} \lambda(t)\ast [B,i,j]
   \in \M(\bv^1,\bw^1)\times\M(\bv^2,\bw^2)
\]
identifies $\Zm(\bv^1,\bw^1;\bv^2,\bw^2)$ with a fiber bundle over
$\M(\bv^1,\bw^1)\times\M(\bv^2,\bw^2)$, where the fiber over
$([B^1,i^1,j^1],[B^2,i^2,j^2])$ is isomorphic to the affine space
given by the direct sum of eigenspaces with positive weights in the
tangent space
\(
    T_{([B^1,i^1,j^1],[B^2,i^2,j^2])} \M(\bv^1+\bv^2,\bw^1+\bw^2)
\).
Similarly the map
\(
   \Zl(\bv^1,\bw^1;\bv^2,\bw^2) \to \La(\bv^1,\bw^1)\times\La(\bv^2,\bw^2)
\)
is a fiber bundle with the same fiber. Both fiber bundles are locally
$G_{W^1}\times G_{W^2}\times\C^*$-equivariant vector bundles.

\rom{(3)} There exists an ordering $<$ on the set
$\{ (\bv^1,\bv^2) \mid \bv^1 + \bv^2 = \bv \}$ \rom(the set of
components of $\M(\bv,\bw)^{\lambda(\C^*)}$\rom) such that
\[
   \bigcup_{(\bv^1,\bv^2)\le (\bv^1_0,\bv^2_0)}
   \Zm(\bv^1,\bw^1;\bv^2,\bw^2),
\qquad
   \left(\text{resp.\ }
   \bigcup_{(\bv^1,\bv^2)\le (\bv^1_0,\bv^2_0)}
   \Zl(\bv^1,\bw^1;\bv^2,\bw^2)\right)
\]
are closed subvarieties of $\Zm$ \rom(resp.\ $\Zl$\rom) for any
fixed $(\bv^1_0,\bv^2_0)$.
\end{Proposition}

See \cite[7.2.5]{Na-qaff}. (The assumption on the properness of the
moment map, which does not hold in the present case, was used to
ensure that the whole space $X$ is a union of $+$-attracting sets.)

\begin{Remark}\label{rem:order}
In fact, the order $<$ in (3) can be described explicitly.
If $f$ is the moment map, then
\[
 (\bv^1,\bv^2) < (\bv^{\prime 1}, \bv^{\prime 2})\Longrightarrow
 f(\M(\bv^1,\bw^1)\times\M(\bv^2,\bw^2)) >
 f(\M(\bv^{\prime 1},\bw^1)\times\M(\bv^{\prime 2},\bw^2)).
\]
With the K\"ahler metric in \cite{Na-quiver}, we have
\[
   f(\M(\bv^1,\bw^1)\times\M(\bv^2,\bw^2))
   = \sum_{k\in I} \zeta^{(k)}_\R \dim V^2_k,
\]
where $\zeta^{(k)}_\R\in \R_{>0}$ ($k\in I$) are parameter for the K\"ahler
metric. The ordering $<$ is independent of the metric, so we can move
the parameters. Therefore we may assume
\(
   (\bv^1,\bv^2) \le (\bv^{\prime 1}, \bv^{\prime 2})
\)
if and only if
\(
    \dim V^2_k \ge \dim V^{\prime 2}_k % \quad\text{for all $k\in I$}.
\)
for all $k\in I$. In particular, $(0,\bv)$ is a minimal element if
$\M(0,\bw^1)\times \M(\bv,\bw^2)$ is nonempty.
Thus $\Zm(0,\bw^1;\bv,\bw^2)$ and $\Zl(0,\bw^1;\bv,\bw^2)$ are
(nonsingular) closed subvarieties.
These are analogue of $\Zm_1$, $\Zl_1$ in \cite[7.10]{Lu-Base2}, and
will play an important role later.
\end{Remark}

\propref{prop:BB} has the following corollary.

\begin{Theorem}\label{thm:free}
\rom{(1)} Both $\Zm$ and $\Zl$ satisfy the properties $(S)$ and
$(T_{G_{W^1}\times G_{W^2}\times\C^*})$.  \rom(See
\cite[\S7.1]{Na-qaff} for the definition.\rom)

\rom{(2)} We have an exact sequence
\[
  0 \to K^{G_{W^1}\times G_{W^2}\times\C^*}(\Zm_{\le(\bv^1_0,\bv^2_0)})
  \to K^{G_{W^1}\times G_{W^2}\times\C^*}(\Zm)
  \to K^{G_{W^1}\times
  G_{W^2}\times\C^*}(\Zm\setminus\Zm_{\le(\bv^1_0,\bv^2_0)})
  \to 0,
\]
where 
\(
  \Zm_{\le (\bv^1_0,\bv^2_0)} \defeq
   \bigcup_{(\bv^1,\bv^2)\le (\bv^1_0,\bv^2_0)}
   \Zm(\bv^1,\bw^1;\bv^2,\bw^2)
\).
The same holds if we replace $\Zm$ by $\Zl$.

\rom{(3)} The direct image maps
\[
   K^{G_{W^1}\times G_{W^2}\times\C^*}(\La(\bw))
   \to
   K^{G_{W^1}\times G_{W^2}\times\C^*}(\Zl)
   \to
   K^{G_{W^1}\times G_{W^2}\times\C^*}(\Zm)
   \to
   K^{G_{W^1}\times G_{W^2}\times\C^*}(\M(\bw))
\]
\rom(induced by the the inclusions $\La(\bw)\subset \Zl\subset
\Zm\subset \M(\bw)$\rom) are injective.
All four modules are free of the same rank.
\end{Theorem}

\begin{proof}(Compare \cite[6.17]{Lu-Base2}.)
(1)(2) The assertion follows from \propref{prop:BB} and results in
\cite[\S7]{Na-qaff}.

(3) We replace the group $G_{W^1}\times G_{W^2}\times\C^*$ by its
maximal torus $H_{\bw^1}\times H_{\bw^2}\times\C^*$. If we tensor the
fraction field of $R(H_{\bw^1}\times H_{\bw^2}\times\C^*)$ to the
above homorphisms, it becomes isomorphisms by \cite{T-loc} since the
$H_{\bw^1}\times H_{\bw^2}\times\C^*$-fixed points are the same on all
four varieties.
Then the assertion for the torus follows from the freeness ((1) and
\cite[7.3.5]{Na-qaff}) of modules. Taking the Weyl group invariant
part, we get the assertion.
\end{proof}

Let $\Zm(\bw^2;\bw^1)$, $\Zl(\bw^2;\bw^1)$ denote the varieties
exactly as above except that the role of $W^1$ and $W^2$ are
exchanged. By the description in the proof of \lemref{lem:closed}, the
intersection $\Zm(\bw^2;\bw^1)\cap \Zl(\bw^1;\bw^2)$ is equal to
$\La(\bw)$. In particular, we can define a bilinear pairing by
\begin{multline}
\label{eq:pair}
   K^{G_{W^1}\times G_{W^2}\times\C^*}(\Zm(\bw^2;\bw^1))
   \times K^{G_{W^1}\times G_{W^2}\times\C^*}(\Zl(\bw^1;\bw^2))
   \ni (F,F')
\\
  \longmapsto p_*(F\otimes_{\M(\bw)}^L F')
   \in R(G_{W^1}\times G_{W^2}\times\C^*),
\end{multline}
where $p\colon \La(\bw)\to\text{point}$ is the projection to the point.
A pairing 
\begin{equation}
\label{eq:pair2}
   H_*(\Zm(\bw^2;\bw^1),\Z)\times H_*(\Zl(\bw^1;\bw^2),\Z)\to \Z
\end{equation}
can be defined in a similar way.
\begin{Theorem}
Pairings \eqref{eq:pair}, \eqref{eq:pair2} are nondegenerate.
\end{Theorem}

The proof is the same as that in \cite[\S7]{Na-qaff}.

We need more precise description of the projection
\(
   \Zm(\bv^1,\bw^1;\bv^2,\bw^2) \to
   \M(\bv^1,\bw^1)\times\M(\bv^2,\bw^2)
\)
later.
By the description of the tangent bundle in the proof of
\lemref{lem:fixed}, the fiber of the projection is isomorphic to \(
\Ker\beta^{21}/\Ima\alpha^{21} \), the cohomology of the complex
\eqref{eq:hecke_complex}.

Fix representatives $(B^1,i^1,j^1)$, $(B^2,i^2,j^2)$ of
a point in $\M(\bv^1,\bw^1)\times \M(\bv^2, \bw^2)$. We take
$(C,I,J)\in\Ker\beta^{21}$.
We define a data $(B,i,j)$ in $\bM(V,W)$ so that its components in
$\bM(V^1,W^1)$, $\bM(V^2,W^2)$,
\(
   \HomE(V^1,V^2) \oplus 
    \HomL(W^1, V^2) \oplus \HomL(V^1,W^2)
\)
are given by $(B^1,i^1,j^1)$, $(B^2,i^2,j^2)$, $(C,I,J)$ respectively.
Then it satisfies $\varepsilon B B + ij = 0$. We claim that $(B,i,j)$
is stable: Suppose that $S$ is contained in $\Ker j$ and invariant
under $B$. Under the projection $V^1\oplus V^2 \to V^1$, the subspace
$S$ define a subspace $S^1\subset V^1$ which is contained in $\Ker
j^1$ and invariant under $B^1$. By the stability condition for
$(B^1,i^1,j^1)$, we have $S^1 = 0$. Therefore $S\subset V^2$. Now the
stability condition for $(B^2,i^2,j^2)$ implies that $S = 0$. Thus we
have a morphism $\Ker\beta^{21} \to \mu^{-1}(0)^{\operatorname{s}}$.
If $(C,I,J)-(C',I',J')\in \Ima\alpha^{21}$, corresponding two data are
in the same $G_V$-orbit.
Thus we have a morphism 
\( 
   \Ker\beta^{21}/\Ima\alpha^{21}\to \Zm(\bv^1,\bw^2;\bv^2,\bw^2)
\).
(The left hand side is the total space of the vector bundle.)

\begin{Proposition}\label{prop:ext}
The morphism
\[
   \Ker\beta^{21}/\Ima\alpha^{21}\to \Zm(\bv^1,\bw^2;\bv^2,\bw^2)
\]
is an isomorphism.
\end{Proposition}

\begin{proof}
The morphism is equivariant under the $\lambda(\C^*)$-action.
Thus it is enough to check the assertion in a neighbourhood of
$\M(\bv^1,\bw^1)\times\M(\bv^2,\bw^2)$. The tangent bundles of these
varieties, restricted to $\M(\bv^1,\bw^1)\times\M(\bv^2,\bw^2)$, are
both given by
\[
   \left(\Ker\beta^{21}/\Ima\alpha^{21}\right)\oplus
   T\left(\M(\bv^1,\bw^1)\times\M(\bv^2,\bw^2)\right),
\]
and the differential of the morphism is the identity. Hence we have
the assertion.
\end{proof}

We also have an isomorphism
\[
   \Ker\beta^{21}/\Ima\alpha^{21}|_{\La(\bv^1,\bw^2)\times
      \La(\bv^2,\bw)}
   \to \Zl(\bv^1,\bw^2;\bv^2,\bw^2),
\]
where the left hand side is the total space of the restriction of the
vector bundle.

\begin{Proposition}\label{prop:lag}
\rom{(1)} The set of irreducible components of $\Zl$
is naturally identified with the sets of irreducible components of
$\La(\bw^1)\times \La(\bw^2)$.

\rom{(2)} $\Zl$ is a lagrangian subvariety of $\M(\bw)$. More
precisely,  $\Zl\cap \M(\bv,\bw)$ is a lagrangian subvariety of
$\M(\bv,\bw)$ for each $\bv$.
\end{Proposition}

\begin{proof}
Recall that $\Zl\cap \M(\bv,\bw)$ is a finite union
$\bigcup_{\bv^1+\bv^2 = \bv} \Zl(\bv^1,\bw^2;\bv^2,\bw^2)$. Hence both
assertions will follow if we show that $\Zl(\bv^1,\bw^2;\bv^2,\bw^2)$
is a lagrangian subvariety.

Let $\rho\colon \Zl(\bv^1,\bw^2;\bv^2,\bw^2)\to \La(\bv^1,\bw^1)\times
\La(\bv^2,\bw^2)$ be the projection in \propref{prop:BB}(2). By
\propref{prop:ext}, we have the following exact sequence of vector
bundles over $\Zl(\bv^1,\bw^2;\bv^2,\bw^2)$:
\begin{equation*}
   0 \to \rho^*\left(\Ker\beta^{21}/\Ima\alpha^{21}\right) \to
   T\Zl(\bv^1,\bw^2;\bv^2,\bw^2) \to 
   \rho^*\left( T\La(\bv^1,\bw^1)\oplus T\La(\bv^2,\bw^2)\right)
   \to 0.
\end{equation*}
(More precisely, we must restrict to the inverse image of the
nonsingular locus of $\La(\bv^1,\bw^1)\times\La(\bv^2,\bw^2)$.)
By the definition of $\Ker\beta^{21}/\Ima\alpha^{21}$,
we have
\[
   \omega\left(\rho^*\left(\Ker\beta^{21}/\Ima\alpha^{21}\right),
      T\Zl(\bv^1,\bw^2;\bv^2,\bw^2)\right) = 0,
\]
and the induced bilinear form on
$\rho^*\left( T\La(\bv^1,\bw^1)\oplus T\La(\bv^2,\bw^2)\right)$
coincides with one induced from the symplectic form on
$\M(\bv^1,\bw^1)\times\M(\bv^2,\bw^2)$.
Since $\La(\bv^1,\bw^1)$, $\La(\bv^2,\bw^2)$ are lagrangian, the
latter vanishes.

It is also clear that the dimension of $\Zl(\bv^1,\bw^2;\bv^2,\bw^2)$
is half of that of $\M(\bv,\bw)$. Just note that
$\Ker\beta^{21}/\Ima\alpha^{21}$ and $\Ker\beta^{21}/\Ima\alpha^{21}$
are dual to each other with respect to the symplectic form.
\end{proof}

\begin{Remark}
By \cite{KN}, $\M(\bv,\bw)$ can be identified with a framed moduli
space of holomorphic vector bundles on an ALE space, when the graph of
type $ADE$, where $\bv$, $\bw$ correspond to the Chern class and the
framing. We have the following geometric description:
\[
   \Zm(\bv^1,\bw^1;\bv^2,\bw^2)
   = \{ \text{an exact sequence $0\to E^2 \to E \to E^1 \to 0$} \},
\]
where $E^1$ (resp.\ $E^2$) has the Chern class and the framing
corresponding to $\bv^1$, $\bw^1$ (resp.\ $\bv^2$, $\bw^2$). Here the
exact sequence is suppose to respect the framing. The inclusion
$\Zm(\bv^1,\bw^1;\bv^2,\bw^2)\to \M(\bv,\bw)$ is
$(0\to E^2\to E\to E^1\to 0)\mapsto E$, and
the projection $\Zm(\bv^1,\bw^1;\bv^2,\bw^2)\to \M(\bv^1,\bw^1)\times
\M(\bv^2,\bw^2)$ is $(0\to E^2\to E\to E^1\to 0)\mapsto (E^1,E^2)$.
The fiber over $(E^1,E^2)$ is the extension group
\(
   \Ext^1(E^2,E^1)
\)
(for the framed vector bundle). This is nothing but
$\Ker\beta^{21}/\Ima\alpha^{21}$. The duality between
$\Ker\beta^{21}/\Ima\alpha^{21}$ and $\Ker\beta^{12}/\Ima\alpha^{12}$
is nothing but the Serre duality
\(
   \Ext^1(E^2,E^1)^* \cong \Ext^1(E^1,E^2)
\),
where the canonical bundle of the ALE space is trivial.
\end{Remark}

\section{Crystal structure}\label{sec:crystal}

Let $\Irr\Zl$ be the set of irreducible components of $\Zl$. It is a
disjoint union of irreducible components of $\Zl\cap\M(\bv,\bw)$ for
various $\bv$.

We define $\wt\colon\Irr\Zl\to P$ by setting
\begin{equation*}
   \wt(X) = \bw-\bv
\end{equation*}
if $X$ is an irreducible component of $\Zl\cap\M(\bv,\bw)$.

Let $X$ be an irreducible component of $\Zl$. Taking a generic element
$[B,i,j]\in X$, we define
\begin{equation*}
   \varepsilon_k(X) \defeq \dim\left(V_k/\Ima\tau_k\right).
\end{equation*}
It defines a function $\varepsilon_k\colon \Irr\Zl\to \Z_{\ge 0}$. We
define $\varphi_k\colon \Irr\Zl\to \Z_{\ge 0}$ by
\begin{equation*}
   \varphi_k(X) = \varepsilon_k(X)+\langle h_k,\wt(X)\rangle.
\end{equation*}
We set $\wt_k(X) = \langle h_k,\wt(X)\rangle$ as before. We have
\begin{equation*}
   \dim \left( \Ker\tau_k/\Ima\sigma_k\right) = \varphi_k(X).
\end{equation*}

Suppose $\varepsilon_k(X)>0$ and $\wt X = \bw - \bv$.
Then we have Grassmann bundles
\begin{gather*}
   p\colon \M_{k;\varepsilon_k(X)}(\bv,\bw) \to
    \M_{k;0}(\bv-\varepsilon_k(X)\alpha_k,\bw),
\\
   p'\colon \M_{k;\varepsilon_k(X)-1}(\bv-\alpha_k,\bw) \to
    \M_{k;0}(\bv-\varepsilon_k(X)\alpha_k,\bw).
\end{gather*}
We define an irreducible component of
$\Zl\cap\M(\bv-\alpha_k,\bw)$ by
\begin{equation*}
   X' \defeq
   \overline{p^{\prime-1}(p(X\cap\M_{k;\varepsilon_k(X)}(\bv,\bw))}.
\end{equation*}
In fact, 
\begin{enumerate}
\item Since $\pi$ factors through $p$,
$p^{\prime-1}(p(X\cap\M_{k;\varepsilon_k(X)}(\bv,\bw))$
is contained in $\Zl$ by \lemref{lem:saturated}.
\item We have $\dim p^{\prime-1}(p(X\cap\M_{k;\varepsilon_k(X)}(\bv,\bw))
= \frac12 \dim \M(\bv-\alpha_k,\bw)$ by \eqref{eq:dim_fiber}.
\end{enumerate}
Thus $X'$ is an irreducible component of
$\Zl\cap\M(\bv-\alpha_k,\bw)$.
We define an operator $\widetilde e_k\colon
\Irr\Zl\to\Irr\Zl\sqcup\{0\}$ by
\begin{equation}
   \widetilde e_k(X) \defeq
   \begin{cases}
      X' & \text{if $\varepsilon_k(X) > 0$}, \\
      0  & \text{otherwise}.
   \end{cases}
\end{equation}

Similarly if $\varphi_k(X)>0$, we consider a Grassmann bundle
\begin{gather*}
   p''\colon \M_{k;\varepsilon_k(X)+1}(\bv+\alpha_k,\bw) \to
    \M_{k;0}(\bv-\varepsilon_k(X)\alpha_k,\bw),
\end{gather*}
and define
\begin{equation*}
    X'' \defeq 
   \overline{p^{\prime\prime-1}(p(X\cap\M_{k;\varepsilon_k(X)}(\bv,\bw))}.
\end{equation*}
(Note that $Q_k(\bv-\varepsilon_k(X)\alpha_k,\bw)$ is a vector bundle
of rank $\wt_k(X)+2\varepsilon_k(X) =
\varepsilon_k(X)+\varphi_k(X)$. So the Grassmann bundle $p''$ is
nonempty by the assumption $\varphi_k(X)>0$.)
We define an operator $\widetilde f_k\colon
\Irr\Zl\to\Irr\Zl\sqcup\{0\}$ by
\begin{equation}
   \widetilde f_k(X) \defeq
   \begin{cases}
      X'' & \text{if $\varphi_k(X) > 0$}, \\
      0  & \text{otherwise}.
   \end{cases}
\end{equation}

Then the following is clear.
\begin{Proposition}
$\Irr\Zl$ together with $\wt$, $\varepsilon_k$, $\varphi_k$,
$\widetilde e_k$, $\widetilde f_k$ is a normal crystal.
\end{Proposition}

\begin{Lemma}\label{lem:gen}
$\Irr\Zl$ is generated by the subset $\{ X \mid
\text{$\varepsilon_k(X) = 0$ for all $k\in I$}\}$ as a crystal.
\end{Lemma}

\begin{proof}
If $\varepsilon_{k_1}(X) > 0$ for some $k_1$, then we have
$X = \widetilde f_{k_1} \widetilde e_{k_1}(X)$. If
$\varepsilon_{k_2}(\widetilde e_{k_1}(X)) > 0$, then 
$X = \widetilde f_{k_2} \widetilde f_{k_1} (\widetilde e_{k_2}
\widetilde e_{k_1} X)$.
We continue this procedure successively. Since the total sum $\sum_k
\dim V_k$ of dimensions decreases under this, it will eventually stop
at
$\widetilde e_{k_N} \cdots \widetilde e_{k_1} X$ with
$\varepsilon_k(\widetilde e_{k_N} \cdots \widetilde e_{k_1} X) = 0$
for all $k\in I$.
\end{proof}

Let $\Irr\La(\bw)$ be the set of irreducible components of
$\La(\bw)$. It is a disjoint union of the set of irreducible
components of $\La(\bv,\bw)$ for various $\bv$.
Since $\La(\bv,\bw)\subset\Zl\cap\M(\bv,\bw)$ and both
have dimension $=\frac12 \dim\M(\bv,\bw)$, we have an inclusion
\[
    \Irr\La(\bw) \subset \Irr\Zl.
\]
We restrict the maps $\wt$, $\varepsilon_k$, $\varphi_k$ to
$\Irr\La(\bw)$. Moreover, it is clear that if $X\in\Irr\La(\bw)$ and
$\widetilde e_k(X)\neq 0$ (resp.\ $\widetilde f_k(X)\neq 0$), then
$\widetilde e_k(X)\in\Irr\La(\bw)$ (resp. $\widetilde
f_k(X)\in\Irr\La(\bw)$). Thus $\Irr\La(\bw)$ inherits the structure
of the crystal from that of $\Irr\Zl$. In fact, this crystal structure 
on $\Irr\La(\bw)$ was introduced in \cite{Na-alg}, and our definition
here is a straightforward modification. The crystal structure in
[loc. cit.] was motivated by a similar construction by
Lusztig~\cite{Lu-crystal}.  

Note that $\M(0,\bw)$ consists of a single point. So $\La(0,\bw) =
\Zl(\bw^1;\bw^2)\cap\M(0,\bw) = \M(0,\bw)$. Let $[0]_\bw$ denote this
point considered as an irreducible component of $\La(\bw)$.

\begin{Proposition}\label{prop:sub}
$\Irr\La(\bw)$ is a highest weight normal crystal, which is the
strictly embedded crystal of $\Irr\Zl$ generated by $[0]_\bw$.
\end{Proposition}

\begin{proof}
By \lemref{lem:gen} $\Irr\La(\bw)$ is generated by
$\{ X\in \Irr\La(\bw)\mid \text{$\varepsilon_k(X) = 0$ for all
$k$}\}$. But this set consists of the single element $[0]_\bw$ as
shown in the proof of \cite[7.2]{Na-alg}.
\end{proof}

By \propref{prop:lag}, we have a natural bijection between sets
\begin{equation*}
   \Irr\La(\bw^1)\times\Irr\La(\bw^2) \cong
   \Irr\Zl.
\end{equation*}
Hereafter these two sets are identified and their element is denoted
by $X^1\otimes X^2$ for $X^1\in\Irr\La(\bw^1)$,
$X^2\in\Irr\La(\bw^2)$. We have
\begin{equation*}
   \wt(X^1\otimes X^2) = \bw - \bv
   = \bw^1 - \bv^1 + \bw^2 - \bv^2 = \wt(X^1) + \wt(X^2),
\end{equation*}
if $X^1$ (resp.\ $X^2$) is an irreducible component of
$\La(\bv^1,\bw^1)$ (resp.\ $\La(\bv^2,\bw^2)$). Thus we have
\eqref{eq:wt_tensor}.

The following is one of the main results in this article.
\begin{Theorem}\label{thm:main1}
The crystal $\Irr\Zl$ is isomorphic to
$\Irr\La(\bw^1)\otimes\Irr\La(\bw^2)$ as a crystal.
\end{Theorem}

Together with \propref{prop:sub}, this theorem implies that
$\{ \Irr\La(\bw) \mid \bw\in P^+\}$ is a closed family of highest
weight normal crystals.
By \propref{prop:closed} we have
\begin{Corollary}\label{cor:cryst}
$\Irr\La(\bw)$ is isomorphic to $\B(\bw)$ as a crystal for any $\bw\in
P^+$.
\end{Corollary}

The proof of \thmref{thm:main1} occupies the rest of this section. A
completely different proof of \corref{cor:cryst} will be given in
\secref{sec:comb}.

Let $X^1\otimes X^2\in \Irr\Zl$. Take a generic element $[B,i,j]$ 
in $X^1\otimes X^2$. By \propref{prop:ext} there exists a direct
sum decomposition $V = V^1\oplus V^2$ such that
\begin{aenume}
\item $V^2$ is invariant under $B$ and satisfies
$j(W^2)\subset V^2$, $i(V^2)\subset W^2$.
\item The data obtained by the restriction of $[B,i,j]$ to $V^2$,
$W^2$ (well-defined thanks to (a)) is contained in $X^2$.
\item The data induced on $V^1 = V/V^2$, $W^1 = W/W^2$ from $[B,i,j]$
(well-defined thanks to (a)) is contained in $X^1$.
\end{aenume}
In particular, the homomorphisms $\sigma_k$, $\tau_k$ decompose into
the following three parts:
\begin{alignat*}{2}
   & \sigma_k^1\colon V^1_k \to \bigoplus_{h:\vin(h)=k}
     V^1_{\vout(h)}\oplus W^1_k, & \qquad
   &  \tau_k^1 \colon \bigoplus_{h:\vin(h)=k} V^1_{\vout(h)}\oplus W^1_k
     \to V^1_k,
\\
   & \sigma_k^{21}\colon V^1_k \to \bigoplus_{h:\vin(h)=k}
     V^2_{\vout(h)}\oplus W^2_k, & \qquad
   &  \tau_k^{21} \colon \bigoplus_{h:\vin(h)=k} V^1_{\vout(h)}\oplus W^1_k
     \to V^2_k,
\\
   & \sigma_k^2\colon V^2_k \to \bigoplus_{h:\vin(h)=k}
     V^2_{\vout(h)}\oplus W^2_k, & \qquad
   &  \tau_k^2 \colon \bigoplus_{h:\vin(h)=k} V^2_{\vout(h)}\oplus W^2_k
   \to V^2_k.
\end{alignat*}
The equation $\tau_k\sigma_k = 0$ is equivalent to
$\tau_k^1\sigma_k^1 = 0$, $\tau_k^2\sigma_k^2 = 0$,
$\tau_k^{21}\sigma_k^1 + \tau_k^2\sigma_k^{21} = 0$. In particular,
$\tau_k^{21}$ induces a homomorphism
$\bigoplus_{h:\vin(h)=k} V^1_{\vout(h)}\oplus W^1_k / \Ima\sigma_k^1
\to V^2_k/\Ima\tau_k$.
Let $\overline{\tau^{21}_k}\colon \Ker\tau_k^1/\Ima\sigma_k^1 \to
V^2_k/\Ima\tau^2_k$ be its restriction.

The inclusion $V^2_k \subset V^1_k\oplus V^2_k$ induces a homomorphism 
\begin{equation*}
   \Ker\tau^2_k / \Ima\sigma^2_k \to
   \Ker\tau_k / \Ima \sigma_k.
\end{equation*}
The projection $V^1_k\oplus V^2_k\to V^1_k$ induces a homomorphism
\begin{equation*}
   \Ker\tau_k / \Ima\sigma_k
   \to \Ker\overline{\tau^{21}_k}.
\end{equation*}
Combining these, we have a three-term sequence of vector bundles:
\begin{equation}
\label{eq:complex}
   0 \to \Ker\tau^2_k / \Ima\sigma^2_k \to
      \Ker\tau_k / \Ima\sigma_k \to \Ker\overline{\tau^{21}_k}
      \to 0.
\end{equation}

\begin{Lemma}\label{lem:exact}
\eqref{eq:complex} is an exact sequence.
\rom(This holds any $[B,i,j]$, not necessarily generic.\rom)
\end{Lemma}

\begin{proof}
By the definition, it is clear that the composite of homomorphisms
is $0$. The exactness at the middle and the right terms are clear.
We now prove the exactness at the left term.

Suppose that $x\in \Ker\tau^2_k$, considered as an element of
$\left(\bigoplus_{h:\vin(h)=k} V^1_{\vout(h)} \oplus
  W^2_k\right)\oplus\Ker\tau^2_k$, is equal to $\sigma_k(y)$ for some
$y = y^1\oplus y^2\in V^1_k\oplus V^2_k$. Since we have
$\sigma^1_k(y^1) = 0$, the injectivity of $\sigma^1_k$ implies that
$y^1 = 0$. Then $\sigma_k(y) = x$ implies that $\sigma_k(y^2) =
x$. Thus we get the exactness.
\end{proof}

\begin{Lemma}\label{lem:ep=0}
Suppose $\varepsilon_k(X^1\otimes X^2) = 0$. Then

\rom{(1)} $\varepsilon_k(X^1) = 0$, $\dim\Ker\overline{\tau^{21}_k}
= \wt_k(X^1) - \varepsilon_k(X^2)$. In particular,
$\wt_k(X^1) \ge \varepsilon_k(X^2)$.

\rom{(2)} We have
\begin{equation*}
{\widetilde f_k}^r (X^1\otimes X^2)
=
\begin{cases}
   {\widetilde f_k}^r X^1 \otimes X^2
   & \text{if $r \le \wt_k(X^1)-\varepsilon_k(X^2)$}, \\
   {\widetilde f_k}^{\wt_k(X^1)-\varepsilon_k(X^2)} X^1
   \otimes {\widetilde f_k}^{r-\wt_k(X^1)+\varepsilon_k(X^2)}X^2
   & \text{otherwise}.
\end{cases}
\end{equation*}
\end{Lemma}

\begin{proof}
(1) By the assumption, $\tau_k$ is surjective. It is true if and only
if the following two statements hold:
\begin{aenume}
\item $\tau_k^1$ is surjective,
\item $\displaystyle
\Ker\tau_k^1\oplus \bigoplus_{h:\vin(h)=k} V^2_{\vout(h)} \oplus
W^2_k \xrightarrow{[\tau^{21}_k|_{\Ker\tau_k^1}\ \tau^2_k]} V^2_k$ is
surjective.
\end{aenume}
Then the second statement is
equivalent to
\begin{itemize}
\item[(b')] The homomorphism $\overline{\tau^{21}_k}$ is surjective.
\end{itemize}
By (a) we have
\(
   \varepsilon_k(X^1) = 0
\).
Therefore
\(
  \dim\Ker\tau^1_k/\Ima\sigma^1_k = \varphi_k(X^1) = \wt_k(X^1)
\).
By (b') we have
\(
   \dim\Ker\overline{\tau^{21}_k}
     = \wt_k(X^1) - \varepsilon_k(X^2).
\)

(2) A subspace $S$ of $\Ker\tau_k/\Ima\sigma_k$ defines subspaces
$S^1$ of $\Ker\overline{\tau^{21}_k}$ and $S^2$ of $\Ker\tau^2_k /
\Ima\sigma^2_k$ with an exact sequence
$0\to S^2 \to S \to S^1 \to 0$. For a generic $S$, we have
\begin{equation*}
\begin{split}
   & \dim S^1 = \min(\dim S, \dim\Ker\overline{\tau^{21}_k}),
\\
   & \dim S^2 = \max(0, \dim S - \dim\Ker\overline{\tau^{21}_k}).
\end{split}
\end{equation*}
Therefore we have the assertion.
\end{proof}

Now \thmref{thm:main1} follows from the following elementary lemma.
\begin{Lemma}
Let $\B_1$, $\B_2$ be normal crystals. Suppose that $\B_1\times\B_2$
has a structure of a normal crystal such that
\begin{aenume}
\item \eqref{eq:wt_tensor} holds for any $b_1\otimes b_2$.
\item If $\varepsilon_k(b_1\otimes b_2) = 0$, then we have
\eqref{eq:ep},\eqref{eq:e} and
\begin{equation*}
{\widetilde f_k}^r (b_1\otimes b_2)
=
\begin{cases}
   {\widetilde f_k}^r b_1 \otimes b_2
   & \text{if $r \le \wt_k(b_1)-\varepsilon_k(b_2)$}, \\
   {\widetilde f_k}^{\wt_k(b_1)-\varepsilon_k(b_2)} b_1
   \otimes {\widetilde f_k}^{r-\wt_k(b_1)+\varepsilon_k(b_2)}b_2
   & \text{otherwise}.
\end{cases}
\end{equation*}
\end{aenume}
\rom(We denote $(b_1,b_2)$ by $b_1\otimes b_2$.\rom)
Then the crystal structure coincides with that of the tensor product.
\end{Lemma}

\begin{proof}
Take $b_1\otimes b_2\in\B_1\times\B_2$ with $\varepsilon_k(b_1\otimes
b_2)=r$. Then we have
\(
  b_1\otimes b_2 = {\widetilde f_k}^r (b'_1\otimes b'_2)
\)
for some $b'_1\otimes b'_2\in\B_1\times\B_2$. We have
\(
  \varepsilon_k(b'_1\otimes b'_2) = 0
\).
Hence the condition~(b) implies
\begin{equation*}
\begin{split}
   & b_1 = 
\begin{cases}
   {\widetilde f_k}^r b'_1
   & \text{if $r \le \wt_k(b'_1)-\varepsilon_k(b'_2)$}, \\
   {\widetilde f_k}^{\wt_k(b'_1)-\varepsilon_k(b'_2)}
   b'_1 
   & \text{otherwise},
\end{cases}
\\
   & b_2 =
\begin{cases}
   b'_2
   & \text{if $r \le \wt_k(b'_1)-\varepsilon_k(b'_2)$}, \\
   {\widetilde f_k}^{r-\wt_k(b'_1)+\varepsilon_k(X^{\prime
   2})}b'_2
   & \text{otherwise}.
\end{cases}
\end{split}
\end{equation*}

First consider the case 
\(
   r < \wt_k(b'_1)-\varepsilon_k(b'_2)
\).
Then we have
\begin{equation}\label{eq:lll}
\begin{gathered}
   \varepsilon_k(b_1) = r, \quad
   \varepsilon_k(b_2) = \varepsilon_k(b'_2),
\\
   \widetilde e_k (b_1\otimes b_2)
   = {\widetilde f_k}^{r-1} (b'_1\otimes b'_2)
   = {\widetilde f_k}^{r-1} b'_1\otimes b'_2
   = \widetilde e_k b_1\otimes b_2,
\\
   \widetilde f_k (b_1\otimes b_2)
   = {\widetilde f_k}^{r+1} (b'_1\otimes b'_2)
   = {\widetilde f_k}^{r+1} b'_1\otimes b'_2
   = \widetilde f_k b_1\otimes b_2.
\end{gathered}
\end{equation}
The first two inequalities give us
\begin{equation}\label{eq:ineq}
   \varepsilon_k(b_2) - \wt_k(b_1) 
  = \varepsilon_k(b'_2) - \wt_k(b'_1) + 2r
 < r = \varepsilon_k(b_1),
\end{equation}
where we have used the assumption in the inequality. Thus we have
checked \eqref{eq:ep},\eqref{eq:e},\eqref{eq:f} in this case.

Next consider the case
\(
   r = \wt_k(b'_1)-\varepsilon_k(b'_2)
\).
We have \eqref{eq:lll} except that the last line is replaced by
\begin{equation*}
   \widetilde f_k (b_1\otimes b_2)
   = {\widetilde f_k}^{r+1} (b'_1\otimes b'_2)
   = {\widetilde f_k}^r b'_1\otimes \widetilde f_k b'_2
   = b_1\otimes \widetilde f_k b_2.
\end{equation*}
The inequality \eqref{eq:ineq} is replaced by
\begin{equation*}
   \varepsilon_k(b_2) - \wt_k(b_1)
   = \varepsilon_k(b_1).
\end{equation*}
We also have \eqref{eq:ep},\eqref{eq:e},\eqref{eq:f} in this case.

Finally consider the case
\(
   r > \wt_k(b'_1)-\varepsilon_k(b'_2)
\).
Then we have
\begin{equation*}
\begin{gathered}
 \varepsilon_k(b_1) 
 = \wt_k(b'_1) - \varepsilon_k(b'_2), \quad
 \varepsilon_k(b_2)
 = r - \wt_k(b'_1) + 2\varepsilon_k(b'_2),
\\
   \widetilde e_k (b_1\otimes b_2)
   = {\widetilde f_k}^{r-1} (b'_1\otimes b'_2)
   = b_1\otimes \widetilde e_k b_2,
\\
   \widetilde f_k (b_1\otimes b_2)
   = {\widetilde f_k}^{r+1} (b'_1\otimes b'_2)
   = b_1\otimes \widetilde f_k b_2.
\end{gathered}
\end{equation*}
We have
\begin{equation*}
\begin{split}
   & \varepsilon_k(b_2) - \wt_k(b_1)
   = r - \wt_k(b'_1) + 2\varepsilon_k(b'_2)
     - \wt_k(b'_1) + 2\left(\wt_k(b'_1) -
   \varepsilon_k(b'_2)\right)
\\
 = \; & r > \wt_k(b'_1)-\varepsilon_k(b'_2) = \varepsilon_k(b_1),
\end{split}
\end{equation*}
where we have used the assumption in the inequality.
Thus we have \eqref{eq:ep},\eqref{eq:e},\eqref{eq:f} in this case.

We have checked \eqref{eq:ep},\eqref{eq:e},\eqref{eq:f} in all cases,
and \eqref{eq:phi} follows from
\(
   \varphi_k(b_1\otimes b_2)
   = \varepsilon_k(b_1\otimes b_2) + \wt_k(b_1\otimes b_2)
\).
\end{proof}

\begin{Remark}
The rule \eqref{eq:ep} is equivalent to that $\overline{\tau_k^{21}}$
is full rank, i.e.,
\begin{equation*}
    \rank\overline{\tau_k^{21}}
    = \max\left(\dim(\Ker\tau^1_k/\Ima\sigma^1_k),
      \dim(V^2_k/\Ima\tau^2_k)\right).
\end{equation*}
\end{Remark}

\begin{Remark}
Saito \cite{Saito} proved \corref{cor:cryst} using the main result of
\cite{KS}. On the other hand, one can show the main result of
\cite{KS} from \corref{cor:cryst}. The detail is left for the reader.
\end{Remark}

\section{$\mathfrak g$-module structure}\label{sec:g-mod}

Let $Z(\bw) \defeq \M(\bw)\times_{\M_0(\infty,\bw)} \M(\bw)$. Let
$H_{\topdeg}(Z(\bw),\Q)$ be the top degree part of the Borel-Moore
homology group of $Z(\bw)$. More precisely, it is the subspace
\[
   \prod_{\bv,\bv'}'
   H_{\dim_\C \M(\bv',\bw)\times\M(\bv,\bw)}
   (Z(\bw)\cap\M(\bv',\bw)\times\M(\bv,\bw),\Q),
\]
of the direct products consisting elements $(F_{\bv,\bv'})$ such
that
\begin{enumerate}
\item for fixed $\bv$, $F_{\bv,\bv'} = 0$ for all but finitely
many choices of $\bv'$,
\item for fixed $\bv'$, $F_{\bv,\bv'} = 0$  for all but finitely
many choices of $\bv$.
\end{enumerate}

By the convolution product (see \cite[\S8]{Na-qaff}), it is an
associative algebra with $1 = \sum_\bv [\Delta(\bv,\bw)]$.
Let $\omega\colon \M(\bv',\bw)\times\M(\bv,\bw)\to
\M(\bv,\bw)\times\M(\bv',\bw)$ be the flip of the components.
By a main result of \cite[\S9]{Na-alg}, the assignment
\begin{gather*}
   P^*\ni h\longmapsto \sum_{\bv} \langle h,\bw-\bv\rangle [\Delta(\bv,\bw)]
\\
   e_k \longmapsto \sum_{\bv} [\Pa_k(\bv,\bw)], \quad
   f_k \longmapsto \sum_{\bv} \pm [\omega\left(\Pa_k(\bv',\bw)\right)]
\end{gather*}
defines an algebra homomorphism
\begin{equation*}
   {\mathbf U}(\mathfrak g) \to H_{\topdeg}(Z(\bw),\Q).
\end{equation*}
Here the sign $\pm$ can be explicitly given by $\bv$, $\bw$. But its
explicit form is not important for our purpose.
(In \cite{Na-alg} the direct sum was used, and ${\mathbf U}(\mathfrak
g)$ was replaced by the modified enveloping algebra.)

Let $H_{\topdeg}(\La(\bw),\Q)$, $H_{\topdeg}(\Zl,\Q)$ be the top
degree part of the Borel-Moore homology group of $\La(\bw)$, $\Zl$,
where `the top degree part' means, as above, the complex dimension
of $\M(\bv,\bw)$. (The sum is the usual direct sum.)
The fundamental classes of irreducible components of $\La(\bw)$
(resp.\ $\Zl$) give a basis of $H_{\topdeg}(\La(\bw),\Q)$ (resp.\ 
$H_{\topdeg}(\Zl,\Q)$). An irreducible component and its fundamental
class is identified hereafter.
Recall that the inclusion $\La(\bw)\subset\Zl$ induces an inclusion
$\Irr\La(\bw)\subset\Irr\Zl$. Hence we have an injection
\begin{equation}
\label{eq:inclu}
  H_{\topdeg}(\La(\bw),\Q) \to H_{\topdeg}(\Zl,\Q).
\end{equation}

Since $\La(\bw)$ and $\Zl$ are $\pi$-saturated
(\lemref{lem:saturated}), the convolution
makes $H_{\topdeg}(\La(\bw),\Q)$ and $H_{\topdeg}(\Zl,\Q)$ into
$H_{\topdeg}(Z(\bw),\Q)$-modules. Moreover, \eqref{eq:inclu} is a
morphism of $H_{\topdeg}(Z(\bw),\Q)$-modules.
Hence they can be considered as $\mathfrak g$-modules.
By \cite[10.2]{Na-alg}, $H_{\topdeg}(\La(\bw),\Q)$ is the simple
$\mathfrak g$-module with highest weight $\bw$. The highest weight
vector is the fundamental class $[0]_\bw$ of $\La(0,\bw) = \M(0,\bw) =
\operatorname{point}$.

\begin{Theorem}\label{thm:main2}
$H_{\topdeg}(\Zl,\Q)$ is isomorphic to $V(\bw^1)\otimes V(\bw^2)$ as a
$\mathfrak g$-module.
\end{Theorem}

\begin{proof}
By the argument in \cite[9.3]{Na-alg}, the module
$H_{\topdeg}(\Zl,\Q)$ is integrable.
Since the category of integrable $\mathfrak g$-modules is completely
reducible, it is enough to show that characters of
$H_{\topdeg}(\Zl,\Q)$, and $V(\bw^1)\otimes V(\bw^2)$ are
equal.
By the definition of $h$, we have
\begin{equation*}
   h\, X = \langle h, \wt X\rangle X \quad\text{for $X\in\Irr\Zl$}.
\end{equation*}
Namely the weight of $X$ as crystal is the same as the weight defined
by the $\mathfrak g$-module structure. Thus \thmref{thm:main1} and
\corref{cor:cryst} imply the equality of the character.

In fact, we do not need the full power of \thmref{thm:main1}. We only
need two things:
(1) there exists an isomorphism of sets
$\Irr\Zl\cong\Irr\La(\bw^1)\times\Irr\La(\bw^2)$ satisfying
$\wt(X^1\otimes X^2) = \wt(X^1)+\wt(X^2)$, and
(2) $H_{\topdeg}(\La(\bw^p),\C)$ is isomorphic to $V(\bw^p)$ ($p=1,2$).
Thus our proof follows from \propref{prop:lag}(1) and \cite{Na-alg}.
\end{proof}

This theorem is abstract. It is desirable to have a concrete
construction of the isomorphism.
For example, the isomorphism can be made so that the injective
homorphism \eqref{eq:inclu} is identified with a homomorphism
\(
   V(\bw^1+\bw^2) \to V(\bw^1)\otimes V(\bw^2),
\)
sending $b_{\bw^1+\bw^2}$ to $b_{\bw^1}\otimes b_{\bw^2}$, where
$b_\lambda$ is the highest weight vector of $V(\lambda)$.
However this condition does not characterize the isomorphism.
In the rest of this section, we study $H_{\topdeg}(\Zl,\Q)$ further
for this desire.

Recall that $\Zl(0,\bw^1;\bv,\bw^2)$ is a (possibly empty) closed
subvariety of $\Zl\cap\M(\bv,\bw)$ (\remref{rem:order}). By
\propref{prop:ext}, it is a vector bundle over
$\La(0,\bw^1)\times\La(\bv,\bw^2)\cong \La(\bv,\bw^2)$. Let
$\Zl_1(\bv) \defeq \Zl(0,\bw^1;\bv,\bw^2)$, $\Zm_1(\bv) \defeq
\Zm(0,\bw^1;\bv,\bw^2)$, $\Zm_1 \defeq \bigsqcup_\bv \Zm_1(\bv)$, and
$\Zl_1 \defeq \bigsqcup_\bv \Zl_1(\bv)$.

\begin{Lemma}\label{lem:intersect}
Let $\widehat\Pa^{(n)}_k(\bv,\bw)$ denote the intersection of
$\Pa^{(n)}_k(\bv,\bw)$ and $\M(\bv-n\alpha_k,\bw)\times\Zm_1(\bv)$
\rom(as submanifolds of $\M(\bv-n\alpha_k,\bw)\times\M(\bv,\bw)$\rom).
Let $p$ \rom(resp.\ $p'$\rom) denote the projection
$\Zm_1(\bv)\to\M(\bv,\bw)$
\rom(resp.\ $\Zm_1(\bv-n\alpha_k)\to\M(\bv-n\alpha_k,\bw)$\rom).

\rom{(1)} $\widehat\Pa^{(n)}_k(\bv,\bw)$ is contained in
$\Zm_1(\bv-n\alpha_k)\times\Zm_1(\bv)$

\rom{(2)} The restriction of $\id\times p\colon
\Zm_1(\bv-n\alpha_k)\times\Zm_1(\bv)\to
\Zm_1(\bv-n\alpha_k)\times\M(\bv,\bw^2)$
to $\widehat\Pa^{(n)}_k(\bv,\bw)$ gives an isomorphism
\[
    \widehat\Pa^{(n)}_k(\bv,\bw) \xrightarrow[\cong]{\id\times p}
    (p'\times\id)^{-1}(\Pa^{(n)}_k(\bv,\bw^2))
    \cong \HomL(W^1,V')|_{\Pa^{(n)}_k(\bv,\bw^2)},
\]
where $\Pa^{(n)}_k(\bv,\bw^2)$ is the Hecke correspondence in
$\M(\bv-n\alpha_k,\bw^2)\times\M(\bv,\bw^2)$, and $V'$ is the
pull-back of the tautological vector bundle of the first factor of
$\M(\bv-n\alpha_k,\bw)\times\M(\bv,\bw)$.

\rom{(3)} The intersection is transverse.
\end{Lemma}

\begin{proof}
By \propref{prop:ext}, $\Zm_1(\bv)$ is the total space of the vector
bundle $\HomL(W^1,V)$ over $\M(0,\bw^1)\times\M(\bv,\bw^2)$ (use
$V^1 = 0$, $V^2 = V$). Then
\begin{equation*}
   \widehat\Pa^{(n)}_k(\bv,\bw)
%   \Pa^{(n)}_k(\bv,\bw)\cap
%   \left(\M(\bv-n\alpha_k,\bw)\times\Zm_1(\bv)\right)
 \cong
 \HomL(W^1,V')|_{\Pa^{(n)}_k(\bv,\bw^2)},
\end{equation*}
where the right hand side is the total space of the restriction of
the vector bundle $\HomL(W^1,V')$.
Since $\Zm_1(\bv-n\alpha_k)$ is the total space of $\HomL(W^1,V')$,
the statements (1) and (2) are clear.

Let us describe the tangent bundles of submanifolds on the
intersection. From now, restrictions or pull-backs of vector bundles
are denoted by the same notation as original bundles.
The tangent bundle of $\M(\bv-n\alpha_k,\bw)\times\M(\bv,\bw)$ appears
in an exact sequence
\begin{multline*}
   0 \to \HomL(V',W^1)\oplus \HomL(W^1,V')\oplus
         \HomL(V,W^1) \oplus \HomL(W^1,V)
\\
   \to T\M(\bv-n\alpha_k,\bw)\oplus T\M(\bv,\bw)
   \to T\M(\bv-n\alpha_k,\bw^2)\oplus T\M(\bv,\bw^2) \to 0.
\end{multline*}
The tangent bundle of $\M(\bv-n\alpha_k,\bw)\times\Zm_1(\bv)$ appears
as
\begin{multline*}
   0 \to \HomL(V',W^1)\oplus \HomL(W^1,V')\oplus \HomL(W^1,V)
\\
   \to T\M(\bv-n\alpha_k,\bw)\oplus T\Zm_1(\bv)
   \to T\M(\bv-n\alpha_k,\bw^2)\oplus T\M(\bv,\bw^2) \to 0.
\end{multline*}
The tangent bundle of $\Pa_k^{(n)}(\bv,\bw)$ appears as
\begin{equation*}
   0 \to \HomL(W^1,V')\oplus \HomL(V,W^1)
   \to T\Pa_k^{(n)}(\bv,\bw)
   \to T\Pa_k^{(n)}(\bv,\bw^2) \to 0.
\end{equation*}
Now the transversality is clear.
\end{proof}

Since $\Zl_1$ is a vector bundle over $\La(\bw^2)$, we have the Thom
isomorphism
\begin{equation}
\label{eq:Thom}
   H_{\topdeg}(\La(\bw^2),\Q) \cong H_{\topdeg}(\Zl_1,\Q).
\end{equation}
Under this isomorphism $\Irr\La(\bw^2)\ni X$ is mapped to
$[0]_{\bw^1}\otimes X\in\Irr\Zl$.

\begin{Proposition}\label{prop:qqq}
\rom{(1)} As a $\mathfrak g$-module, $H_{\topdeg}(\Zl,\Q)$ is
generated by its subspace $H_{\topdeg}(\Zl_1,\Q)$.

\rom{(2)} The subspace $H_{\topdeg}(\Zl_1,\Q)$ is invariant under the
action of ${\mathbf U}(\mathfrak g)^+$, and the Thom
isomorphism~\eqref{eq:Thom} is compatible with the ${\mathbf
U}(\mathfrak g)^+$-module structures.
\end{Proposition}

\begin{proof}
(1) The proof of \cite[10.2]{Na-alg} shows that $H_{\topdeg}(\Zl,\Q)$
is generated by elements $X^1\otimes X^2\in \Irr\Zl$ with
$\varepsilon_k(X^1\otimes X^2) = 0$ for all $k\in I$. By
\lemref{lem:ep=0}(1), we have $\varepsilon_k(X^1) = 0$ for all $k\in
I$. As we already used in the proof of \propref{prop:sub}, this
implies that $X^1$ is the highest weight vector $[0]_{\bw^1}$. Hence
$X^1\otimes X^2$ is contained in $H_{\topdeg}(\Zl_1,\Q)$.

(2) The first statement follows from \lemref{lem:intersect}(1), i.e.,
$\Pa_k(\bw)\cap \left(\M(\bw)\times \Zm_1\right) = \widehat\Pa_k(\bw)
\subset \Zm_1\times \Zm_1$.
The rest of proof is based on the arguement in
\cite[\S8]{Na-qaff}. Let $i\colon \Zm_1\to \M(\bw)$ be the
inclusion. By the pull-back with support map with respect to
$\id\times i\colon \M(\bw)\times \Zm_1\to \M(\bw)\times\M(\bw)$
(\cite[6.4]{Na-qaff}), we have
%\begin{multline*}
%   H_{\topdeg}(\Pa_k(\bw),\Q) \cong
%   H^{\topdeg}(\M(\bw)\times\M(\bw), \Pa_k(\bw),\Q)
%\\
%   \xrightarrow{(\id\times i)^*}
%   H^{\topdeg}(\M(\bw)\times\Zm_1, \widehat\Pa_k(\bw),\Q)
%   \cong
%   H_{\topdeg}(\widehat\Pa_k(\bw),\Q).
%\end{multline*}
\begin{equation*}
   H_{\topdeg}(\widehat\Pa_k(\bw),\Q)
   \xrightarrow{(\id\times i)^*}
   H_{\topdeg}(\widehat\Pa_k(\bw),\Q).
\end{equation*}
By \cite[8.2.3]{Na-qaff} this map is compatible with the convolution
product, that is, the following is a commutative diagram:
\begin{equation*}
\begin{CD}
  H_{\topdeg}(\Pa_k(\bw),\Q) \otimes H_{\topdeg}(\Zl_1,\Q)
      @>>> H_{\topdeg}(\Zl_1,\Q)
\\
  @V{(\id\times i)^*\otimes \id}VV @|
\\
  H_{\topdeg}(\widehat\Pa_k(\bw),\Q) \otimes H_{\topdeg}(\Zl_1,\Q)
      @>>> H_{\topdeg}(\Zl_1,\Q),
\end{CD}
\end{equation*}
where the upper horizontal arrow is the convolution relative to
$\M(\bw)$, and the lower horizontal arrow is the convolution relative
to $\Zm_1$. Furthermore, \lemref{lem:intersect}(3) implies that
$(\id\times i)^*[\Pa_k(\bw)] = [\widehat\Pa_k(\bw)]$.

By \lemref{lem:intersect}(2), we have an isomorphism
\begin{equation*}
    H_{\topdeg}(\widehat\Pa_k(\bw),\Q)
    \xrightarrow[\cong]{\left((p'\times\id)^*\right)^{-1}(\id\times p)_*}
    H_{\topdeg}(\Pa_k(\bw^2),\Q),
\end{equation*}
where $(p'\times\id)^*$ is the Thom isomorphism. By
\cite[8.3.5]{Na-qaff}, it is compatible with the convolution product,
that is, the following is commutative:
\begin{equation*}
\begin{CD}
  H_{\topdeg}(\widehat\Pa_k(\bw),\Q) \otimes H_{\topdeg}(\Zl_1,\Q)
      @>>> H_{\topdeg}(\Zl_1,\Q)
\\
  @V{\left((p'\times\id)^*\right)^{-1}(\id\times p)_*\otimes 
  (p^*)^{-1}}VV
  @VV{(p^*)^{-1}}V
\\
  H_{\topdeg}(\Pa_k(\bw^2),\Q) \otimes H_{\topdeg}(\La(\bw^2),\Q)
      @>>> H_{\topdeg}(\La(\bw^2),\Q),
\end{CD}
\end{equation*}
where the horizotal arrows are the convolution product relative to
$\Zm_1$ and $\M(\bw^2)$, and $p^*$ is the Thom isomorphism \eqref{eq:Thom}.
Moreover we have
\(
   \left((p'\times\id)^*\right)^{-1}(\id\times p)_*
   [\widehat\Pa_k(\bw)] = [\Pa_k(\bw^2)]
\)
by \lemref{lem:intersect}(2). Combining two commutative diagrams, we
get the assertion.
\end{proof}

\begin{Conjecture}
There exists a unique $\mathfrak g$-module isomorphism
\[
   H_{\topdeg}(\La(\bw^1),\Q)\otimes H_{\topdeg}(\La(\bw^2),\Q)
   \to H_{\topdeg}(\Zl,\Q)
\]
such that its restriction to
$[0]_{\bw^1}\otimes H_{\topdeg}(\La(\bw^2),\Q)$ is \eqref{eq:Thom},
composed with the inclusion $H_{\topdeg}(\Zl_1,\Q)\to
H_{\topdeg}(\Zl,\Q)$.
(The uniqueness is clear from \propref{prop:qqq}(1).)
\end{Conjecture}

The rest of this section is devoted to the proof of the conjecture
when $\mathfrak g$ is of type $ADE$.
When $\mathfrak g$ is of type $ADE$, $V(\bw^2)$ contains the lowest
weight vector $q_{{\mathbf m}^2}$, where $\mathbf m^2$ is the weight
of $q_{{\mathbf m}^2}$. We denote by ${}_{\mathbf m^2}[q]$ the
corresponding element in $\Irr\La(\bw^2)\subset
H_{\topdeg}(\La(\bw^2),\Q)$.
For a later purpose, we denote by $[0]_{\bw^1}\circ{}_{\mathbf
m^2}[q]$ the irreducible component of $\Zl$, corresponding to the
tensor product of $[0]_{\bw^1}$ and ${}_{\mathbf m^2}[q]$ under the
identification $\Irr\Zl = \Irr\La(\bw^1)\otimes\Irr\La(\bw^2)$.

\begin{Lemma}\label{lem:def_tensor}
For each $k\in I$, we have
\begin{subequations}
\begin{align}
  & e_k^{-\langle h_k, \mathbf m^2\rangle + 1}
    \left([0]_{\bw^1}\circ{}_{\mathbf m^2}[q]\right) = 0,
\\
  & f_k^{\langle h_k,\bw^1\rangle + 1}
    \left([0]_{\bw^1}\circ{}_{\mathbf m^2}[q]\right)
  = 0, \label{eq:mid}
\\
  & h\left([0]_{\bw^1}\circ{}_{\mathbf m^2}[q]\right) =
  \langle h,\bw^1+\mathbf m^2\rangle[0]_{\bw^1}\circ{}_{\mathbf m^2}[q].
\end{align}
\end{subequations}
\end{Lemma}

\begin{proof}
The first equation follows from \propref{prop:qqq}(2)
and
\(
   e_k^{-\langle h_k, \mathbf m^2\rangle + 1}\; {}_{\mathbf m^2}[q] = 0
\)
(a well-known property of the lowest weight vector).
The last equation follows from \thmref{thm:main1} and the
compatibility of weight structures for the crystal and the $\mathfrak
g$-module.

Let us show \eqref{eq:mid}. Assume that the left hand side is nonzero.
Then $\Pa^{(\langle h_k,\bw^1\rangle+1)}_k(\bw)$ intersects with
$p_2^{-1}([0]_{\bw^1}\circ{}_{\mathbf m^2}[q])$, where $p_2\colon
\M(\bw)\times\M(\bw)\to \M(\bw)$ is the second projection.  Take a
point $[B,i,j]$ in the image under $p_2$ of a point in the
intersection and consider the exact sequence \eqref{eq:complex}.  Then
$p_2^{-1}([B,i,j])\cap\Pa^{(\langle h_k,\bw^1\rangle+1)}_k(\bw)$ is
identified with the Grassmann variety of subspaces $S$ in
$\Ker\tau_k/\Ima\sigma_k$ with $\dim S = \langle h_k,\bw^1\rangle+1$
by \cite[4.5]{Na-alg} or \cite[5.4.3]{Na-qaff}.

Since $\mathbf m^2$ is the lowest weight, we have 
\(
   \widetilde f_k\; {}_{\mathbf m^2}[q] = 0
\).
Therefore,
\begin{equation*}
   0 = \varphi_k({}_{\mathbf m^2}[q]) 
     = \dim \left(\Ker\tau^2_k/\Ima\sigma^2_k\right).
\end{equation*}
We must be careful with the second equality since it holds only for a
{\it generic\/} $[B,i,j]$ in general. However, $\M(\bw^2-\mathbf
m^2,\bw^2) = \La(\bw^2-\mathbf m^2,\bw^2) = {}_{\mathbf m^2}[q]$ is
isomorphic to $\M(0,\bw^2)$, which is a single point, as seen by the
Weyl group symmetry \cite{Na-refl}.
So the equality holds for our $[B,i,j]$.
Thus we have $\Ker\tau_k/\Ima\sigma_k \cong\Ker\overline{\tau^{21}_k}$ 
by the exactness of \eqref{eq:complex}. But
\begin{equation*}
   \dim\Ker\overline{\tau^{21}_k}
   \le \dim \left(\Ker\tau_k^1/\Ima\sigma_k^1\right)
   = \dim W^1_k = \langle h_k,\bw^1\rangle.
\end{equation*}
So the Grassmann variety is empty. This contradiction comes from our
assumption. Hence we have \eqref{eq:mid}.
\end{proof}

\begin{Theorem}\label{thm:main3}
There exists a unique $\mathfrak g$-module isomorphism
\(
  H_{\topdeg}(\La(\bw^1),\Q)\otimes H_{\topdeg}(\La(\bw^2),\Q)\to
  H_{\topdeg}(\Zl,\Q)
\)
sending $[0]_{\bw^1}\otimes {}_{{\mathbf m}^2}[q]$ to
$[0]_{\bw^1}\circ{}_{\mathbf m^2}[q]$.
Moreover, its restriction to $[0]_{\bw^1}\otimes
H_{\topdeg}(\La(\bw^2),\Q)$ is \eqref{eq:Thom}, composed with the
inclusion $H_{\topdeg}(\Zl_1,\Q)\to H_{\topdeg}(\Zl,\Q)$.
\end{Theorem}

\begin{proof}
We identify $H_{\topdeg}(\La(\bw^1),\Q)\otimes
H_{\topdeg}(\La(\bw^2),\Q)$ with the tensor product $V(\bw^1)\otimes
V(\bw^2)$.

By \cite[23.3.6]{Lu-book} the assignment
\(
   \mathbf U(\mathfrak g)\ni u \mapsto
   u([0]_{\bw^1}\otimes {}_{{\mathbf m}^2}[q])
\)
is a surjective homomorphism
\(
  \mathbf U(\mathfrak g)\to H_{\topdeg}(\La(\bw^1),\Q)\otimes
  H_{\topdeg}(\La(\bw^2),\Q)
\)
with kernel
\[
  \sum_k \mathbf U(\mathfrak g) f_k^{w_k^1+1}
  +
  \sum_k \mathbf U(\mathfrak g) e_k^{-m_k^2+1}.
\]
Therefore \lemref{lem:def_tensor} implies that
\(
   \mathbf U(\mathfrak g)\ni u \mapsto
   u([0]_{\bw^1}\circ {}_{\mathbf m^2}[q])
\)
factors through a $\mathfrak g$-module homomorphism
\(
   H_{\topdeg}(\La(\bw^1),\Q)\otimes
  H_{\topdeg}(\La(\bw^2),\Q) \to H_{\topdeg}(\Zl,\Q)
\)
sending $[0]_{\bw^1}\otimes {}_{{\mathbf
m}^2}[q]$ to $[0]_{\bw^1}\circ{}_{\mathbf m^2}[q]$. The uniqueness is
clear.

Since $H_{\topdeg}(\La(\bw^2),\Q) = \mathbf U(\mathfrak
g)^+q_{{\mathbf m}^2}$, \propref{prop:qqq}(2) implies
\[
   H_{\topdeg}(\Zl_1,\Q) = \mathbf U(\mathfrak g)^+
   \left([0]_{\bw^1}\circ {}_{\mathbf m^2}[q]\right).
\]
Together with \propref{prop:qqq}(1), we have the surjectivity of
the homomorphism.

Now we compare the dimensions of the domain and the target. We have
\begin{equation*}
\begin{split}
   & \dim H_{\topdeg}(\La(\bw^1),\Q)\otimes
      H_{\topdeg}(\La(\bw^2),\Q)
   = \# \Irr\La(\bw^1) \# \Irr\La(\bw^2)
\\
  = \; & \# \Irr \Zl = \dim H_{\topdeg}(\Zl,\Q).
\end{split}
\end{equation*}
Thus the surjective homomorphism must be an isomorphism.

The second statement follows from \propref{prop:qqq}(2).
\end{proof}

\section{$\Ul$-module structure}\label{sec:Ul}

In this section, we assume $\mathfrak g$ is of type $ADE$.

For $p=1,2$, let $H_{W^p}$ be a maximal torus of $G_{W^p}$.  Let
$H = H_{W^1}\times H_{W^2}$, $\wH_{W^p} = H_{H^p}\times\C^*$,
$\wH = H\times\C^*$.
The representation ring $R(\C^*)$ of $\C^*$ is $\Z[q,q^{-1}]$, where
$q$ is the class of the canonical $1$-dimensional representation of
$\C^*$.
Considering $H$ as a subgroup of $G_W$, we make $\wH$ act on
$\M(\bw)$. The action preserves $\M(\bw^1)\times \M(\bw^2)$, $\Zm$ and
$\Zl$.
Let $K^{\wH}(Z(\bw))$ be the equivariant $K$-homology group of
$Z(\bw)$. It is an associative $R(\wH) \cong R(H)[q,q^{-1}]$-algebra
with unit under the convolution product.
Let $\Delta$ denote the diagonal embedding $\M(\bv,\bw)\to
\M(\bv,\bw)\times\M(\bv,\bw)$.
Let $\iota\colon\Pa_k(\bv,\bw)\to
Z(\bw)\cap\M(\bv',\bw)\times\M(\bv,\bw)$ denote the inclusion.
By \cite{Na-qaff}, the assignment
\begin{gather*}
   q^{h} \longmapsto \sum_\bv q^{\langle h, \bw-\bv\rangle}
     \Delta_* \shfO_{\M(\bv,\bw)},
\\
  \psi_k^\pm(z)
  \longmapsto
  \sum_{\bv} q^{\rank C_k^\bullet(\bv,\bw)}\Delta_*
      \left(\frac{\Wedge_{-1/qz} (C_k^\bullet(\bv,\bw))}
                 {\Wedge_{-q/z} (C_k^\bullet(\bv,\bw))}\right)^\pm,
\\ 
  e_{k,r} \longmapsto 
  \sum_{\bv} \pm \iota_*\left(q^{-1}V/V'\right)^{\otimes r-s}
  \otimes\mathcal L,
\\
  f_{k,r} \longmapsto 
  \sum_{\bv} \pm' \omega_*\iota_*\left(q^{-1}V/V'\right)^{\otimes r-s'}
  \otimes\mathcal L'
\end{gather*}
defines an algebra homomorphism
\begin{equation*}
   \Uli\otimes_{\Z[q,q^{-1}]} R(\wH) \to
   K^{\wH}(Z(\bw))/\operatorname{torsion}.
\end{equation*}
Here $\pm$, $\pm'\in \{1,-1\}$, $s$, $s'\in\Z$ can be given in terms
of $\bv$, $\bw$, $k$ as \cite[9.3.2]{Na-qaff}, but their explicit
forms are not important for our later discussion. Similarly $\mathcal L$,
$\mathcal L'$ are line bundles whose explicit forms are not important.
Those terms are independent of $r$.

The convolution make $K^{\wH}(\La(\bw))$, $K^{\wH}(\Zl)$,
$K^{\wH}(\Zm)$, $K^{\wH}(\M(\bw))$ into
$K^{\wH}(Z(\bw))/\operatorname{torsion}$-modules. (Note all these are
free.) Moreover, the inclusions in \thmref{thm:free}(3) respect the
module structure.
Let $[0]_\bw$ be the class represented by the structure sheaf of
$\M(0,\bw) = \La(0,\bw) = \text{point}$. By
\cite[12.3.2, 13.3.1]{Na-qaff}, $[0]_\bw$ is an {\it l\/}--highest
weight vector with Drinfeld polynomial $P_k(u) = \Wedge_{-u} q^{-1}
W_k$, that is,
\begin{gather*}
   e_{k,r}[0]_\bw = 0\quad\text{for any $k\in I$, $r\in \Z$},
\\
   K^{\wH}(\La(\bw)) = 
   \left(\Uli^-\otimes_{\Z[q,q^{-1}]} R(\wH)\right)[0]_{\bw},
\\
   \psi_k^\pm(z)[0]_\bw
   = \left( \Wedge_{-1/qz} (q^{-1} - q) W_k\right)^\pm
   [0]_\bw.
\end{gather*}
This $K^{\wH}(\La(\bw))$, more precisely, its Weyl group invariant
part $K^{G_W\times \C^*}(\La(\bw))$, is the universal standard
module $M(\bw)$ mentioned in the introduction and
\subsecref{subsec:qloop}.

Let us take a closed subvariety $\Zl_1$ be as in
\secref{sec:g-mod}. We have the Thom isomorphism in the $K$-theory:
\begin{equation}
\label{eq:Thom2}
   K^{\wH}(\La(\bw^2)) \cong K^{\wH}(\Zl_1),
\end{equation}
where $H_{W^1}$ acts trivially on $\La(\bw^2)$.
By \thmref{thm:free}(2) the inclusion $\Zl_1\subset \Zl$ induces an
injective $R(\wH)$-homomorphism
\begin{equation}
\label{eq:Zl1_incl}
   K^{\wH}(\Zl_1) \hookrightarrow K^{\wH}(\Zl).
\end{equation}
The following can be proved exactly as in \propref{prop:qqq}.
(For (1), we use \cite[12.3.2]{Na-qaff} instead of \cite[10.2]{Na-alg}.
And for (2), we must use $\Pa_k^{(n)}(\bw)$ corresponding to divided powers.)

\begin{Proposition}\label{prop:qqq2}
\rom{(1)} As a $\Uli\otimes_{\Z[q,q^{-1}]}R(\wH)$-module,
$K^{\wH}(\Zl)$ is generated by $K^{\wH}(\Zl_1)$.

\rom{(2)} $K^{\wH}(\Zl_1)$ is invariant under the
action of $\Uli^+$, and the Thom
isomorphism~\eqref{eq:Thom2} is compatible with the $\Uli^+$-module
structures \rom(up to shift $e_{k,r}\to e_{k,r+s}$ and invertible
elements in $R(\wH)$\rom).
\end{Proposition}

Let $\mathfrak R(\wH)$ be the fraction field of the $R(\wH)$. If $M$
is an $R(\wH)$-module, $M\otimes_{R(\wH)}\mathfrak R(\wH)$ is denoted
by $M_{\mathfrak R}$.
By the localization theorem~\cite{T-loc}, we have
\begin{equation*}
   K^{\wH}(\ast)_{\mathfrak R}
   \xrightarrow[\cong]{i^*}
   K^{\wH}(\ast^{\wH})_{\mathfrak R},
\end{equation*}
where $\ast = Z(\bw),\La(\bw),\Zl,\Zm$, or $\M(\bw)$, $\ast^{\wH}$
denotes the fixed point set, and $i$ denotes the inclusion
$\M(\bw)^{\wH}\times\M(\bw)^{\wH}\to \M(\bw)$ (if $\ast = Z(\bw)$) or
$\M(\bw)^{\wH}\to \M(\bw)$ (if $\ast\neq Z(\bw)$).
For $Z(\bw)$, we replace $i^*$ by $r = 1\boxtimes (\Wedge_{-1}
N^*)^{-1} i^*$, where $N$ is the normal bundle of $\M(\bw)^{\wH}$ in
$\M(\bw)$ as in \cite[5.11]{Gi-book}. Then it is compatible with the
convolution product, i.e., $r$ is an algebra homomorphism, and the
following is commutative:
\begin{equation*}
\begin{CD}
   K^{\wH}(Z(\bw))_{\mathfrak R}\otimes_{\mathfrak R(\wH)}
     K^{\wH}(\ast)_{\mathfrak R}
    @>>>
   K^{\wH}(\ast)_{\mathfrak R}
\\
   @V{r\otimes i^*}V{\cong}V @V{\cong}V{i^*}V
\\
   K^{\wH}(Z(\bw)^{\wH})_{\mathfrak R}\otimes_{\mathfrak R(\wH)}
    K^{\wH}(\ast^{\wH})_{\mathfrak R}
    @>>> 
   K^{\wH}(\ast^{\wH})_{\mathfrak R},
\end{CD}
\end{equation*}
where $\ast = \La(\bw),\Zl,\Zm$, or $\M(\bw)$, and the horizontal
arrows are convolution relative to $\M(\bw)$ and $\M(\bw)^{\wH}$
respectively. Therefore $K^{\wH}(\ast^{\wH})_{\mathfrak R}$ has a
structure of a $\Ul\otimes_{\Q(q)}\mathfrak R(\wH)$-module.

On the other hand, $K^{\wH}(\La(\bw^1))_{\mathfrak R}\otimes
K^{\wH}(\La(\bw^2))_{\mathfrak R}$ can be considered as a
$\Ul\otimes_{\Q(q)}\mathfrak R(\wH)$-module by the
comultiplication~\eqref{eq:comul}.
The following lemma first appeared in \cite[14.1.2]{Na-qaff}.

\begin{Lemma}\label{lem:general_point}
There exists a unique $\Ul\otimes_{\Q(q)}\mathfrak R(\wH)$-module
isomorphism
\begin{equation*}
   \Phi_{\mathfrak R} \colon
   K^{\wH}(\La(\bw^1))_{\mathfrak R}\otimes_{\mathfrak R(\wH)}
   K^{\wH}(\La(\bw^2))_{\mathfrak R}
     \to
   K^{\wH}(\Zl)_{\mathfrak R},
\end{equation*}
sending $[0]_{\bw^1}\otimes [0]_{\bw^2}$ to $[0]_\bw$.
\end{Lemma}

\begin{proof}
We have (cf.\ \cite[4.2.2]{Na-qaff})
\begin{equation}\label{eq:fixed}
%\begin{gathered}
   \La(\bw)^{\wH} = \Zl^{\wH} = \Zm^{\wH} =
   \M(\bw)^{\wH} 
   = \M(\bw^1)^{\wH}\times \M(\bw^2)^{\wH}
   = \La(\bw^1)^{\wH}\times \La(\bw^2)^{\wH}
.
%,
%\\
%   Z(\bw)^{\wH} 
%    = Z(\bw^1)^{\wH}\times Z(\bw^2)^{\wH}
%    = \left(\M(\bw^1)^{\wH}\right)^2\times
%    \left(\M(\bw^2)^{\wH}\right)^2.
%\end{gathered}
\end{equation}
By the argument in \cite[14.1.2]{Na-qaff}, this implies that
\(
%  K^{\wH}(\La(\bw)^{\wH})_{\mathfrak R} = 
%  K^{\wH}(\M(\bw)^{\wH})_{\mathfrak R} =
    K^{\wH}(\Zl^{\wH})_{\mathfrak R}
    = K^{\wH}(\Zl)_{\mathfrak R}
\)
is a simple $\Ul\otimes_{\Q(q)}\mathfrak R(\wH)$-module. Its $k$th
Drinfeld polynomial is given by
\begin{equation*}
   \Wedge_{-u} q^{-1} W_k 
   = \Wedge_{-u} q^{-1} W^1_k \otimes \Wedge_{-u} q^{-1} W^2_k.
\end{equation*}
(Here Drinfeld polynomials have values in $R(\wH)$.)
It is equal to the product of the Drinfeld polynomials of
$K^{\wH}(\La(\bw^1))_{\mathfrak R}$ and
$K^{\wH}(\La(\bw^2))_{\mathfrak R}$. Therefore, it is a subquotient of 
the tensor product module
\begin{equation*}
   K^{\wH}(\La(\bw^1))_{\mathfrak R}\otimes_{\mathfrak R(\wH)}
   K^{\wH}(\La(\bw^2))_{\mathfrak R}.
\end{equation*}
(A well-known argument based on \lemref{lem:comult}(1).)  By the
localization theorem, \eqref{eq:fixed} and the K\"unneth isomorphism
(\thmref{thm:Kunneth}), we know that the dimensions of the both hand
sides are equal. (See also the remark below.) Hence we have the unique
$\Ul\otimes_{\Q(q)}\mathfrak R(\wH)$-isomorphism $\Phi_{\mathfrak R}$
sending $[0]_{\bw^1}\otimes[0]_{\bw^2}$ to $[0]_{\bw}$.
\end{proof}

\begin{Remark}
The left hand side of \eqref{eq:PhiR} is isomorphic to
\(
   K^{\wH}(\La(\bw^1)^{\wH})_{\mathfrak R}\otimes_{\mathfrak R(\wH)}
   K^{\wH}(\La(\bw^2)^{\wH})_{\mathfrak R}
\)
by the localization theorem. Combining it with \eqref{eq:fixed} and
the K\"unneth isomorphism, we have an isomorphism (of $\mathfrak
R(\wH)$-modules) between the left-hand side and the right-hand
side. But it does not respect $\Ul$-module structures.
In fact, a computation in \cite[7.4]{VV} implies that the isomorphism
respects $\Ul$-module structures, if we endow the left hand side with
a $\Ul$-module structure given by Drinfeld's new comultiplication
(after an explicit twist). Thus our map $\Phi_{\mathfrak R}$ should be
given explicitly by factors of the universal $R$-matrix as in
\cite{KT}.
However, it will become difficult (at least for the author) to show
the commutatibity of the diagram \eqref{eq:PhiR} below, in the formula
in \cite{KT}. This is the reason why we do not use the explicit form
of $\Phi_{\mathfrak R}$, unlike \cite{VV}.
\end{Remark}

Let $\La(\bw^1)^{\wH} = \M(\bw^1)^{\wH} = \bigsqcup_{\rho}
\M(\rho;\bw^1)$ and $\La(\bw^2)^{\wH} = \M(\bw^2)^{\wH} =
\bigsqcup_{\rho'} \M(\rho';\bw^2)$ be decomposition into connected
components. Then we have direct sum decomposition
\begin{equation*}
    K^{\wH}(\La(\bw^1))_{\mathfrak R}
    = \bigoplus_\rho K^{\wH}(\M(\rho;\bw^1))_{\mathfrak R},
\quad 
    K^{\wH}(\La(\bw^2))_{\mathfrak R}
    = \bigoplus_{\rho'} K^{\wH}(\M(\rho';\bw^2))_{\mathfrak R}.
\end{equation*}
By \cite[13.4.5]{Na-qaff}, these are {\it l\/}--weight space
decomposition. Here {\it l\/}--weights are elements in $\mathfrak
R(\wH)[[z^\pm]]^I$. Under $\mathfrak R(\wH) = \mathfrak
R(\wH_{W^1})\otimes_{\Q(q)}\mathfrak R(\wH_{W^2})$, they have forms
of $f\otimes 1$ and $1\otimes g$ respectively.
Therefore, by \lemref{lem:comult}(1),
\begin{equation}\label{eq:wt_decomp}
  K^{\wH}(\La(\bw^1))_{\mathfrak R}\otimes_{\mathfrak R(\wH)}
  K^{\wH}(\La(\bw^2))_{\mathfrak R}
 =
  \bigoplus_{\rho,\rho'}
  K^{\wH}(\M(\rho;\bw^1))_{\mathfrak R} \otimes_{\mathfrak R(\wH)}
  K^{\wH}(\M(\rho';\bw^2))_{\mathfrak R}
\end{equation}
is the {\it l\/}--weight space decomposition. The {\it l\/}--weight is 
the form of $f\otimes g$. Therefore each summand has distinct {\it
l\/}--weights.

\begin{Lemma}
The following diagram is commutative:
\begin{equation}\label{eq:PhiR}
\begin{CD}
   K^{\wH}(\La(\bw^2))
   @>{\eqref{eq:Thom2}}>>
   K^{\wH}(\Zl_1)
\\
   @V{\heartsuit}VV @VV\eqref{eq:Zl1_incl}V
\\
   K^{\wH}(\La(\bw^1))\otimes_{R(\wH)} K^{\wH}(\La(\bw^2))
   @.
   K^{\wH}(\Zl)
\\
   @V\otimes_{\mathfrak R(\wH)}VV @VV\otimes_{\mathfrak R(\wH)}V
\\
   K^{\wH}(\La(\bw^1))_{\mathfrak R}\otimes_{\mathfrak R(\wH)}
   K^{\wH}(\La(\bw^2))_{\mathfrak R}
     @>\cong>\Phi_{\mathfrak R}>
   K^{\wH}(\Zl)_{\mathfrak R},
\end{CD}
\end{equation}
where $\heartsuit$ is the inclusion
\begin{equation*}
   K^{\wH}(\La(\bw^2)) \ni E \mapsto [0]_{\bw^1}\otimes E
    \in
   K^{\wH}(\La(\bw^1))\otimes_{R(\wH)} K^{\wH}(\La(\bw^2)).
\end{equation*}
\end{Lemma}

\begin{proof}
Let ${}_{\mathbf m^2}[q]$ be the class represented by the structure
sheaf of $\M(\bw^2-\mathbf m^2,\bw^2) = \La(\bw^2-\mathbf m^2,\bw^2) = 
\text{point}$, as in the previous section (the lowest weight vector).
We have $f_{k,r}\,{}_{\mathbf m^2}[q] = 0$ for any $k\in I$, $r\in\Z$.
Consider the element $T_{w_0}$ of the Braid group corresponding to the
longest element $w_0$ of the Weyl group (of $\mathfrak g$). It acts on
$K^{\wH}(\La(\bw^2))$ by \cite[Part VI]{Lu-book}. Since $w_0 \bw =
\mathbf m^2$, $T_{w_0}$ maps $[0]_{\bw^2}$ to $\alpha\,{}_{\mathbf
m^2}[q]$, where $\alpha$ is an invertible element in $R(\wH)$.
We have 
\[
   \left(\Uli\otimes_{\Z[q,q^{-1}]} R(\wH)\right){}_{\mathbf m^2}[q]
   = \left(\Uli^+\otimes_{\Z[q,q^{-1}]} R(\wH)\right){}_{\mathbf m^2}[q]
   = K^{\wH}(\La(\bw^2)).
\]
Combining this with \propref{prop:qqq2}, we understand that it is
enough to check the commutativity of the diagram for the element
${}_{\mathbf m^2}[q]$. Let ${}'([0]_{\bw^1}\otimes {}_{\mathbf
  m^2}[q])$ be the image of ${}_{\mathbf m^2}[q]$ under the
composition of \eqref{eq:Zl1_incl} and \eqref{eq:Thom2}.
We want to show
\begin{equation}
  \label{eq:want}
   \Phi_{\mathfrak R}([0]_{\bw^1}\otimes{}_{\mathbf m^2}[q])
   = {}'([0]_{\bw^1}\otimes {}_{\mathbf m^2}[q]).
\end{equation}

In the {\it l\/}--weight space decomposition \eqref{eq:wt_decomp},
both hand sides of \eqref{eq:want} is contained in the summand
\begin{equation*}
  K^{\wH}(\M(0,\bw^1))_{\mathfrak R} \otimes_{\mathfrak R(\wH)}
  K^{\wH}(\M(\bw^2-\mathbf m^2,\bw^2))_{\mathfrak R}.
\end{equation*}
Recall both $\M(0,\bw^1)$ and $\M(\bw^2-\mathbf m^2,\bw^2)$ are a
single point. Therefore the {\it l\/}--weight space is
$1$-dimensional. Thus \eqref{eq:want} holds up to a nonzero constant
in ${\mathfrak R}(\wH)$.

By \cite[39.1.2]{Lu-book}, we have
\begin{equation*}
   [0]_{\bw^2} = \alpha T_{w_0}\,{}_{\mathbf m^2}[q]
   = \alpha e_{k_1,0}^{(a_1)} e_{k_2,0}^{(a_2)}
   \cdots e_{k_N,0}^{(a_N)}\,{}_{\mathbf m^2}[q],
\end{equation*}
where $s_{k_1}s_{k_2}\cdots s_{k_N}$ is a reduced expression of $w_0$, 
and $a_i\in\Z_{\ge 0}$ is an explicitly computable natural number.
By \lemref{lem:comult}(2) we have
\begin{equation*}
   [0]_{\bw^1}\otimes [0]_{\bw^2} 
   = \alpha e_{k_1,0}^{(a_1)} e_{k_2,0}^{(a_2)} \cdots e_{k_N,0}^{(a_N)}
   \left([0]_{\bw^1}\otimes{}_{\mathbf m^2}[q]\right)
\end{equation*}
in $K^{\wH}(\La(\bw^1))\otimes_{R(\wH)} K^{\wH}(\La(\bw^2))$.
Therefore,
\begin{equation*}
   \alpha e_{k_1,0}^{(a_1)} e_{k_2,0}^{(a_2)} \cdots e_{k_N,0}^{(a_N)}
   \left(\Phi_{\mathfrak R}([0]_{\bw^1}\otimes{}_{\mathbf m^2}[q])\right)
 = [0]_{\bw}.
\end{equation*}
On the other hand, by \propref{prop:qqq2}(2) we have
\begin{equation*}
   \alpha e_{k_1,0}^{(a_1)} e_{k_2,0}^{(a_2)} \cdots e_{k_N,0}^{(a_N)}
   \; {}'([0]_{\bw^1}\otimes {}_{\mathbf m^2}[q])
   = [0]_{\bw}.
\end{equation*}
We have used the commutativity of the diagram for the element
$[0]_{\bw^2}$, which is obvious from the definition.
Therefore the constant must be $1$.
\end{proof}

\begin{Theorem}\label{thm:main4}
$\Phi_{\mathfrak R}$ induces an
$\Uli\otimes_{\Z[q,q^{-1}]}R(\wH)$-module isomorphism
\begin{equation*}
   \Phi\colon
   K^{\wH}(\La(\bw^1))\otimes_{R(\wH)}
   K^{\wH}(\La(\bw^2))
\xrightarrow{\cong}
   K^{\wH}(\Zl).
\end{equation*}
\end{Theorem}

\begin{proof}
Since $K^{\wH}(\La(\bw^1))$ is generated by $[0]_{\bw^1}$, 
$K^{\wH}(\La(\bw^1))\otimes_{R(\wH)}K^{\wH}(\La(\bw^2))$ is generated
by $[0]_{\bw^1}\otimes K^{\wH}(\La(\bw^2))$.
Therefore, $K^{\wH}(\La(\bw^1))\otimes_{R(\wH)} K^{\wH}(\La(\bw^2))$
is mapped to $\Uli K^{\wH}(\Zl_1)$ under $\Phi_{\mathfrak R}$.
But it is equal to $K^{\wH}(\Zl)$ by \propref{prop:qqq2}(1).
\end{proof}

As an application, we have a new proof of \cite[7.12]{VV}.
\begin{Corollary}\label{cor:VV}
Let $\varepsilon\in\C^*$.
Suppose two $I$-tuple polynomials $P^1 = (P^1_k)_{k\in I}$, $P^2 =
(P^2_k)_{k\in I}$ satisfy that
\begin{equation*}
   \alpha/\alpha' \notin \{ \varepsilon^n \mid n\in \Z, n \ge 2 \}
   \quad
   \text{for any pair $(\alpha,\alpha')$ with $P^1_k(\alpha) = 0$,
   $P^2_{k'}(\alpha') = 0$ \rom($k,k'\in I$\rom)}.
\end{equation*}
Then we have a unique $\Ule$-isomorphism
\begin{equation*}
   M({P^1P^2}) \cong M(P^1)\otimes_\C M(P^2),
\end{equation*}
sending $[0]_{\bw}$ to $[0]_{\bw^1}\otimes [0]_{\bw^2}$.
\end{Corollary}

\begin{proof}
Take diagonal matices $s^1\in H_{W^1}$, $s^2\in H_{W^2}$ whose
entries are roots of $P^1$, $P^2$ respectively.
The evaluation at $a = (s^1,s^2,\varepsilon)\in H_{W^1}\times
H_{W^2}\times\C^* = \wH$ defines a homomorphism
\(
   R(\wH) = R(H)[q,q^{-1}] \to \C
\).
Then we set
\begin{equation*}
\begin{gathered}
   M(P^1P^2) = K^{\wH}(\La(\bw)) \otimes_{R(\wH)} \C,
\\
   M(P^1)\otimes_\C M(P^2)
   = K^{\wH}(\La(\bw^1))\otimes_{R(\wH)}
     K^{\wH}(\La(\bw^2))\otimes_{R(\wH)}\C
   = K^{\wH}(\Zl)\otimes_{R(\wH)}\C.
\end{gathered}
\end{equation*}
The inclusion in \thmref{thm:free}(3) induces a $\Ule$-homomorphism
\(
   M(P^1P^2) \to M(P^1)\otimes_\C M(P^2)
\),
sending $[0]_{\bw}$ to $[0]_{\bw^1}\otimes [0]_{\bw^2}$.
If we denote by $X_{a}$ the localization of an $R(\wH)$-module $X$ at
the maximal ideal corresponding to $a$, the above equalities factor
through as
\begin{equation*}
   M(P^1P^2) = K^{\wH}(\La(\bw))_a \otimes_{R(\wH)_a} \C,
   \quad
   M(P^1)\otimes M(P^2)
   = K^{\wH}(\Zl)_a\otimes_{R(\wH)_a}\C.
\end{equation*}
By the localization theorem \cite{T-loc}, we have
\begin{equation*}
   K^{\wH}(\La(\bw))_a \cong K^{\wH}(\La(\bw)^a)_a,\quad
   K^{\wH}(\Zl)_a \cong K^{\wH}(\Zl^a)_a
\end{equation*}
where $\ast^a$ denotes the set of points fixed by $a$. So the result
follows from the following.
\begin{Claim}
$\La(\bw)^a = \Zl^a$. 
\end{Claim}
Let $W(\alpha) = W^1(\alpha)\oplus W^2(\alpha)$ be the eigenspace of
$s_1\oplus s_2$ with eigenvalue $\alpha$.
Let $[B,i,j]\in\Zl^a$ and 
\[
   f\defeq j_{\vin(h_N)} B_{h_N} B_{h_{N-1}}\cdots B_{h_1}
   i_{\vout(h_1)}\colon W_{\vout(h_1)} \to W_{\vin(h_N)},
\]
where $h_1$, \dots, $h_N$ is a path in our graph. Since $[B,i,j]$ is
fixed by $a$, $f$ maps $W_{\vout(h_1)}(\alpha)$ to
$W_{\vin(h_N)}(\varepsilon^{-2-N}\alpha)$. On the other hand, the
condition $[B,i,j]\in\Zl$ implies that $f$ maps $W^2$ to $0$ and $W^1$
to $W^2$ as in the proof of \lemref{lem:closed}.
Therefore, $f$ must be $0$ by the assumption. Since the path is
arbitrary, it means that $\pi([B,i,j]) = 0$, i.e.,
$[B,i,j]\in\La(\bw)$.
\end{proof}

\begin{Remark}
In our proof of \thmref{thm:main4}, the assumption that $\mathfrak g$
is of type ADE is used for the existence of the comultiplication
$\Delta$. If one can prove the existence of it such that
\lemref{lem:comult} still holds, then our proof goes well for a
general Kac-Moody Lie algebra $\mathfrak g$. Or, if one can give a
different proof of the existence of the isomorphism $\Phi$ in
\thmref{thm:main4} as merely $R(\wH)$-modules, without using
$\Uli$-module structures, then it means that one can consider $\Phi$
as a `definition' of the tensor product module.
\end{Remark}

\section{General case}\label{sec:general}

Almost all results in previous sections can be generalized to the case
of a tensor product of more than two modules.

Let us suppose a direct sum decomposition $W=W^1\oplus W^2\oplus
\cdots \oplus W^N$ of $I$-graded vector spaces is given.
Let $H_{W^p}\subset G_{W^p}$ be the maximal torus of diagonal
matrices, and let $H \defeq H_{W^1}\times H_{W^2}\times \cdots \times
H_{W^N}$, and let $\wH \defeq H\times\C^*$.
We choose a one-parameter subgroup $\lambda\colon \C^*\to
G_{W^1}\times G_{W^2}\times\cdots\times G_{W^N}$ given by
\begin{equation*}
   \lambda(t) = t^{m_1}\id_{W^1}\oplus t^{m_2} \id_{W^2} \oplus\cdots 
   \oplus t^{m_N} \id_{W^N},
\end{equation*}
with $m_1 < m_2 < \cdots < m_N$. Moreover we take generic $m_i$'s and 
assume that the fixed point set $\M(\bw)^{\lambda(\C^*)}$ is
$\M(\bw^1)\times \M(\bw^2)\times\cdots\times \M(\bw^N)$.
We define
\begin{gather*}
   \Zm(\bw^1;\bw^2;\cdots;\bw^N)
   \defeq \left\{ [B,i,j]\in \M(\bw) \left| \;
   \text{$\lim_{t\to 0} \lambda(t)\ast[B,i,j]$ exists} \right\}\right.,
\\
   \Zl(\bw^1;\bw^2;\cdots;\bw^N)
   \defeq \left\{ [B,i,j]\in \M(\bw)  \left| \;
   \lim_{t\to 0} \lambda(t)\ast[B,i,j]\in
   \La(\bw^1)\times\cdots\times\La(\bw^N)
   \right\}\right..
\end{gather*}
Equivalently, we can define inductively
\begin{gather*}
\begin{split}
   & \Zm(\bw^1;\bw^2;\cdots;\bw^N)
\\
   & \qquad = \left\{ [B,i,j]\in \M(\bw) \left| \;
   \lim_{t\to 0} \lambda'(t)\ast[B,i,j]\in
   \Zm(\bw^1;\bw^2;\cdots;\bw^{N-1})\times \M(\bw^N)
   \right\}\right.,
\end{split}
\\
\begin{split}
   & \Zl(\bw^1;\bw^2;\cdots;\bw^N)
\\
   & \qquad
   = \left\{ [B,i,j]\in \M(\bw)  \left| \;
   \lim_{t\to 0} \lambda'(t)\ast[B,i,j]\in
   \Zl(\bw^1;\bw^2;\cdots;\bw^{N-1})\times\La(\bw^N)
   \right\}\right.,
\end{split}
\end{gather*}
where
\begin{equation*}
   \lambda'(t) = \id_{W^1}\oplus\cdots \oplus\id_{W^{N-1}}
   \oplus t\id_{W^N}.
\end{equation*}
These are closed subvarieties and $\Zl(\bw^1;\bw^2;\cdots;\bw^N)$ is
lagrangian. We have
\begin{enumerate}
\item $\Irr\Zl(\bw^1;\bw^2;\cdots;\bw^N)$ has a structure of a crystal 
isomorphic to
$\B(\bw^1)\otimes\cdots\otimes\B(\bw^N)$.
(\thmref{thm:main1})

\item $H_{\topdeg}(\Zl(\bw^1;\bw^2;\cdots;\bw^N),\Q)$ is isomorphic to 
$V(\bw^1)\otimes\cdots\otimes V(\bw^N)$ as a
$\mathfrak g$-module. (\thmref{thm:main2})

\item (When $\mathfrak g$ is of type $ADE$)
$K^{\wH}(\Zl(\bw^1;\bw^2;\cdots;\bw^N))$ is isomorphic to
$K^{\wH}(\La(\bw^1))\otimes_{R(\wH)}\cdots\otimes_{R(\wH)}
K^{\wH}(\La(\bw^N))$ as a
$\Uli\otimes_{\Z[q,q^{-1}]}R(\wH)$-module. (\thmref{thm:main4})
\end{enumerate}

\thmref{thm:main3} depends on \cite[23.3.6]{Lu-book}, which seems
difficult to generalize. This is the only reason why the author does
not know the generalization of \thmref{thm:main3}.

\section{Combinatorial description of the crystal}\label{sec:comb}

In this section, we give a combinatorial description of the cristal
$\Irr\Zl(\bw^1;\dots;\bw^N)$.

Since $\Irr\Zl(\bw^1;\dots;\bw^N) =
\Irr\La(\bw^1)\otimes\dots\otimes\Irr\La(\bw^N)$, it is enough to give 
a description of $\Irr\La(\bw^p)$.
However, we study a slightly general situation.

For a given one-parameter subgroup $\rho_0\colon\C^*\to G_W$, we
define a $\C^*$-action on $\bM$ and $\M$ by
\begin{equation*}
    B_h \mapsto
    \begin{cases}
      B_h & \text{if $h\in\Omega$},
      \\
      t B_h & \text{if $h\in\overline\Omega$},
    \end{cases}
\quad    
    i \mapsto \rho_0(t)\ast i,
\quad
    j\mapsto  \rho_0(t)\ast (tj).
\end{equation*}
We denote this $\C^*$-action by
$(B,i,j)\mapsto t\diamond (B,i,j)$ and $[B,i,j]\mapsto t\diamond
[B,i,j]$.
If $\rho_0(t)=\id_W$, this is the $\C^*$-action considered in
\cite[\S5]{Na-quiver}. Its crucial property is that the symplectic
form $\omega$ is transformed as $t\omega$. We define
\begin{equation*}
   \Zl^\diamond\defeq
   \left\{ x\in \M(\bw) \left|\; \text{$\lim_{t\to \infty} t\diamond x$
   exists}\right\}\right..
\end{equation*}
%($\Zl^\diamond$ depends on $\Omega$ and $\rho_0$, as well as $\bw$.)
Note that we consider the limit for $t\to\infty$ while we have studied 
the limit ${t\to 0}$ in previous sections.
When $\Omega$ contains no cycle and $\rho_0(t) = \id_W$, we have
$\Zl^\diamond = \La(\bw)$ \cite[5.3(2)]{Na-quiver}.

It is easy to show that the set $\Irr\Zl^\diamond$ of irreducible
components of $\Zl^\diamond$ has a structure of a normal crystal as in
\secref{sec:crystal}.

Take a fixed point $x$ of the $\C^*$-action and its representative
$(B,i,j)$. Then there exists a unique homomorphism $\rho\colon\C^*\to
G_V$ such that
\begin{equation}\label{eq:equation}
   t \diamond (B,i,j) = \rho(t)^{-1}\cdot (B,i,j).
\end{equation}
(The uniqueness follows from the freeness of the action of $G_V$ on
the set of stable points.) Moreover, the conjugacy class of $\rho$ is
independent of the choice of the representative $(B,i,j)$ of $x$.
Thus we have a decomposition
\begin{equation}\label{eq:fix}
   \M(\bw)^{\C^*} = \bigsqcup_\rho \mathfrak F^\diamond(\rho),
\end{equation}
where $\mathfrak F^\diamond(\rho)$ is the set of $[B,i,j]$ satisfying
\eqref{eq:equation} up to conjugacy.
It is clear that the conjugacy class of $\rho$ is constant on each
connected component.
Furthremore, by the argument in \cite[5.5.6]{Na-qaff}, we can show
that each summand $\mathfrak F^\diamond(\rho)$ is the connected
component of $\M(\bw)^{\C^*}$ (if it is nonempty). (The $\C^*$-action
used in [loc.\ cit.] is different from the above one. But the argument
still works. Moreover the assumption $-(\alpha_k,\alpha_l)\le 1$ there
becomes unnecessary for the above $\C^*$-action.)

Combining with the argument in \cite[5.8]{Na-quiver}, we get the
following.
\begin{Proposition}
We have a decomposition
\begin{equation*}
   \Zl^\diamond
   = \bigsqcup_\rho \Zl^\diamond(\rho); \qquad
  \Zl^\diamond(\rho)\defeq
   \left\{ x\in \M(\bw) \left|\; \lim_{t\to \infty} t\diamond x\in
       \mathfrak F^\diamond(\rho)\right\}\right..
\end{equation*}
The irreducible components of $\Zl^\diamond$ is the closure of each
summand $\Zl^\diamond(\rho)$, and they are all lagrangian
subvarieties.
\end{Proposition}

Here the property $t^*\omega = t\omega$ played the crusial role.

The conjugacy class of $\rho$ corresponds bijectively to the
dimensions of its weight space in $V_k$. Thus we have an injective map
\begin{equation*}
   \Irr\Zl^\diamond
   \to \Z^{I\times\Z};
\qquad
   \Zl^\diamond(\rho) \mapsto \left(\dim V_k^p\right)_{k\in I,p\in\Z},
\end{equation*}
where
\(
   V^p = \{ v\in V \mid \rho(t)\cdot v = t^p v \}
\).

We set
\(
   W^p = \{ w\in W \mid \rho_0(t)\ast w = t^p w \}
\).
Then \eqref{eq:equation} is equivalent to
\begin{gather*}
   B_h(V_{\vout(h)}^p)
   \subset
   \begin{cases}
      V_{\vin(h)}^p & \text{if $h\in\Omega$},
      \\
      V_{\vin(h)}^{p-1} & \text{if $h\in\overline\Omega$},
   \end{cases}
\qquad
%\\
    i_k(W_k^p)\subset V_k^p, \quad
    j_k(V_k^p)\subset W_k^{p-1}.
\end{gather*}

By the same formula in \eqref{eq:hecke_complex}, we define an
analogous complex of vector bundles over $\mathfrak F^\diamond(\rho)$
for each $p,q\in\Z$:
\begin{equation*}
  \HomL(V^p, V^q)
  \overset{\alpha^{qp}}{\longrightarrow}
  \begin{matrix}
    \HomE_\Omega(V^p,V^q)
    \\
    \oplus
    \\
    \HomE_{\overline\Omega}(V^p,V^{q-1})
  \end{matrix}
  \oplus
  \begin{matrix}
  \HomL(W^p, V^q)
  \\
  \oplus
  \\
  \HomL(V^p,W^{q-1})
  \end{matrix}
  \overset{\beta^{qp}}{\longrightarrow}
  \HomL(V^p, V^{q-1}),
\end{equation*}
where
\begin{equation*}
  \HomE_{\Omega}(V^p,V^q) =
  \bigoplus_{h\in\Omega} \Hom(V^p_{\vout(h)}, V^q_{\vin(h)}),
\qquad
  \HomE_{\overline\Omega}(V^p,V^{q-1}) =
  \bigoplus_{h\in\overline\Omega} \Hom(V^p_{\vout(h)}, V^{q-1}_{\vin(h)})
\end{equation*}
The restriction of the tangent bundle $T\M$ to $\mathfrak F^\diamond(\rho)$
decomposes as
\(
   \bigoplus_{p,q}  \Ker\beta^{qp}/\Ima\alpha^{qp}
\).
Since the tangent space to the fixed point set is the $0$-weight
space,
\begin{equation*}
   T\mathcal F^\diamond(\rho) \cong
   \bigoplus_p \Ker\beta^{pp}/\Ima\alpha^{pp}.
\end{equation*}
Moreover, as in \propref{prop:ext} we have
\begin{equation*}
   \Zl^\diamond(\rho) \cong \bigoplus_{q < p} \Ker\beta^{qp}/\Ima\alpha^{qp},
\end{equation*}
where the right hand side is the total space of the vector bundle.
The natural map $\Zl^\diamond(\rho)\to \mathfrak F^\diamond(\rho)$ is
identified with the projection map of the vector bundle. 

Take a point $[B,i,j]\in\Zl^\diamond(\rho)$ and consider the complex
$C_k^\bullet$ \eqref{eq:taut_cpx}. The homomorphisms $\sigma_k$,
$\tau_k$ have the matrix expression
$\sigma_k = (\sigma_k^{qp})_{p,q}$, $\tau_k = (\tau_k^{qp})_{p,q}$, where
\begin{gather*}
   \sigma_k^{qp}\colon V_k^{p} \longrightarrow
  \bigoplus_{h\in\Omega: \vin(h) = k} V_{\vout(h)}^{q-1}
  \oplus
  \bigoplus_{h\in\overline\Omega: \vin(h) = k} V_{\vout(h)}^{q}
    \oplus W_k^{q-1},
\\
  \tau^{qp}_k\colon
  \bigoplus_{h\in\Omega: \vin(h) = k} V_{\vout(h)}^{p-1}
  \oplus
  \bigoplus_{h\in\overline\Omega: \vin(h) = k} V_{\vout(h)}^{p}
    \oplus W_k^{p-1}
    \longrightarrow
  V_k^{q-1}.
\end{gather*}
Components $\sigma_k^{qp}$, $\tau_k^{qp}$ vanish if $q > p$. From the 
equation $\tau_k\sigma_k = 0$, we have
\(
   \sum_{r:p \le r \le q} \tau_k^{qr}\sigma^{rp}_k = 0
\).
In particular, we have the complex
\begin{equation*}
   C_k^{p\bullet}\colon
   V_k^{p}
\xrightarrow{\sigma_k^{pp}}
   \bigoplus_{h\in\Omega: \vin(h) = k} V_{\vout(h)}^{p-1}
   \oplus
   \bigoplus_{h\in\overline\Omega: \vin(h) = k} V_{\vout(h)}^{p}
    \oplus W_k^{p-1}
\xrightarrow{\tau_k^{pp}}
   V_k^{p-1}.
\end{equation*}

\begin{Lemma}
Let
\begin{equation*}
   \overline\varepsilon_k^p \defeq - \sum_{q:q > p} \rank C_k^{q\bullet},
\qquad
   \overline\varphi_k^p \defeq \sum_{q:q\le p} \rank C_k^{q\bullet}.
\end{equation*}
We have
%\begin{gather*}
\begin{equation*}
   \varepsilon_k(\Zl^\diamond(\rho)) = 
   \max_{p\in\Z} \overline\varepsilon_k^p,
\quad
   \varphi_k(\Zl^\diamond(\rho)) = 
   \max_{p\in\Z} \overline\varphi_k^p.
\end{equation*}
Let $\rho'$ \rom(resp.\ $\rho''$\rom) is the one parameter subgroup,
obtained from $\rho$, with $\dim V_k^{p}$ decreased \rom(resp.\ 
increased\rom) by $1$, and other components unchanged, where
\begin{gather*}
   p = \min\left\{ q \left| \overline\varepsilon_k^q
     = \varepsilon_k(\Zl^\diamond(\rho))\right\}\right.
\\
   \left(\text{resp.\ }
   p = \max\left\{ q \left| \overline\varphi_k^q
     = \varphi_k(\Zl^\diamond(\rho))\right\}\right.
   \right).
\end{gather*}
Then
\begin{equation*}
   \widetilde e_k(\Zl^\diamond(\rho)) =
   \begin{cases}
     0 & \text{if $\varepsilon_k(\Zl^\diamond(\rho)) = 0$},\\
     \Zl^\diamond(\rho') & \text{otherwise},   
   \end{cases}
\quad
   \widetilde f_k(\Zl^\diamond(\rho)) =
    \begin{cases}
     0 & \text{if $\varphi_k(\Zl^\diamond(\rho)) = 0$},\\
     \Zl^\diamond(\rho'') & \text{otherwise}.
   \end{cases}
\end{equation*}
\end{Lemma}

\begin{proof}
The situation is almost the same as that studied in
\secref{sec:crystal}. So $\varepsilon_k(\Zl^\diamond(\rho))$ is given
by the same formula as in the tensor product crystal, if we know the
codimension of $\Ima\tau^{pp}_k$ in $V_k^{p-1}$.
In fact, it was given in \cite[5.5.5]{Na-qaff}. We have
\begin{equation*}
%\begin{split}
%   & \Ker\sigma^{pp}_k = 0
%   \quad
%      \text{on every point $[B,i,j]\in\mathfrak F^\diamond(\rho)$},
%\\
%   & 
   \dim V_k^{p-1}/\Ima\tau^{pp}_k = \max(0,-\rank C_k^{p\bullet})
   \quad
   \text{on a generic point $[B,i,j]\in\mathfrak F^\diamond(\rho)$}.
%\end{split}
\end{equation*}
(Although the $\C^*$-action in \cite{Na-qaff} is different from our
$\C^*$-action, the argument works.)

Now we repeat the argument in \secref{sec:crystal}, where $\wt_k(b_p)$
(resp.\ $\varepsilon_k(b_p)$) is replaced by $\rank C_k^{p\bullet}$
(resp.\ $\max(0, -\rank C_k^{p\bullet})$), and the order of the tensor
product is reversed.
Therefore, we get
{\allowdisplaybreaks
\begin{equation*}
\begin{split}
  \varepsilon_k(\Zl^\diamond(\rho))
   & = \max_{p\in\Z}
   \left(
     \max(0,-\rank C_k^{p\bullet}) - \sum_{q:q>p} \rank C_k^{q\bullet}
   \right)
\\
   & = \max_{p\in\Z}
   \left(
     - \sum_{q:q\ge p} \rank C_k^{q\bullet}
   \right)
   = \max_{p\in\Z} \overline\varepsilon_k^p,
\end{split}
\end{equation*}
and other formulas}.
\end{proof}

Let $\bw^p = \sum_k \dim W_k^p\, \Lambda_k$, $\bv^p = \sum_k \dim
V_k^p\, \alpha_k$. We define a crystal $\widetilde T_p$ by
\begin{equation*}
  \widetilde C_p \defeq
    T_{\bw^p}\otimes \bigotimes_{k\in I} \B_k,
\end{equation*}
where $\B_k$ is the crystal in Example~\ref{exp:crystal}(1), and we
have used the numbering of $I$ to determine the order of the tensor
product.
Let $S_0$ be the crystal consisting of a single element $s_0$ with
$\wt s_0 = 0$, $\varepsilon_k(s_0) = \varphi_k(s_0) = 0$, $\widetilde
e_k s_0 = \widetilde f_k s_0 = 0$.
We define
\begin{equation*}
   \widetilde C \defeq 
   \cdots\otimes\widetilde C_{p+1}
   \otimes \widetilde C_p\otimes \widetilde C_{p-1} \otimes \cdots.
\end{equation*}

\begin{Theorem}
We have a strict embedding of crystal
\begin{alignat*}{2}
   \psi\colon & \Irr\Zl^\diamond &\to &S_0\otimes\widetilde C\otimes S_0,
\\
   & \Zl^\diamond(\rho) & \mapsto
   & s_0\otimes\left(\cdots\otimes
   \left(t_{\bw^p}\otimes \otimes_k b_k(-\dim V_k^p) \right)
   \otimes\cdots\right)\otimes s_0.
\end{alignat*}
\end{Theorem}

\begin{proof}
It is clear that $\psi$ commutes with $\wt$.
%\begin{equation*}
%   \wt(\Zl^\diamond(\rho)) = \bw - \bv
%   = \sum_p \bw^p - \bv^p
%   = \sum_p \wt\left(t_{\bw^p}\otimes \otimes_k b_k(-\dim V_k^p) \right).
%\end{equation*}

Let $\psi'\colon \Irr\Zl^\diamond\to \widetilde C$ be the map defined
by omitting the first and last $s_0$ from the above formula.
Let
\(
   b = \cdots\otimes b_{p+1}\otimes b_p\otimes b_{p-1}\otimes\cdots
   \in \widetilde C.
\)
We define $\varepsilon_k^p$, $\varphi^p_k$ as in \eqref{eq:vep}.
If $b = \psi'(\Zl^\diamond(\rho))$, then
\begin{equation*}
\begin{split}
   \varepsilon_k^p
   = & \dim V_k^p - \sum_{q\ge p} \dim W_k^q
   + \sum_{\substack{l\in I:l < k\\ q\ge p}}
   \langle h_k,\alpha_l\rangle \dim V_l^q
\\
   & \qquad\qquad\qquad
   + \sum_{\substack{l\in I:l > k\\ q > p}}
   \langle h_k,\alpha_l\rangle \dim V_l^q
   + 2 \sum_{q > p} \dim V_k^q
\\
  = & - \sum_{q > p} \rank C_k^{q\bullet}
  = \overline\varepsilon_k^p,
\end{split}
\end{equation*}
where we have used that the property $h\in\Omega\Rightarrow \vout(h) <
\vin(h)$, $h\in\overline\Omega\Rightarrow\vout(h) > \vin(h)$.

Thus $\psi'$ commutes with $\varepsilon_k$. From the formula
$\varphi_k^p = \varepsilon_k^p + \wt_k(\cdots\otimes
b_p\otimes\cdots)$, it also commutes with $\varphi_k$.
Now we get (\ref{eq:mor}b, c) for $\psi'$ by \eqref{eq:sub} and the
preceding lemma.

Then we get
\begin{gather*}
   \varepsilon_k(\psi(\Zl^\diamond(\rho))) =
%   \varepsilon_k(s_0\otimes\psi'(\Zl^\diamond(\rho))\otimes s_0) =
   \max(0,\varepsilon_k(\Zl^\diamond(\rho)), -\wt_k(\Zl^\diamond(\rho)))
   = \varepsilon_k(\Zl^\diamond(\rho)),
\\
   \varphi_k(\psi(\Zl^\diamond(\rho))) =
%   \varphi_k\left(s_0\otimes\psi'(\Zl^\diamond(\rho))\otimes s_0\right) =
   \max(0,\varphi_k(\Zl^\diamond(\rho)),\wt_k(\Zl^\diamond(\rho)))
   = \varphi_k(\Zl^\diamond(\rho)),
\end{gather*}
where we have used $\varepsilon_k(\Zl^\diamond(\rho))$,
$\varphi_k(\Zl^\diamond(\rho))\ge 0$, which is
clear from the definition.
We have
\begin{equation*}
   \widetilde e_k(\psi(\Zl^\diamond(\rho)))
   = 
   \begin{cases}
     s_0\otimes \psi'\widetilde e_k(\Zl^\diamond(\rho))\otimes s_0 
      & \text{if $\varepsilon_k(\Zl^\diamond(\rho)) \neq 0$},
      \\
     0 & \text{if $\varepsilon_k(\Zl^\diamond(\rho)) = 0$}.
   \end{cases}
\end{equation*}
Since $\varepsilon_k(\Zl^\diamond(\rho)) = 0\Leftrightarrow
\widetilde e_k(\Zl^\diamond(\rho)) = 0$, $\widetilde e_k$ commutes
with $\psi$. Similarly, $\widetilde f_k$ commutes with $\psi$.
\end{proof}

Let us take $\rho_0(t) = \id_W$ and hence $\Zl^\diamond =
 \La(\bw)$. 
Then the strictly embedded crystal generated by
\(
   s_0\otimes
   \cdots\otimes(t_{\bw^p}\otimes \otimes_k b_k(0))\otimes\cdots
   \otimes s_0
\)
in the right hand side is the same as a combinatorial description of
the crystal $\B(\bw)$ in \cite{Kas-Dem}. Thus we obtain a different
proof of \corref{cor:cryst}.

\section{Examples}\label{sec:example}

In this section, we give examples of $\Zl$.

\subsection{}
Suppose the graph is of type $A_n$. We number the vertices as
$$
\def\longlongrightarrow{\relbar\joinrel\relbar\joinrel}
\catcode`\@=11
\newbox\tempbox \newbox\bulletbox
\def\numberedbullet#1{\setbox\bulletbox=\hbox{$\bullet$}
  \setbox\tempbox=\hbox{$\scriptstyle #1$}
  \kern .5\wd\bulletbox \kern -.5\wd\tempbox
  \raise1.5ex\copy\tempbox \kern -.5\wd\tempbox \kern -.5\wd\bulletbox
  \mathord\bullet}
  \m@th\numberedbullet{1}\longlongrightarrow
  \numberedbullet{2}\longlongrightarrow
  \numberedbullet{3}\longlongrightarrow
  \cdots\longlongrightarrow
  \numberedbullet{n-2}\longlongrightarrow
  \numberedbullet{n-1}\longlongrightarrow
  \numberedbullet{n}
\catcode`\@=\active
$$
We take
$\bw = r\Lambda_1$, $\bv = \sum_{k=1}^n v_k \alpha_k$ with
$r \ge v_1 \ge v_2 \ge \dots \ge v_n \ge 0$. By \cite[7.3]{Na-quiver}
$\M(\bv,\bw)$ is isomorphic to the cotangent bundle $T^*\mathcal F$ of
the partial flag variety $\mathcal F$ consisting of all sequences
\(
  \phi = 
  (\C^r = V_0 \supset V_1 \supset \cdots \supset V_n \supset V_{n+1} = 0) 
\)
with $\dim V_k = v_k$. The correspondence is given by
\begin{equation*}
   [B,i,j] \longmapsto
   (\phi,\xi) = 
   \left((W_1 \supset \Ima j_1 \supset \Ima(j_1 B_{1,2}) \supset \cdots
   \supset \Ima(j_1 B_{1,2}\cdots B_{n-1,n}) \supset 0),
   j_1i_1\right),
\end{equation*}
where $\xi$ is a cotangent vector, i.e., an endomorphism of $\C^r$
with $\xi(V_k)\subset V_{k+1}$ ($k = 0,\dots, n$).

Let $W = W^1\oplus \cdots \oplus W^N$ be a direct sum
decomposition. We have the associated flag
\begin{equation*}
   W = \breve{W}^0 \supset \breve{W}^1 = \bigoplus_{p>1} W^p
   \supset \breve{W}^2 = \bigoplus_{p>2} W^p
   \supset \cdots \supset \breve{W}^{N-1} = W^N \supset \breve{W}^N =
   \{ 0\}.
\end{equation*}
Then we have
\begin{equation*}
\begin{split}
   \Zl & = \{(\phi,\xi)\in T^*\mathcal F\mid
     \xi(\breve{W}^p)\subset \breve{W}^{p+1} \}
\\
   \Zl(\bv^1,\bw^1;\cdots;\bv^N,\bw^N) &= 
    \left\{ (\phi,\xi)\in\Zl \left|
    \dim\left(V_k\cap\breve{W}^p\right)
    = \sum_{q:q>p} v_k^q \right\}\right.,
\end{split}
\end{equation*}
where $v_k^q = \langle h_k,\bv^q\rangle$.
Therefore each stratum of $\Zl$ is the conormal bundle of a Schubert cell.

\subsection{}
Again suppose the graph is of type $A_n$. We take
$\bw = \Lambda_1 + \Lambda_n$, $\bv = \sum_{k=1}^n \alpha_k$. By a
work of Kronheimer, $\M(\bv,\bw)$ is the minimal resolution of the
simple singularity $\C^2/\Z_{n+1}$ (see
\cite[Chapter~4]{Lecture}). The lagrangian subvariety $\La(\bv,\bw)$
is the exceptional set. It is a union of $n$ projective lines.
The intersectection graph is of type $A_n$.

Take coordinates $(x,y)$ of $\C^2$ and suppose the action of
$\Z_{n+1}$ is given by $(x,y)\mapsto (\zeta x,\zeta^{-1} y)$ where
$\zeta$ is a primitive $(n+1)$th root of unity. We have a $\C^*$-action
given by $(x,y)\mapsto (t^{-1}x,ty)$ commuting with the
$\Z_{n+1}$-action. This action lifts to an action on $\M(\bv,\bw)$ and 
coincides with the action considered in \secref{sec:vardef} with
$\bw^1 = \Lambda_1$, $\bw^2 = \Lambda_n$
(after composed with $t\mapsto t^{n+1}$). Then we have
\begin{equation*}
   \Zm = \Zl = (\text{the exceptional set})\cup
   (\text{the strict transform of the $y$-axis}).
\end{equation*}
Each stratum $\Zl(\bv^1,\bw^1;\bv^2,\bw^2)$ is isomorphic to the
affine line $\C$. The intersectection graph of closures of stratum is
of type $A_{n+1}$, where the strict transform of $y$-axis is the last vertex.

\end{document}